\documentclass{article}

\usepackage{phd_package}

\usepackage{phd_command}

\author{Sixtine Michel\footremember{cea}{CEA CESTA, 15 av. des Sabli\`eres, 33114 Le Barp, France.}, Davide Torlo\footremember{Trieste}{Mathematics Area, SISSA Mathlab, SISSA, via Bonomea, 265, 34136 Trieste, Italy.}, Mario Ricchiuto\footremember{inria}{Team CARDAMOM, Inria Bordeaux sud-Ouest, 200 av.  de la vieille tour, 33405 Talence, France.}\,, R\'emi Abgrall\footremember{zurich}{Institut f\"ur Mathematik, Winterthurstrasse 190, CH 8057 Z\"urich, Switzerland.}    }

\title{Spectral analysis of high order continuous FEM for hyperbolic PDEs on triangular meshes: influence of approximation, stabilization, and time-stepping}
\date{\today}

\begin{document}
\graphicspath{{figures/}}

\maketitle

\begin{abstract}
In this work we study various continuous finite element discretization for two dimensional hyperbolic partial differential equations, varying the polynomial space (Lagrangian on equispaced, Lagrangian on quadrature points (\textit{Cubature}) and Bernstein), the stabilization techniques (streamline-upwind Petrov-Galerkin, continuous interior penalty, orthogonal subscale stabilization) and the time discretization (Runge-Kutta (RK), strong stability preserving RK and deferred correction). This is an extension of the one dimensional study by Michel S. et al \textit{J. Sci. Comput.} (2021), whose results do not hold in multi-dimensional frameworks.
The study ranks these schemes based on efficiency (most of them are mass-matrix free), stability and dispersion error, providing the best CFL and stabilization coefficients. 
The challenges in two-dimensions are related to the Fourier analysis. Here, we perform it on two types of periodic triangular meshes varying the angle of the advection, and we combine all the results for a general stability analysis.
Furthermore, we introduce additional high order viscosity to stabilize the discontinuities, in order to show how to use these methods for tests of practical interest.
All the theoretical results are thoroughly validated numerically both on linear and non-linear problems, and  error-CPU time curves are provided. Our final conclusions suggest that \textit{Cubature} elements combined with SSPRK and OSS stabilization is the most promising combination.

\end{abstract}


%
\section{Introduction}\label{sec:intro}
We study several continuous finite element formulations to approximate the solution of the two dimensional hyperbolic conservation laws
\begin{equation}\label{eq:conservation_law2D}
	\partial_t u (x,t) + \nabla \cdot f(u(x,t)) = 0 \quad x\in \Omega \subset \R, \, t\in \R^+,
\end{equation}
where $\Omega\subset \R^2$ is the domain, $f:\R^D\to\R^{2\times D}$ is the flux function and $u:\Omega\to \R^D$ is the unknown of the system of equations. 

The  largest part of the paper is dedicated to the two-dimensional spectral analysis of different stabilized approaches applied to the scalar ($D=1$) transport equations 
 obtained for
\begin{equation}\label{eq:conservation_law1}
 f(u(x,t)) = \mathbf{a} u(x,t)\, \qquad \mathbf{a} \in \R^2 .
\end{equation}
One of the main objectives of this paper is to identify strategies to build (linearly) stable fully explicit high order continuous finite element schemes
to discretize  \eqref{eq:conservation_law2D} on triangulations of the spatial domain $\Omega$. To this end we
will vary the basis functions, the stabilization technique and the time discretization.
In general, the standard Finite Element Method (FEM) derived by this approach require the inversion of a large sparse mass matrix.
This procedure can be expensive as either inverting the mass matrix and, hence, the matrix multiplications must be repeated for every time step or the linear solver must be applied at each time step.
Various techniques have been introduced to overcome the mass matrix inversion while keeping the high order accuracy of the scheme. 

The first strategy we study is the one proposed in \cite{DeC_2017}. In the reference it is suggested to combine mass lumping with a 
deferred correction (DeC) iterative time integration method  allowing to introduce appropriate  corrections in the
right--hand side in order to recover the original order of accuracy. This approach can only be used in combination with finite elements whose basis functions have  positive
integrals.
Another approach is based on a 
careful choice of approximation points defining sufficiently accurate quadrature formulas with all positive weights. If the variational form is evaluated with
this underlying quadrature, as in spectral element methods, we obtain
  a diagonal mass matrix without loosing the order of accuracy. We refer to this case as
  \textit{cubature} elements \cite{LIU2017cubature9order}. For this choice, the classical use of Runge--Kutta methods will provide the high order accuracy also for the time discretization.

Secondly, we will study the influence of the stabilization strategy. When solving \eqref{eq:conservation_law2D} with continuous finite elements
some additional stabilization operator is necessary
to enforce the $\mathbb L_2$ stability. Several stabilization techniques can be devised to introduce a level of
dissipation  
comparable to  that of 
discontinuous Galerkin  methods with upwind fluxes 
\cite{moura2020EigenanalysisGJP,moura2020SpatialEigenanalysis}.
Three approaches will be studied:
the streamline upwind Petrov--Galerkin (SUPG) stabilization \cite{article_supg1,article_supg2}, which is strongly consistent, but it is also introducing new terms in the mass 
matrix; 
 the
 continuous interior penalty (CIP) method \cite{article_cip5,article_cip3,article_cip4}, consisting in adding edge penalty terms 
 proportional to the jump of the first derivative of the solution; 
 the orthogonal subscale stabilization (OSS) \cite{CODINA20001579}, a term that penalizes the $\mathbb{L}_2$ projection of the gradient of the error within the elements. 
As the CIP stabilization, this technique does not affect the mass matrix, but it requires the solution of another linear system for the $\mathbb{L}_2$ projection.
In this respect, the choice of the finite element space  and of the quadrature  have enormous impact on the  cost of the method.
Note that the strategy to impose boundary  conditions also plays a major role in ensuring stability \cite{abgrall2020analysis,abgrall2021analysis},
but this will not be considered  here.

Our objective is to perform a fully discrete spectral  analysis on triangulations of the spatial domain to characterize the 
stability and accuracy of different combinations of approximation, quadrature, stabilization, and time stepping. In the linear case, this allows
to propose  optimal values of the CFL and stabilization parameters. Moreover, we provide some heuristic strategy to include in the analysis the impact
of residual based high order diffusion operators aiming at stabilizing discontinuities.
%
%
%
Numerical simulations for both linear and non-linear scalar problems, and for the shallow water system 
confirm the theoretical results, and allow to further investigate the impact of the discretization choices  on the performance of the schemes
and on their cost. 

The paper is organized as follows. In \cref{sec:numericalDiscretization} we describe the continuous Galerkin discretization, the stabilization techniques, the basis functions and the time integration techniques. In \cref{sec:fourierAnalysis} we introduce the Fourier analysis space definitions that lead to von Neumann analysis, we discuss some technical details on the passage from physical functions to Fourier modes for different meshes and we find the parameters for which the schemes are stable for some mesh configurations. In \cref{sec:app_fourier_visco}, we also propose to introduce a viscosity term in order to enforce stability when the previous von Neumann analysis reveals instabilities. In \cref{sec:validation_num-meshX} and \cref{sec:numericalSimulations} we test the found parameters on some linear and nonlinear problems, checking the order of accuracy and the computational times. Finally, in \cref{sec:conclusion2D} we derive some conclusions on the presented schemes and possible applications of the found results.

%
\section{Numerical discretization}\label{sec:numericalDiscretization}
In this section we describe the discretization of the hyperbolic conservation law \eqref{eq:conservation_law2D}. We consider a tessellation of the spatial domain $\Omega$ consisting of non overlapping (triangular) cells, which we denote by $\Omega_h\subset \R^2$. The generic element of the tessellation $\Omega_h$ will be denoted by $K$,  so that 
$\Omega_h=\bigcup K$. We denote the set of internal element boundaries (edges) of $\Omega_h$ by $\mathcal{F}_h$, using ${\sf f}$ for a general element. $h$ denotes the characteristic mesh size of $\Omega_h$. 
Despite of the fact that most of the discussion is performed for  the scalar case, most of it generalizes readily to systems. If a significant difference arises in this generalization, this will be explicitly discussed.

The discrete solution is sought in a continuous finite element space  $V_h^p = \lbrace v_h \in \mathcal{C}^0 (\Omega_h) : \, \restriction{v_h}{K} \in \P_p(K),\, \forall K\in \Omega_h \rbrace $. 
We will use nodal and modal finite elements, and  we will denote by $\varphi_j$ the basis functions associated to the degree of freedom $j$, so that  $V_h^p=\text{span}\left\{\varphi_j\right\}_{j\in\Omega_h}$ and we can write $$u_h(x)=\sum_{j\in\Omega_h} u_j \varphi_j(x),$$ where, with an abuse of notation, with $j\in\Omega_h$ we mean the set of degrees of freedom with support in $\Omega_h$. With a similar meaning, we will also use the notation $j\in K$ to mean the degrees of freedom with support on the cell $K$. 

The unstabilized CG approximation of \eqref{eq:conservation_law2D} reads: find $u_h\in V_h^p$ such that for any $v_h\in W_h\subset \mathbb{L}_2(\Omega_h):=\lbrace v: \Omega_h \to \R: \int_{\Omega_h} |v|^2 < \infty \rbrace$
\begin{equation}
\int_{\Omega_h} v_h \partial_t u_h dx  - \int_{\Omega_h} \nabla  v_h f(u_h)\; dx + \int_{\partial \Omega_h} v_h f(u_h)\cdot {\bf n} d \Gamma =0, \label{remi:1}
\end{equation}
where ${\bf n}$ is the normal to the boundary facing outward the domain.
The choice of $W_h$ will be based on $V_h$, but it may 
take different forms for different stabilizations.  

As already said, we will consider several stabilized variants of \cref{remi:1} which can be all formulated in the form:
find $u_h\in V_h^p$ that satisfies  
\begin{equation}
\int_{\Omega} v_h ( \partial_t u_h   + \nabla \cdot f(u_h)) dx + S(v_h,u_h)=0, \quad \forall v_h \in V^p_h \label{remi:2}
\end{equation}
where the flux term is written before the integration by part as  we will consider only continuous piecewise polynomials approximations, whose derivatives are integrable. Here, $S$ denotes a bilinear stabilization operator defined on $V^p_h\times V_h^p$. Several different choices for $S$ exist, and are discussed in detail in the following  sections.

\subsection{Stabilization Terms}\label{sec:stabilization}

\subsubsection{Streamline-Upwind/Petrov-Galerkin - SUPG} \label{SUPG}
The SUPG method was introduced in  \cite{article_supg01} (see also \cite{hst10, article_supg2}  and references therein) and is strongly consistent in the sense that it vanishes when replacing the discrete solution with the exact one.   
It can be written as a Petrov-Galerkin method replacing $v_h$ in  \cref{remi:1} with a test function belonging to the space 
\begin{equation}
W_h := \{ w_h:\quad w_h=v_h+\tau_K \nabla_u f(u_h) \cdot \nabla v_h; \quad v_h \in V_h^p \} . \label{eq_supg0}
\end{equation}
Here, $\nabla_u f(u_h) \in \R^{D\times D \times 2}$ is the Jacobian of the flux, $D$  the dimensions of the system, $\tau_K$ denotes a positive definite stabilization parameter with the dimensions of $D\times D$ 
that we will assume to be  constant for every element.
Although other definitions are possible, here  we will evaluate this parameter as
\begin{equation}
\tau_K = \delta h_K ( J_K ) ^{-1} 
\end{equation}
where $h_K$ is the cell diameter and  $J_K$  represents the norm of the flux Jacobian on a reference value of the element $K$. In the scalar case, $J_K = ||\nabla_u f(u)  ||_K$.

The final stabilized variational formulation of \eqref{remi:2} reads
\begin{equation}
\int_{\Omega} v_h \partial_t u_h \; dx + \int_{\Omega} v_h \nabla \cdot f(u_h) \; dx + \underbrace{\sum_{K \in \Omega} \int_{K} \big( \nabla_u f(u_h) \cdot \nabla v_h \big) \tau_K 
\left(\partial_t u_h + \nabla \cdot f(u_h) \right)\; dx}_{S(v_h,u_h)} = 0.
\label{eq_supg1}
\end{equation}

The main problem of this stabilization method is that it depends on the time derivative of $u$ and, hence, it does not maintain the structure of the mass matrix in most cases.

To characterize the accuracy of the  method, we can use the consistency analysis discussed \textit{inter alia}  in \cite[\S3.1.1 and \S3.2]{AR:17}. In particular, of a finite element polynomial approximation of degree $p$ we can 
easily  show that given a smooth exact solution $u^e(t,x)$, replacing formally  $u_h$ by the projection of $u^e$ on the finite element space,
we can write 
\begin{equation}
\begin{split}
\epsilon(\psi_h) :=& \Big|
\int_{\Omega} \psi_h \partial_t (u_h^e - u^e) \; dx - \int_{\Omega} \nabla \psi_h \cdot (\nabla f(u_h^e)-\nabla f(u^e))\; dx \\ + &
\sum_{K \in \Omega}\sum\limits_{l,m \in K} \dfrac{\psi_l - \psi_m}{k+1}  \int_{K} \big (\nabla_u f(u_h) \cdot \nabla \varphi_i) \tau_K
	\left(\partial_t (u_h^e - u^e) + \nabla \cdot ( f(u_h^e) -f(u^e)) \right)\; dx 
	\Big| \le C h^{p+1},
\end{split}
\label{eq_supg2}
\end{equation}
with $C$ a constant independent of $h$, for all  functions $\psi$ of class at least $\mathcal{C}^1(\Omega)$, of which  $\psi_h$  denotes the finite element projection.
A key point in this estimate is the strong consistency of the method allowing to subtract its formal  application to the exact solution (thus subtracting zero),
and obtaining the above expression featuring differences between the exact solution/flux and its evaluation on the finite element space. 
Preserving this error estimate  precludes the possibility of lumping the mass matrix, and in particular  the entries associated to the stabilization term.
This makes the scheme relatively inefficient when using standard explicit time stepping.

As a final note, for a linear flux \cref{eq:conservation_law1}, 
exact integration, with  $\tau_K = \tau$ and in the time continuous case,
a classical result is obtained for homogeneous boundary conditions by testing  with $v_h =u_h + \tau\, \partial_t u_h$ \cite{article_supg2}:
\begin{equation}
\begin{split}
\int\limits_{\Omega_h}\partial_t\left(\dfrac{u^2_h}{2}+\tau^2\dfrac{(\mathbf{a} \cdot  \nabla u_h)^2}{2}\right) +
\int\limits_{\Omega_h} \mathbf{a} \cdot  \nabla \left( \dfrac{u^2_h}{2}+\tau^2\dfrac{( \partial_t u_h)^2}{2}\right) = -\int\limits_{\Omega_h}\tau (\partial_tu_h+\mathbf{a} \cdot  \nabla u_h)^2.
\end{split}
\label{eq_supg3}
\end{equation}
For periodic, or homogeneous boundary conditions, this 
 shows that the  norm $|||u|||^2 :=\int_{\Omega_h} \dfrac{u^2_h}{2}+\tau^2\dfrac{(\mathbf{a} \cdot  \nabla u_h)^2}{2} dx$ is non-increasing.
The interested reader can refer to  \cite{article_supg2} for  the analysis of some  (implicit) fully discrete schemes.

\subsubsection{Note on the SUPG technique applied to non scalar problems}
The extension of the SUPG method to a non scalar problem is not straightforward. Here we used the following formulation. First, we
define the following system of dimension $D$:
\begin{equation}
    \left  \{
    \begin{array}{ll}
	    & \partial_t U + \nabla \cdot \mathcal{F}(U) = \mathbf{S}(U) \\
    	& \mathcal{F}=(F_1,F_2)
	\end{array}
    \right .    
    \label{eq:nonscalar_system}
\end{equation}
with $U \in \R^D$, $\mathcal{F}(U) \in \R^{2 \times D}$ and $\mathbf{S}(U)\in \R^D $. 
For example, in the results section we will consider the shallow water equations with $D =3$ which read
\begin{equation*}
	U=\begin{pmatrix}
	h \\
	hu \\
	hv
	\end{pmatrix}
	\quad
	F_1(U)= \begin{pmatrix}
	hu \\
	hu^2 +g\frac{h^2}{2} \\
	huv
	\end{pmatrix} \quad
	F_2(U)= \begin{pmatrix}
	hv \\
	huv	\\
	hv^2 +g\frac{h^2}{2} 
	\end{pmatrix}
	\quad \mbox{and} \quad \mathbf{S}(U)=\begin{pmatrix}
	0 \\
	-gh b_x \\
	-gh b_y
	\end{pmatrix}
\end{equation*}
where $ \mathbf{S}(U)$ is the source term given by a topography term.
Equation \eqref{eq:nonscalar_system} can also be written in its quasi-linear form
\begin{equation}
\partial_t U + \nabla_U \mathcal{F}(U) \cdot \nabla U = \mathbf{S}(U),
\end{equation}
where $\nabla_U \mathcal{F}(U_h) \in \R^{D\times D \times 2}$ is the Jacobian of the flux $\mathcal{F}(U_h)$.


Following the definition of the SUPG method and \cite[sec.~5]{RICCHIUTO20091071} 
we define 
a positive-definite stabilization matrix $\mathbf{\tau_K} \in \R^{D\times D}$
constant for every element $K$. Here this matrix is evaluated as \cite{RICCHIUTO20091071}
\begin{equation}
{\tau_K} = \delta h_K \left( \sum_{j\in S_K} \left| \nabla_U \mathcal{F}(\bar{U}_K) \cdot n_j \right| \right)^{-1}, 
\end{equation}
with $S_K$ the set of vertices of $K$, and $n_j$ the outward normal of the edge opposite to the vertex $j\in S_K$. $h_K$ is the cell diameter and $\nabla_u  \mathcal{F}(\bar{U}_K)$ represents the flux Jacobian of the  the average value of $U_h$ on the element $K$. \\


The SUPG stabilized formulation reads, for each equation of the system $i=1,\dots,D$ 
\begin{equation}
\int_{\Omega} v_h \left( \partial_t U_h + \nabla \cdot \mathcal{F}(U_h)-\mathbf{S}(U_h) \right)_{i} +
 \underbrace{ \left(\sum_{K \in \Omega} \int_{K} \big( \nabla v_h  \cdot \nabla_U \mathcal{F}(U_h) \big)  {\tau_K} \left(\partial_t U_h + \nabla \cdot \mathcal{F}(U_h)-\mathbf{S}(U_h) \right)\; dx \right)_{i}}_{S(v_h,U_h)_i}
 = 0,
\label{eq_supgSW1}
\end{equation}
where $(V)_{i}$ denotes the $i$-{th} component of a vector $V\in \R^D$.

\subsubsection{Continuous Interior Penalty - CIP}  \label{CIP}
Another stabilization technique, which maintains sparsity and symmetry of the Galerkin matrix, is the continuous interior penalty (CIP) method. It was developed by Burman and Hansbo originally in \cite{article_cip6} and then in a series of works \cite{article_cip5,article_cip3,article_cip4}. It can also be seen as a variation of the method proposed by Douglas and Dupont \cite{inbook_cip7}. 

The method stabilizes the Galerkin formulation by adding 
edge penalty terms 
proportional to the jump of the gradient of the derivatives of the solution across the cell interfaces. The CIP introduces high order viscosity to the formulation, allowing the solution to tend to the vanishing viscosity limit. This term is  
does not affect the structure of the mass 
matrix. The method reads
\begin{equation}
\int_{\Omega_h} v_h \partial_t u_h \; dx + \int_{\Omega_h} v_h \nabla \cdot f(u_h)\; dx+  \underbrace{ \sum_{{\sf f} \in\mathcal{F}_h}  \int_{\sf f} \tau_{\sf f} [\![n_{\sf f} \cdot \nabla v_h]\!] \cdot [\![n_{\sf f} \cdot \nabla u_h]\!] \; d\Gamma}_{S(v_h,u_h)} = 0,    
\label{eq_cip0}
\end{equation}
where $[\![\cdot]\!]$ denotes the jump of a quantity across a face $\sf f$, $n_{\sf f}$ is a normal to the face $\sf f$ and where $\mathcal{F}_h$ is the collection of internal boundaries,  and ${\sf f}$ are its elements. 
Although other definitions are possible, we evaluate the scaling parameter in the  stabilization as 
\begin{equation}
\tau_{\sf f} =  \delta   \,h_{\sf f}^2 \| \nabla_uf\|_{\sf f},
\label{eq_cip1}
\end{equation}
where $\| \nabla_uf\|_{\sf f}$ a reference value of the norm of the flux Jacobian on ${\sf f}$ and $h_{\sf f}$ a characteristic size of the mesh neighboring $\sf f$.

As stated above, a clear advantage of CIP is that it does not modify the mass matrix, allowing to obtain efficient schemes if a mass lumping strategy can be devised. 
On the other side, the stencil of the scheme increases as the jump of a degree of freedom interacts with cells which are not next to the degree of freedom itself (up to 2 cells distance).
Note that for higher order approximations  \cite{Burman2020ACutFEmethodForAModelPressure,larson2019stabilizationHighOrderCut}  suggest the use of jumps in higher derivatives to improve the stability of the method. However, here we consider the jump in the first derivatives in order to be able to apply the stability analysis and to study the influence of $\delta$ on the stability of the method. Some results might be definitely improved adding these stabilizations on higher derivatives.


The accuracy of CIP can be assessed with a consistency analysis as discussed in \cite[\S3.1.1 and \S3.2]{AR:17}. 
This consists in, formally substituting $u_h$ by the projection onto the finite element polynomial of degree $p$ space of $u^e$, a given smooth exact solution $u^e(t,x)$, we can show that for all  functions $\psi$ of class at least $\mathcal{C}^1(\Omega)$, of which  $\psi_h$  denotes the finite element projection, we have the truncation error estimate
\begin{equation}
\begin{split}
\epsilon(\psi_h) := \Big|
\int_{\Omega} \psi_h \partial_t (u_h^e - u^e) \; dx -& \int_{\Omega} \nabla \psi_h \cdot ( f(u_h^e)- f(u^e))\; dx \\ 
              +&  \sum\limits_{{\sf f}\in\mathcal{F}_h} \int\limits_{\sf f}\tau_{\sf f} [\![n_f \cdot \nabla \psi_h]\!] \cdot [\![n_f \cdot \nabla (u_h^e-u^e)]\!] 
	\Big| \le C h^{p+1},
\end{split}
\label{eq_cip2}
\end{equation}
with $C$ a constant independent of $h$. The estimate can be derived from standard approximation results applied to $u^e_h-u^e$ and to its derivatives, noting that $\tau_{\sf f}$ is an $\mathcal{O}(h^2)$, which allows to obtain the estimation with the right order.

The symmetry of the stabilization  allows to easily derive an energy stability estimate for the space discretized scheme only. In particular, for periodic boundary conditions and a linear flux
we can easily show that 
\begin{equation}
\begin{split}
\int\limits_{\Omega_h}\partial_t\dfrac{u^2_h}{2}= - \sum\limits_{{\sf f}\in\mathcal{F}_h}\int\limits_{\sf f} \tau_{\sf f} [\![n_f \cdot \nabla u_h]\!]^2,
\end{split}
\label{eq_cip3}
\end{equation}
which gives a bound in time on the $\mathbb{L}_2$ norm of the solution.  \\

Note that for higher than second order it may be relevant to consider additional penalty terms based on higher derivatives (see e.g. \cite{Burman2020ACutFEmethodForAModelPressure,burman2021weighted,paola_svetlana}). We did not do this in this work.

\subsubsection{Orthogonal Subscale Stabilization - OSS} \label{OSS}
Another symmetric stabilization approach is the  Orthogonal Subscale Stabilization (OSS) method. 
Originally introduced as Pressure Gradient Projection (GPS) in \cite{CODINA1997373} for Stokes equations,
it was extended to the OSS method in \cite{CODINA20001579,OSSCodinaBadia} for different problems with numerical instabilities, such as convection--diffusion--reaction problems.
This stabilization penalizes the fluctuations of the gradient of the solution with a projection of the gradient onto the finite element space. The method applied to \cref{remi:1} reads: find $u_h\in V_h^p$ such that $\forall v_h \in V_h^p$

\begin{equation}
    \left  \{
    \begin{array}{ll}
    	& \int_{\Omega_h} v_h \partial_t u_h \; dx  + \int_{\Omega_h} v_h \nabla \cdot f(u_h) \; dx  + \underbrace{\sum_{K \in \Omega_h} \int\limits_{K} \tau_K   \nabla  v_h \cdot(\nabla  u_h - w_h) \; dx}_{S(v_h ,u_h)}= 0,  \\
    	& \int_{\Omega_h} v_h   w_h\; dx - \int_{\Omega_h} v_h \nabla  u_h\; dx = 0. 
	\end{array}	
    \right .\label{eq_oss0}
\end{equation}

For this method, the stabilization parameter is evaluated as 
\begin{equation}
\label{eq_oss1}
\tau_K = \delta  h_K  \|\nabla_u f\|_K  .
\end{equation}
The drawback of this method, with respect to CIP, is the requirement of a matrix inversion to project the gradient of the solution in the second equation of \eqref{eq_oss0}. This cost can be alleviated by the choice of elements and quadrature rules if they result in a diagonal mass matrix, as it will be the case for \textit{Cubature} elements that we will describe below.

As before we can easily characterize the accuracy of this method.  The truncation error estimate for a polynomial approximation of degree $p$ reads in this case
	\begin{equation}
\begin{split}
\epsilon(\psi_h) := \Big|
\int_{\Omega_h} &\psi_h \partial_t (u_h^e - u^e) \; dx - \int_{\Omega_h} \nabla \psi_h \cdot (f(u_h^e)-f(u^e))\; dx \\ 
              +& \sum\limits_{K\in \Omega_h}\tau_K \int\limits_{K} \nabla \psi_h \cdot \nabla(  u^e_h - u^e  )
                          +  \sum\limits_{K \in\Omega_h} \tau_K\int\limits_{K} \nabla \psi_h \cdot (  \nabla u^e - w_h^e  )
	\Big| \le C h^{p+1},
\end{split}
\label{eq_oss2}
\end{equation}

where the last term is readily estimated using the projection error and the boundness of $\psi_h$ as
$$
\int_{\Omega_h} \psi_h  ( w^e_h-\nabla u^e)\; dx = \int_{\Omega_h} \psi_h (\nabla  u_h^e -    \nabla u^e ) =  \mathcal{O}(h^p).
$$

Finally, for a linear flux, periodic boundaries and taking $\tau_K=\tau$ constant along the mesh, we can test with  $v_h=u_h$ in the first equation of \eqref{eq_oss0}, and with $v_h=\tau w_h$ in the second one and sum up the result to get
\begin{equation}
\begin{split}
\int\limits_{\Omega_h}\partial_t\dfrac{u^2_h}{2}= - \sum\limits_{K} \int\limits_{K} \tau_K ( \nabla  u_h -  w_h)^2,
\end{split}
\label{eq_oss3}
\end{equation}
which can be integrated in time to obtain a bound on the $\mathbb{L}_2$ norm of the solution. \\

	The truncation consistency error analysis presented above for the three stabilization terms is completely formal and it does not comprehend an entire classical error analysis. These estimations tell us that the stabilization terms that we introduced are of the wanted order of accuracy and that they are usable to aim at the prescribed order of accuracy. This type of analysis has been already done for multidimensional problems inter alia in \cite{abgrallEntropyConservative}. More rigorous proof of error bounds with $h^{p+\frac12}$ estimates can be found in \cite{burman2021weighted} for the CIP.
We did not consider in this work projection stabilizations involving higher derivatives.

\subsection{Finite Element Spaces and Quadrature Rules} \label{sec:discretization2D}
In this section we describe three finite element polynomial approximation strategies used in the paper.
%
In particular, on 
a triangular element $K$ of $\Omega_h$,
we define in this section the restriction of the basis functions of $V_h^p$ on each element $K$, which are polynomials of degree at most $p$. We denote by $\{\varphi_1, \ldots, \varphi_N\}$  the basis functions and they will have degree at most $p$, and their definitions amounts to describe the degrees of freedom, i.e., the dual basis.

\subsubsection{\textit{Basic} Lagrangian equispaced elements}
On triangles, we consider Lagrange polynomials with degrees at most $p$: $\mathbb P^p=\lbrace \sum_{\alpha+\beta\leq p } c_{\alpha,\beta}x^\alpha y^\beta  \rbrace$. We define the barycentric coordinates $\lambda_i(x,y)$ which are affine functions on $\R^2$ verifying the following relations 
\begin{equation}
	\lambda_i(v_j)=\delta_{ij}, \quad \forall i,j=1,\dots,3,
\end{equation}
where $v_j=(x_j,y_j)$ are the vertexes of the triangle and, with an abuse of notation, they can be written in barycentric coordinates as $v_j=(\delta_{1j},\delta_{2j},\delta_{3j})$.
Using these coordinates, we can define the Lagrangian polynomials on equispaced points on triangles. The equispaced points are defined on the intersection of the lines $\lambda_j=\frac{k}{p}$ for $k=0,\dots,p$. 
A way to define the basis functions corresponding to the point $(x_\alpha,y_\alpha)=(\alpha_1/p,\alpha_2/p,\alpha_3/p)$ in barycentric coordinates, with $\alpha_i\in \{  0,\dots,  p\}$ and $\sum_i \alpha_i =1$, is in \cref{algo:baryCoordsLagrange}.
\begin{algorithm}
\begin{algorithmic}
	\Require Point $(x_\alpha,y_\alpha)=(\alpha_1/p,\alpha_2/p,\alpha_3/p)$ in barycentric coordinates
	\State $\varphi_\alpha(x)\gets 1 $
	\For{$i=1,2,3$}
	\For{$z=0,\dots,a_i$}
	\State $\varphi_\alpha(x,y) \gets \varphi_\alpha(x,y) \cdot (\lambda_i(x,y)-\frac{z}{p})$
	\EndFor
	\EndFor
\end{algorithmic}
	\caption{Lagrangian basis function in barycentric coordinates\label{algo:baryCoordsLagrange}}
\end{algorithm}

The polynomials so defined in a triangle form a partition of unity, but they have also negative values. This leads to negative or zero values of their integrals. This is problematic for some time discretization and we will see why. We will use these polynomials in combination with exact Gauss--Lobatto quadrature formulae for such polynomials and we will refer to them as \textit{Basic} elements.

\subsubsection{Bernstein polynomials}
Bernstein polynomials are as well a basis of $\mathbb P^p$ but they are not Lagrangian polynomials, hence, there is not a unique correspondence between point values and coefficients of the polynomials. Anyway, there exist a geometrical identification with the Greville points $(x_\alpha,y_\alpha)= (\alpha_1/p,\alpha_2/p,\alpha_3/p)$. Given a triplet $\alpha\in \N^3$ with $\alpha_i\in \llbracket 0,\dots,  p\rrbracket$ and $\sum_i \alpha_i=p$, the Bernstein polynomials are defined as 
\begin{equation}\label{eq:BernsteinPoly}
	\varphi_\alpha(x,y) = p!\prod_{i=1}^3 \frac{\lambda_i^{\alpha_i}(x,y)}{\alpha_i!}.
\end{equation}

Bernstein polynomials verify additional properties besides the one already cited for Lagrangian points. As before, they form a partition of unity, the basis functions are nonnegative in any point of the triangle, and so
their integrals are strictly positive. More precisely $$\int_K \varphi_\alpha = \frac{|K|}{S}, \qquad S= \# \left\lbrace \alpha \in \N^3 : |\alpha|_1 = p\right\rbrace. $$ These properties lead also to the fact that the value at each point is a convex combination of the coefficients of the polynomials, so that
it is easy to bound minimum and maximum of the function by the minimum and maximum of the coefficients. This has been used in different techniques to preserve positivity of the solution \cite{bacigaluppi2019posteriori,kuzmin2020subcellBernstein}. We will use these polynomials with corresponding high order accurate quadrature formulae. We will denote these elements with the symbol $\mathbb B^p$ and we refer to them as \textit{Bernstein} elements.

\subsubsection{\textit{Cubature} elements} \label{sec:discretization_cubature2D}
\begin{figure}[H]
\centering
\includegraphics[width=0.6\linewidth]{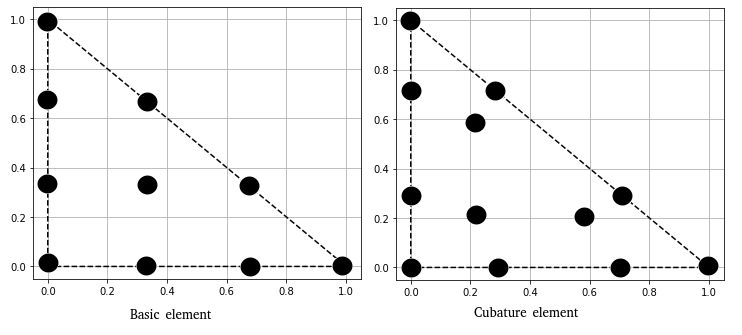}
\caption{Comparison of the equispace repartition at left and the cubature repartition at right for elements of degree $p=3$.} \label{fig:P3-basicVScub}
\end{figure}

Contrary to the work done in 1D \cite{michel2021spectral}, the extension of Legendre--Gauss--Lobatto points which minimize the interpolation error do not exist for the triangle. They have to be computed numerically such as \textit{Fekete} points \cite{article_fekete1, article_fekete2_SHERWIN1995189, article_fekete}. The problem of this approach is that it requires as classical finite elements the inversion of a sparse global mass matrix.\\
\textit{Cubature} elements were introduced by G. Cohen and P. Joly in $2001$ \cite{article_cubature_2001} for the wave equation (second order hyperbolic equation), and are an extension of Lagrange polynomials with the goal of optimizing the underlying quadrature formula error. We will denote the with the symbol $\tilde {\mathbb P}^p$ and they will be contained in another larger space of Lagrange elements, i.e., $\mathbb P ^p \subseteq \tilde {\mathbb P}^p \subseteq \mathbb P^{p'}$, with $p'$ the smallest possible integer.
Similar techniques have been used to minimize the interpolation error \cite{article_fekete1, article_fekete2_SHERWIN1995189, article_fekete}. The objective of these polynomials is to use the points of the Lagrangian interpolation of the polynomials as quadrature points. This means that the obtained quadrature is $\int_K f(x,y) = \sum_{\alpha} \omega_\alpha f(x_\alpha,y_\alpha)$, where $\int_K \varphi_\alpha = \omega_\alpha$ and $\varphi_\alpha(x_\beta,y_\beta) = \delta_{\alpha \beta}$. This approach can be considered an extension of the Gauss--Lobatto quadrature in 1D for non Cartesian meshes. The biggest advantage of this approach is to obtain a diagonal mass matrix. 
The drawback is that one needs to increase the number of basis function inside one element to obtain an accurate enough quadrature rule. In our work, we propose to extend this approach to first order hyperbolic equations. A successful extension to elliptic problem is proposed in \cite{pasquetti:hal-01954133}.
A comparison between the equispace repartition and the \textit{Cubature} repartition for elements of degree $p=3$ is shown in \cref{fig:P3-basicVScub}.

{For completeness we detail further the construction of the basis functions.} The challenges of this approach are the following:
\begin{itemize}
	\item Obtain a quadrature which is highly accurate, at least $p+p'-2$ order accurate \cite{book_ciarlet};
	\item Obtain positive quadrature weights $\omega_\alpha>0$ for stability reasons \cite{tordj1995};  
	\item Minimize the number of basis functions of $\tilde{\mathbb{P}^p}$;
	\item The set of quadrature points has to be $\tilde{\P}^p$-unisolvent;
	\item The number of quadrature points of edges as to be sufficient ensure the conformity of the finite element.
\end{itemize}
The optimization procedure that lead to these elements consists of several steps where the different goals are optimized one by one. The optimization strategy exploits heavily the symmetry properties that the quadrature point must have. 

For $p=1$ the \textit{Cubature} elements do not differ from the \textit{Basic} elements but in the quadrature formula. For $p=2$ the \textit{Cubature} elements introduce an other degree of freedom at the center of the triangle, leading to 7 quadrature points and basis functions per element. For $p=3$ the additional degree of freedom in the triangle are 3, leading to 13 basis functions per triangle.
All the details of such elements can be found in \cite{article_cubature_2001,jund2007arbitrary}. We provide 
in \cref{sec:appendix_cohen_bf} the detailed expressions
of the polynomials used in this work. We will use the symbol $\tilde{\mathbb{P}}^p$ and the name \textit{Cubature} elements to refer to them. \\

Other elements such as \textit{Fekete-Gauss} points \cite{article_cubature_2006,pasquetti:hal-01589136} exist in the literature. They are optimized to interpolate and integrate with high accuracy. However, it is shown that they require more computing time to achieve similar results than cubature points for high order of accuracy.

\subsection{Time integration}
The spatial discretization leads to
a coupled system of ordinary differential equation which can be written as 
\begin{equation}\label{eq:linear_system}
  \mass  \dfrac{dU}{dt}  = \rhs(t) 
\end{equation}
where $U$ 
is the vector of all the degrees of freedom on all the domain, $\mass$ and $\rhs$ are the global mass matrix and right-hand side terms obtained through the discretization of the previous section with some finite elements and stabilization terms.
We remark that $\mass$ is diagonal only in the case of the \textit{Cubature} elements without the SUPG stabilization, while, for all other choices, it is a sparse non--diagonal matrix.

In the following, we describe two different time integration method:
explicit Runge--Kutta (RK) methods and their strong stability preserving (SSP) variants; and the Deferred Correction (DeC) algorithm, which allows to avoid the mass matrix inversion through the correction iterations. 
\subsubsection{Explicit Runge--Kutta and Strong Stability Preserving Runge--Kutta schemes} 
Runge--Kutta time integration methods are one step methods consisting in $S$ stages defined by 
\begin{equation}
\begin{split}
	&U^{(0)}:=U^n,\\
	&U^{(s)}:=U^n + \Delta t \sum_{j=0}^{s-1}\alpha_j^s \mass^{-1} \rhs(U^{(j)})\quad s=1,\dots, S,\\
	&U^{n+1}:= U^n +  \Delta t  \sum_{s=0}^S \beta_s \mass^{-1} \rhs(U^{(s)}).
	\end{split}\label{eq:RK}
\end{equation}
Here, we use for the solution the superscript $n$ to indicate the timestep and the superscript in brackets $(s)$ to denote the stage of the method. The coefficients $\alpha_j^s$ and $\beta_j^s$ can be defined in many different ways.
In particular, we will refer to Heun's method with RK2, to Kutta's method with RK3 and the original Runge--Kutta fourth order method as RK4. The respective Butcher tables can be found in \cref{sec:timeCoefficients} in \cref{tab:Butcher}, see \cite{butcher08numODE}.

A subset of the RK methods are the SSPRK introduced in \cite{shu-1988}. They consist in convex combinations of forward Euler steps, and can be rewritten as follows
\begin{equation}\label{eq:SSPRKformula}
	\begin{split}
	&U^{(0)}:=U^n,\\
	&U^{(s)}:=\sum_{j=0}^{s-1} \left( \gamma_j^s U^{(j)} + \Delta t \mu_j^s \mass^{-1} \rhs(U^{(j)}) \right) \quad s=1,\dots, S,\\
	&U^{n+1}:= U^{(S)} ,
	\end{split}
\end{equation}
with $\gamma_j^s, \mu_j^s\geq 0$ for all $j,s=1,\dots, S$. 
We will consider here the second order 3 stages SSPRK(3,2) presented by Shu and Osher in \cite{shu-1988}, the third order SSPRK(4,3) presented in \cite[Page 189]{Ruuth-2006}, and the  fourth order SSPRK(5,4) defined in \cite[Table 3]{Ruuth-2006}.
For complete reproducibility of the results, we put all their Butcher's tableaux in \cref{sec:timeCoefficients} in \cref{tab:ButcherSSPRK}.

\subsubsection{The \textit{Deferred Correction} scheme}
Deferred Correction methods were  introduced in \cite{ecdd148d37dd401ca8d415ae18d3ecf2} as explicit time integration methods for ODEs, but soon implicit \cite{DeC_ODE}, linearly implicit positivity preserving \cite{oeffner_torlo_2019_DeCPatankar} versions and extensions to PDE solvers \cite{DeC_2017} were studied.
In particular, in \cite{DeC_2017,DeC_AT,paola_svetlana,abgrall2019reinterpretation} the DeC is used in a different formulation for finite element methods and it introduces two operator through which it is possible to use a diagonal mass matrix without losing the accuracy order.
This is only achievable when the lumped matrix (defined as the sum on the rows of the full mass matrix) has only positive values on its diagonal. Hence, the use of \textit{Bernstein} polynomials is recommended in \cite{DeC_2017}, but also \textit{Cubature} elements can serve the purpose.

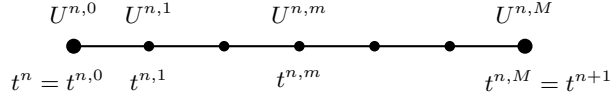
\begin{figure}[h]
	\small	\centering
	\begin{tikzpicture}
	\draw [thick]   (0,0) -- (6,0) node [right=2mm]{};
	\fill[black]    (0,0) circle (1mm) node[below=2mm] {$t^n=t^{n,0} \,\, \quad$} node[above=2mm] {$U^{n,0}$}
	(1,0) circle (0.7mm) node[below=2mm] {$t^{n,1}$} node[above=2mm] {$U^{n,1}$}
	(2,0) circle (0.7mm) node[below=2mm] {}
	(3,0) circle (0.7mm) node[below=2mm] {$t^{n,m}$} node[above=2mm] {$U^{n,m}$}
	(4,0) circle (0.7mm) node[below=2mm] {}
	(5,0) circle (0.7mm) node[below=2mm] {}
	(6,0) circle (1mm) node[below=2mm] {$\qquad t^{n,M}=t^{n+1}$} node[above=2mm] {$U^{n,M}$};
	\end{tikzpicture}
	\caption{Subtimesteps inside the time step $[t^n,t^{n+1}]$} \label{fig_DeC-time_disc}
\end{figure}

Consider a discretization of each timestep into $M$ subtimesteps
as in  \cref{fig_DeC-time_disc}.
For each subtimestep we define a high order approximation of the integral form of the ODE \eqref{eq:linear_system} from $t^{n,0}$ to $t^{n,m}$, i.e.,
\begin{equation}\label{eq:L2}
	\mass \left( U^{n,m} -U^{n,0} \right) - \int_{t^{n,0}}^{t^{n,m}} \rhs (U(s)) ds \approx \L^2(\UU)^m:= \mass \left( U^{n,m} -U^{n,0} \right) - \Delta t \sum_{z \in \llbracket 0, M \rrbracket } \rho_{z}^m \rhs (U^{n,z}) = 0,
\end{equation}
with $\UU=\left( U^{n,0},\dots, U^{n,M} \right)$. Moreover, the quadrature rule in time uses the subtimesteps $t^{n,m}$ as quadrature points. The corresponding weights $\rho^{m}_z$ for every different subinterval are defined by Lagrangian basis functions in these subtimesteps (see  \cite{DeC_2017,DeC_AT,paola_svetlana} for details). 
The algebraic system $\L^2(\UU^*)=0$ is in general implicit and nonlinear 
and, in order not to recast to nonlinear solvers, the DeC procedure approximates the solution of $\L^2(\UU^*)=0$ by successive iterations relying on a low order easy--to--invert operator $\L^1$.  This operator is typically a first order forward Euler approximation with a lumped mass matrix, i.e.,
\begin{equation}\label{eq:L1}
\mass \left( U^{n,m} -U^{n,0} \right) - \int_{t^{n,0}}^{t^{n,m}} \rhs (U(s)) ds \approx \L^1(\UU)^m:= \diag \left( U^{n,m} -U^{n,0} \right) - \Delta t \beta^m \rhs (U^{n,0}) = 0.
\end{equation}
Here, $\diag$ denotes a diagonal matrix obtained from the lumping of $\mass$, i.e., $\diag_{ii}:=\sum_{j} \mass_{ij}$, and $\beta^m:= \frac{t^{n,m}-t^{n,0}}{t^{n+1}-t^n}$. The values of the coefficients $\beta^m$ and $\rho^m_z$ for equispaced subtimesteps can be found in \cref{sec:timeCoefficients}. 
Denoting with the superscript $(k)$ index the iteration step, we describe the DeC algorithm as
\begin{subequations}
\begin{align}
	&U^{n,m,(0)}:=U^n & m=0,\dots,M,\\
	&U^{n,0,(k)}:=U^n  & k=0,\dots, K,\\
	&\L^1(\UU^{(k)})=\L^1(\UU^{(k-1)})-\L^2(\UU^{(k-1)})& k=1,\dots, K, \label{eq:DeCUpdate}\\
	&U^{n+1}:=U^{n,M,(K)}.&
\end{align}
\end{subequations}

It has been proven \cite{DeC_2017} that if $\L^1$ is coercive, $\L^1-\L^2$ is Lipschitz with a constant $\alpha_1 \Delta t >0$ and the solution of $\L^2(\UU^*)=0$ exists and is unique, then, the method converges with an error of $\mathcal{O}(\Delta t^K)$. Hence, choosing $K=M+1$ we obtain a $K$-th order accurate scheme.

Relying only on the inversion of the low order operator, the method gets rid of the computational costs of the solution of the linear systems, leaving in the right hand side the mass matrix of the $\L^2$ operator, that should not be inverted.
The only requirement that is necessary for the DeC approach is the invertibility of the lumped mass matrix $\diag$, 
which limits its application to spatial elements which possess this property.\textit{Basic} Lagrange polynomials do not guarantee such constraint already for degree 2. Hence, only other polynomials as \textit{Bernstein} and \textit{Cubature} can be used in combination with DeC. 

Finally, for the following analysis we note that the DeC method can be cast in a form similar to a Runge--Kutta method by rewriting \cref{eq:DeCUpdate} as
\begin{equation}\label{eq:DeCasSSPRK}
	U^{n,m,(k+1)}=U^{n,m,(k)} - \diag^{-1} \mass \left(U^{n,m,(k)}-U^{n,0,(k)}\right) +\sum_{j=0}^M \Delta t \rho_{j}^m \diag^{-1}\rhs(U^{n,j,(k)}). 
\end{equation}
Comparing with the system of equations \eqref{eq:SSPRKformula}, we can immediately define the SSPRK coefficients associated to DeC as $\gamma^{m,(k+1)}_{m,(k)}=\mathbb{I}-\diag ^{-1} \mass$ with $\mathbb{I}$ the identity matrix,  $\gamma^{m,(k+1)}_{0,(0)}=\diag ^{-1} \mass$, $\mu^{m,(k+1)}_{r,(k)}=\rho^m_r$ for $m,r=0,\dots,M$ and $k=0,\dots,K-1$ and instead of the mass matrix, we use the diagonal one.

\begin{remark}[DeC with SUPG]
	The iterative procedure of the DeC method allows even to overcome the difficulties that some implicit stabilization as the SUPG has. Indeed, the SUPG stabilization term can be added only to the $\L^2$ operator, keeping the high order accuracy of this operator. Since the $\L^2$ operator is applied to the previously computed iteration, all the terms of the SUPG, included the time derivative of $u$ in \cref{eq_supg1}, can be explicitly computed on $U^{(k-1)}$, keeping then the diagonal mass matrix for the whole scheme.
\end{remark}
%
\section{Fourier analysis}\label{sec:fourierAnalysis}
\subsection{Preliminaries and time continuous analysis} \label{sec:fourier}

In order to study the stability and the dispersion properties of the previously presented numerical schemes, we will perform a dispersion analysis on the linear advection problem with periodic boundary conditions: 
\begin{equation}
    \partial_t u(t,\mathbf{x}) + \mathbf{a} \cdot \nabla u(t,\mathbf{x}) = 0, \quad \mathbf{a}\in \R^2, \quad (t,\mathbf{x}) \in \R^+ \times \Omega,
    \label{eq_disp_1}
\end{equation}
with $\Omega = [0,1]\times[0,1]$. For simplicity, we consider $\mathbf{a} = (\cos( \Phi), \sin ( \Phi))$ with $\Phi \in [0,2\pi]$. 
We then introduce the ansatz
\begin{align}
    & u_h(\mathbf{x} , t) = Ae^{i(\mathbf{k} \cdot \mathbf{x} - \xi t)} = Ae^{i(\mathbf{k}\cdot \mathbf{x}-\omega t)}e^{\epsilon t} \label{eq_disp_uex1}  \\
     \mbox{with} \quad & \xi = \omega + i \epsilon, \quad i=\sqrt{-1}, \quad \mathbf{k}=(k_x,k_y)^T. 
\end{align}
Here, $\epsilon$ denotes the damping rate, while the wavenumbers are denoted by $\mathbf{k}=(k_x,k_y)$, with  $k_x=2\pi/L_x$ and $k_y=2\pi/L_y$ with $L_x$ and $L_y$ the wavelengths in $x$ and $y$ directions respectively. 
The phase velocity $\mathbf{c}$ can be defined from 
\begin{equation}
    \mathbf{c}\cdot\mathbf{k} =   \omega 
    \label{eq_PDE_velocity}
\end{equation}
and represents the celerity with which waves propagate in space. It is in general a function of the wavenumber. Substituting 
\eqref{eq_disp_uex1} in the advection equation \eqref{eq_disp_1} for an exact solution we obtain that
\begin{align}
	\label{eq:Phi_def}
	\omega = \mathbf{k}\cdot \mathbf{a} \,,\quad
    \mathbf{c} =  \mathbf{a} \quad \mbox{and} \quad \epsilon = 0.
\end{align}
In other words
\begin{align}
     u_h(\mathbf{x} , t) = Ae^{i\mathbf{k} \cdot (\mathbf{x} - \mathbf{a} t)}  \,.
\end{align}
The objective of the next sections is to provide the semi- and fully-discrete equivalents of the above relations for the finite element methods introduced earlier.
We will consider  polynomial degrees up to 3, for all combinations  of stabilization methods and time integration techniques. This will also allow to investigate the parametric stability with respect to
the time step (through the \CFL~number) and stabilization parameter $\delta$. In practice, for each choice we will evaluate the accuracy of the discrete
approximation   of  $\omega$ and $\epsilon$,  and we will provide conditions for the non-positivity of the damping $\epsilon$, i.e., the von Neumann stability of the method.

\subsection{The eigenvalue system}
\label{sec:fourier_space}

\begin{figure}
	\centering
	\includegraphics[width=0.3\linewidth]{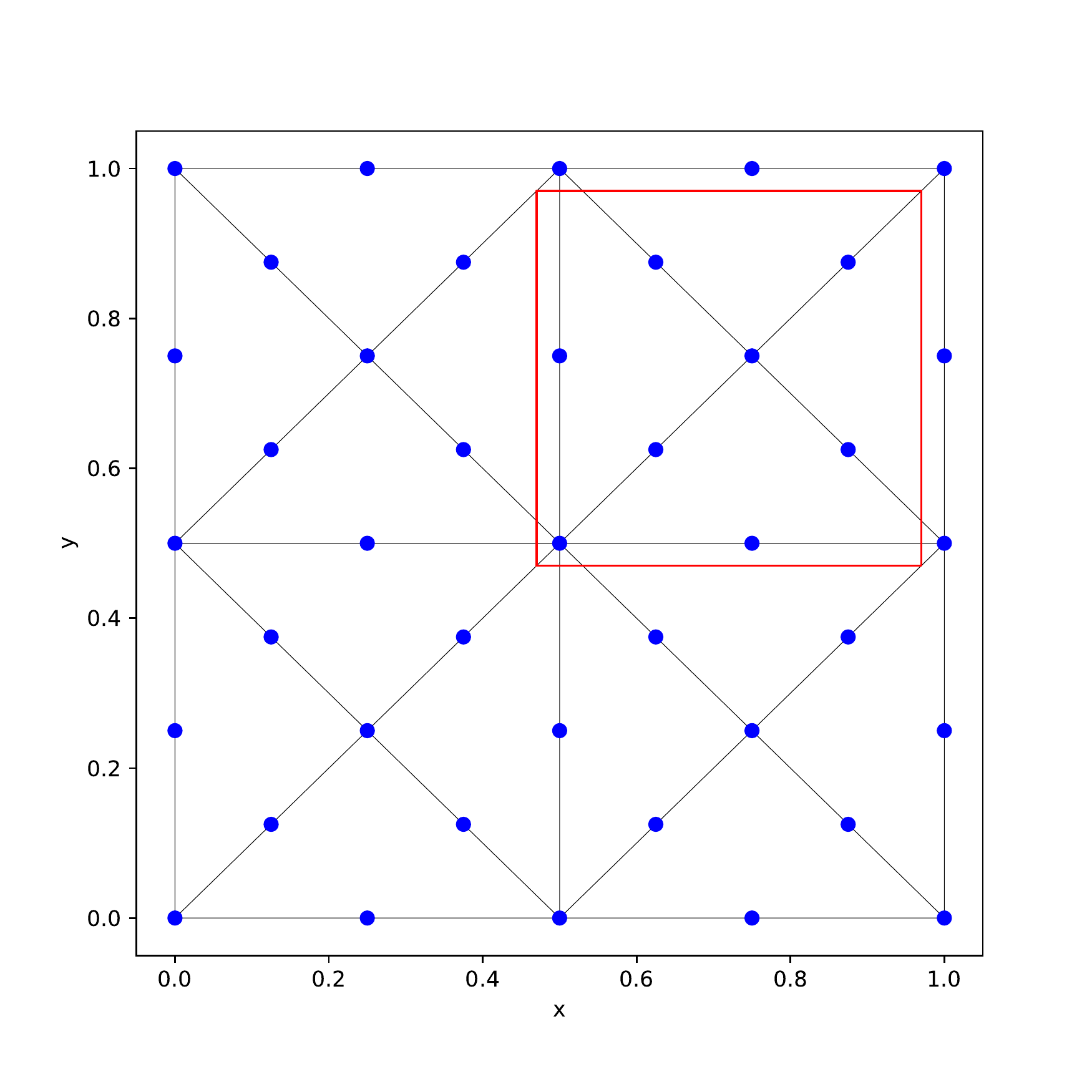}
	\includegraphics[width=0.3\linewidth]{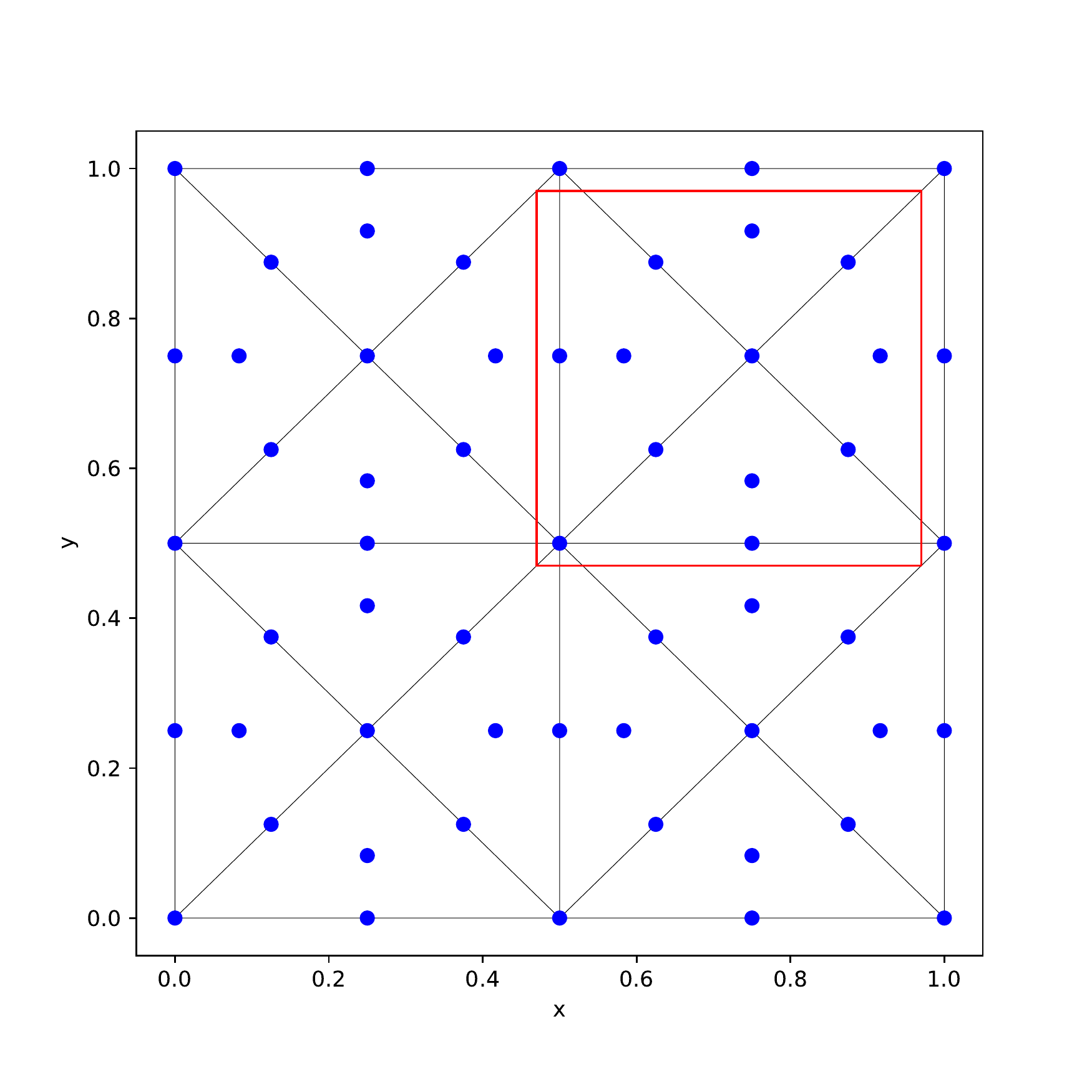}
	\caption{The \textit{X} type triangular mesh. At left, the \textit{Basic} finite element discretisation with $\P_2$ elements. At right, the grid configuration for $\Tilde{\P}_2$ \textit{Cubature} elements. The red square represents the periodic elementary unit that contains the degrees of freedom of interest for the Fourier analysis} \label{fig_Xmesh_P2}
\end{figure}

The Fourier analysis for numerical schemes on the periodic domain  is based on a discrete Parseval theorem. 
%
%
Thanks to this theorem, we can study the amplification and the dispersion of the basis functions of the Fourier space. 
The key ingredient of this study is the repetition of the stencil of the scheme from one cell to another one. In particular, using the ansatz \eqref{eq_disp_uex1} we can write local equations coupling degrees of freedom belonging to neighbouring cells through a multiplication by factors $e^{i\theta_x}$ and $e^{i\theta_y}$ representing the shift in space along the oscillating solution. The dimensionless  coefficient
\begin{equation}\label{eq_theta}
\theta_x:=  k_x\Delta x\, \text{ and }\, \theta_y:=  k_y\Delta y
\end{equation}
are the discrete reduced wave numbers which naturally appear all along the analysis. Here, $\Delta x$ and $\Delta y$ are defined by the size of the elementary periodic unit that is highlighted with a red square as an example in \cref{fig_Xmesh_P2}.

\begin{figure}
	\centering
	\subfigure[\textit{Basic} $\P_2$ elements]{
		\includegraphics[width=0.35\linewidth]{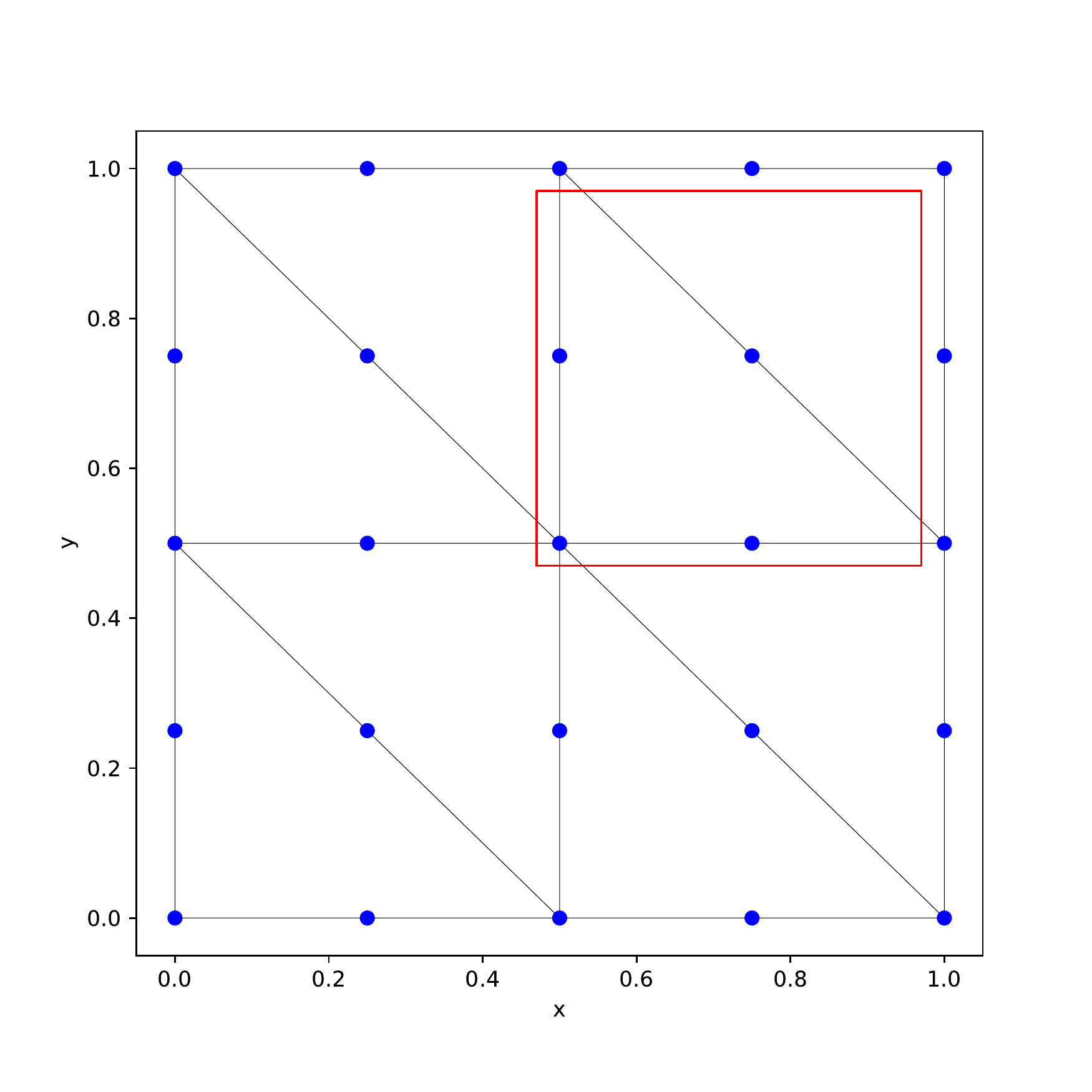}}
	\subfigure[\textit{Cubature} $\Tilde{\P}_2$ elements]{
		\includegraphics[width=0.35\linewidth]{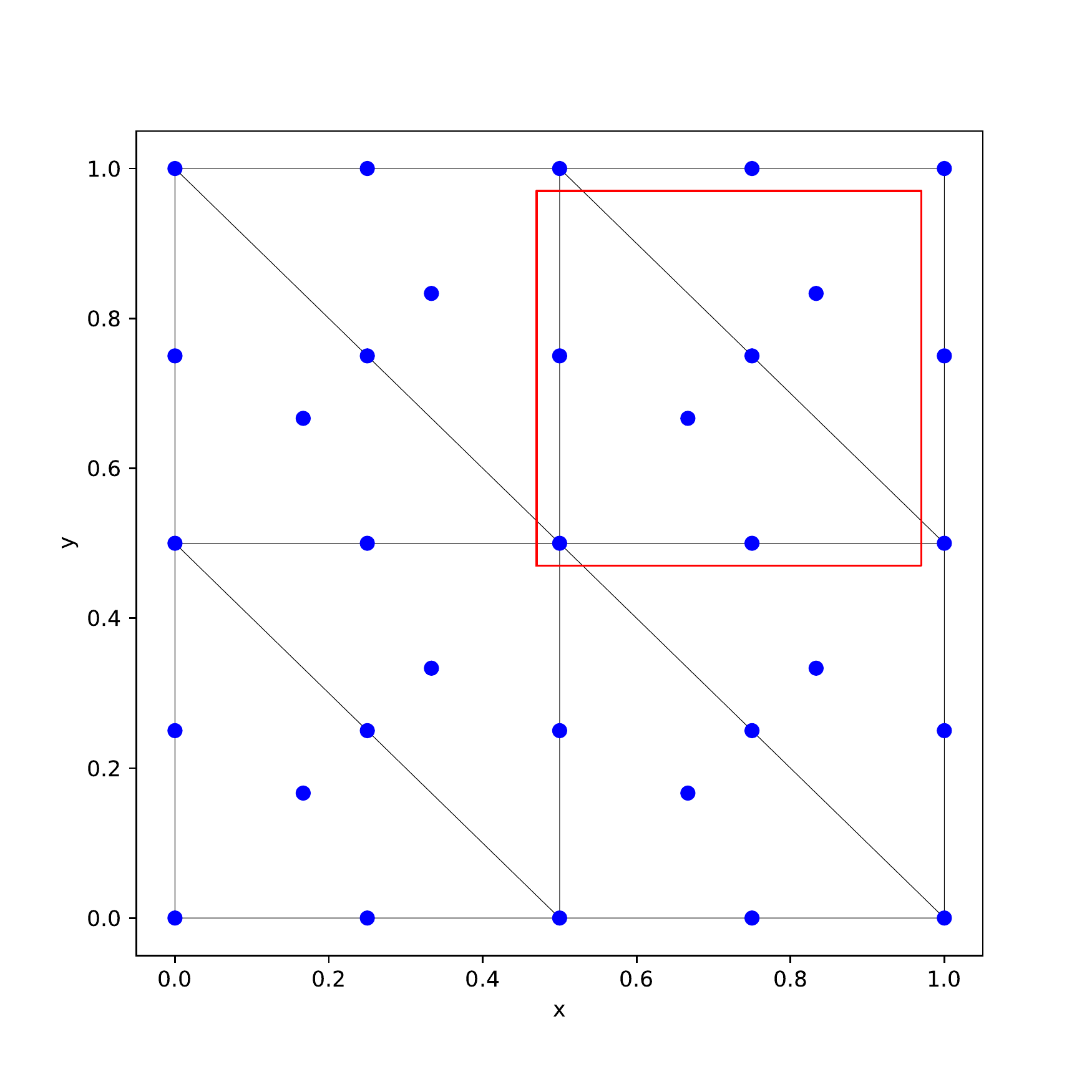}}
	\caption{The \textit{T} type triangular mesh with degrees of freedom in blue and periodic unit in the red square for the Fourier analysis} \label{fig_Tmesh_P2}
\end{figure}

Formally replacing the ansatz in the scheme we end up with a dense  algebraic problem of dimension $N_{dof}$, where $N_{dof}$ is the number of all the degrees of freedom in the mesh. The obtained system with dimension $N_{dof}$ in the time continuous case reads
\begin{equation}
 \text{Equations \eqref{eq_disp_1}  and \eqref{eq_disp_uex1}} \quad \Rightarrow \quad - i\xi \mass \mathbf{U}    + \mathbf{a} \cdot ( \mathcal{K}_x \mathbf{U} ,\mathcal{K}_y \mathbf{U})  + \delta \mathbb{S}\mathbf{U} = 0 \label{eq:eigenvalueProblem}
\end{equation}
\begin{equation}
\hspace*{-1cm}
    \mbox{with} \quad (\mass)_{ij} = \int_{\Omega} \phi_i \phi_j dx, \qquad (\mathcal{K}_{x})_{ij} = \int_{\Omega} \phi_i \partial_x \phi_j dx, \qquad (\mathcal{K}_{y})_{ij} = \int_{\Omega} \phi_i \partial_y \phi_j dx
\end{equation} 
with $\phi_j$ being any finite element basis functions, $ \mathbf{U} $ the array of all the degrees of freedom and $\mathbb{S}$ being the stabilization matrix defined through one of the stabilization techniques of \cref{sec:stabilization}. Although system \eqref{eq:eigenvalueProblem}   is in general a global eigenvalue problem, we can reduce its complexity by exploiting more explicitly the ansatz \eqref{eq_disp_uex1}.
The choice of the mesh is crucial in order to exploit the ansatz and to find a unit block that repeats periodically in space. Hence, we must consider structured periodic meshes and we will focus, in particular, on two types of meshes. The first one is the $X$-mesh that is depicted in \cref{fig_Xmesh_P2} and the second one is the $T$-mesh depicted in \cref{fig_Tmesh_P2}. In those pictures also the distribution of some $\P_2$ elements are represented as an example. 

More precisely, as it is done in \cite{article_fekete2_SHERWIN1995189} we can introduce elemental vectors of unknowns $\widetilde{\mathbf{U}}_{Z_{ij}}$, where $Z_{ij}$ is the stencil denoted by the red square in \cref{fig_Xmesh_P2}, which repeats periodically on the domain. So that $\widetilde{\mathbf{U}}_{Z_{ij}}$, for continuous finite elements, is an array of $d$ degrees of freedom inside a periodic unitary block $Z_{ij}$, excluding two boundaries (one on the top and one on the right  for example). This number depends on the chosen (periodic) mesh type and on the elements. As an example, in \cref{fig_Xmesh_P2} we display for the \textit{X }type mesh the periodic elementary unit (in the red square) with \textit{Basic} and cubature degrees of freedom with $p=2$.
In the \textit{X} mesh for \textit{Basic} elements $p=2$ we have $d=8$, while for \textit{Cubature} $p=2$ we have $d=12$. 
Using the periodicity of the solution and the ansatz \eqref{eq_disp_uex1} and denoting by $Z_{i\pm1,j\pm1}$ the neighboring elementary units, we can write the neighboring degrees of freedom by
\begin{equation}
\widetilde{\mathbf{U}}_{Z_{i\pm 1,j}} = e^{\pm\theta_x}\widetilde{\mathbf{U}}_{Z_{i,j}}, \qquad \widetilde{\mathbf{U}}_{Z_{i,j\pm 1}} = e^{\pm\theta_y}\widetilde{\mathbf{U}}_{Z_{i,j}},
\label{eq_ad_fourier0}
\end{equation}
and by induction all other degrees of freedom of the mesh.
This allows to show that the system \eqref{eq:eigenvalueProblem} is equivalent to a  compact system of dimension $d$ (we drop the subscript $_K$ as they system is equivalent for all cells)
\begin{equation}
-i\xi \widetilde{\mass}   \widetilde{\mathbf{U}}  + a_x \widetilde{\mathcal{K}}_x \widetilde{\mathbf{U} } +a_y \widetilde{\mathcal{K}}_y \widetilde{\mathbf{U} } + \delta \widetilde{\mathbb{S}} \widetilde{\mathbf{U} }  =0,
\label{eq_ad_fourier}
\end{equation}
where the matrices $\widetilde{\mass}$, $ \widetilde{\mathcal{K}}_x$, $ \widetilde{\mathcal{K}}_y$ and $\widetilde{\mathbb{S}}$ are readily obtained from the elemental discretization matrices by using  Equations \eqref{eq_ad_fourier0}.

For the discrete Parseval theorem, we know that the norm or the reduced variable $\widetilde{\mathbf{U}}$ is equivalent to the norm of the discrete vector $\mathbf{U}$. Hence, studying the amplification factor of the two is equivalent.

We apply the same analysis to stabilized methods. The interested reader can access all 2D dispersion plots online \cite{TorloMichel2021git}. From the plot we can see that the increase in polynomial degree
provides the  expected large reduction in dispersion error, while retaining a small amount of numerical dissipation, which permits the damping of \textit{parasite} modes.

\begin{figure}[H]
     \centering
     \includegraphics[width=0.9\textwidth]{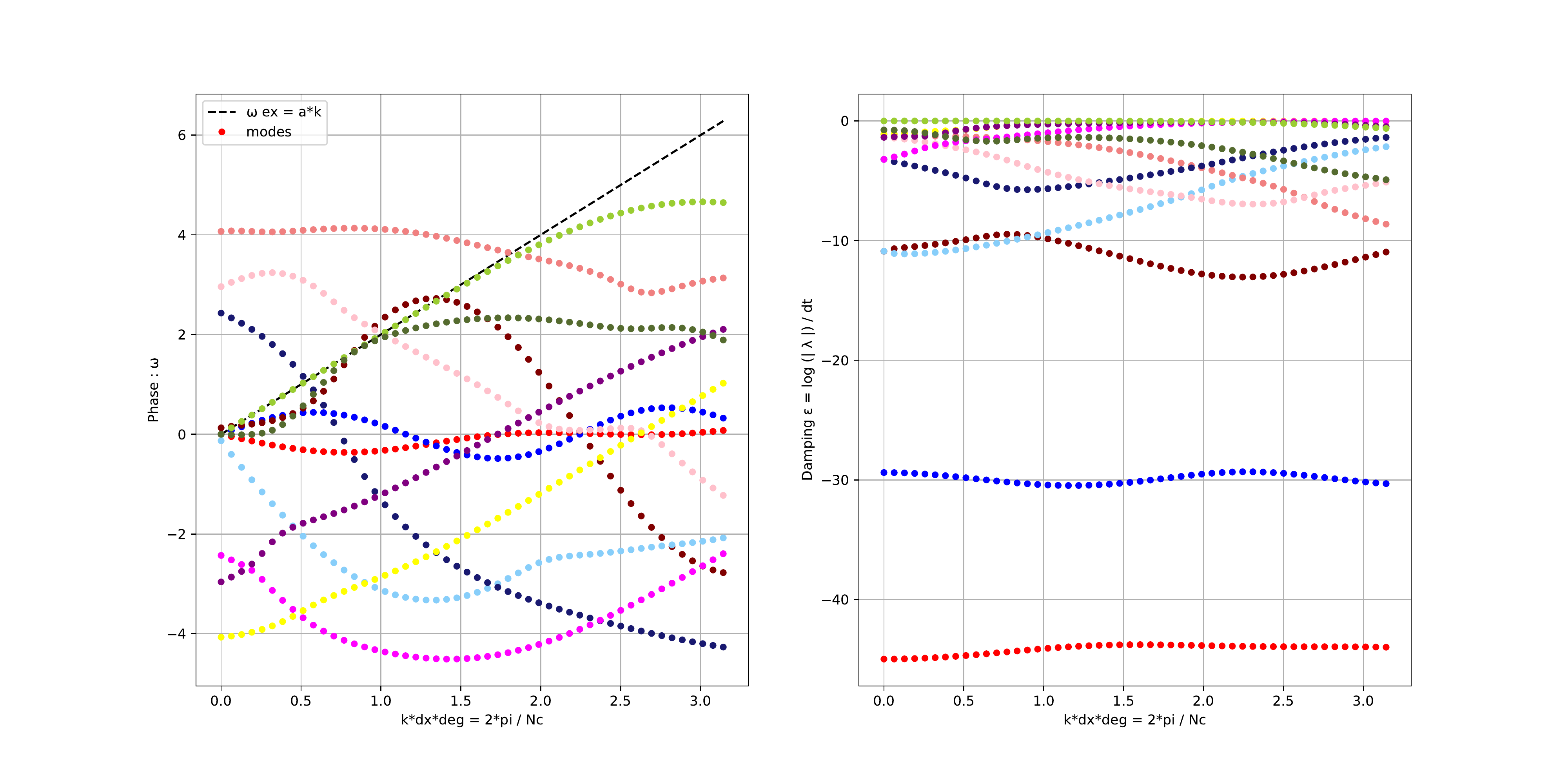}
     \caption{Dispersion curves related to the 12 modes of $\widetilde{\mathbf{U}}_{Z_{ij}}$ of the system given by \textit{Cubature} $\TP_2$ elements, the CIP stabilization technique, and a wave angle $\theta = 5 \pi / 4$ on an \textit{X} mesh.\\
     Phases $\omega$ (left) and amplifications $\epsilon$ (right).}
     \label{fig:disp_cohen_P2_CIP}
\end{figure}

An example of dispersion curves is given in \cref{fig:disp_cohen_P2_CIP}. The method used \textit{Cubature} $\TP_2$ elements, the CIP stabilization technique, and a wave angle $\theta = 5 \pi / 4$. We here show all 12 \textit{parasite} modes (see \cref{fig_Xmesh_P2}). The \textit{principal} mode of this system is represented in green. This figure also show the complexity of the analysis because of the number of modes to consider. 

We summarize the number of modes for the \textit{X} mesh in \ref{tab:number_of_modes_2D-meshX}. A representation of each mesh is done in \ref{sec:app_mesh4fourier} for element of degree $p=2$ and $3$.
\begin{table}[H] 
\small  
 \begin{center} 
		\begin{tabular}{| c || c | c | c | }  
	     \hline 
	     Element  & $\mathbb{P}_1$ & $\mathbb{P}_2$ & $\mathbb{P}_3$    \\ \hline \hline 
       Cub.              &  2 & 12 & 26  \\   \hline 
        Basic.             &  2 & 8 & 18  \\   \hline 
       Bern.               &  2 & 8 & 18  \\   \hline 
        \end{tabular} 
    \end{center} 
     \caption{\textit{X} mesh: Summary table of number of modes per systems.} \label{tab:number_of_modes_2D-meshX}
\end{table}%

\subsection{The fully discrete analysis} \label{sec:spectralAnalysis}
We analyze now  
the fully discrete  schemes obtained  using the RK, SSPRK and DeC time marching methods. 
Let us consider as an example the SSPRK schemes. 
If we define as $A:=\mass^{-1} (a_x\mathcal{K}_x+a_y\mathcal{K}_y+\delta  {\mathbb{S}})$ we can write the schemes as follows
\begin{equation}
\left  \{
    \begin{array}{ll}
    	\mathbf{U}^{(0)}:= & \mathbf{U}^n \\
        \mathbf{U}^{(s)} := & \sum_{j=0}^{s-1} \left( \gamma_{sj} \mathbf{U}^{(j)} + \Delta t \mu_{sj} A \mathbf{U}^{(j)} \right), \quad s \in \llbracket 1,S \rrbracket ,\\
	    \mathbf{U}^{n+1}:=&\mathbf{U}^{(S)}.
	\end{array}
    \right .
    \label{eq:discreteSSPRK}
\end{equation}

Expanding all the stages, we can obtain the following representation of the final stage:
\begin{equation}
\mathbf{U}^{n+1} = \mathbf{U}^{(0)} + \sum_{j=1}^{S} \nu_{j} \Delta t^jA^j  \mathbf{U}^{(0)} = \left(\mathcal{I} + \sum_{j=1}^{S} \nu_{j}\Delta t^j A^j  \right) \mathbf{U}^n, \label{RK_gamma0}
\end{equation}
where coefficients $\nu_j$ in \cref{RK_gamma0} are obtained as combination of coefficient $\gamma_{sj}$ and $\mu_{sj}$ in \cref{eq:discreteSSPRK} and $\mathcal{I} $ is the identity matrix. For example, coefficients of the fourth order of accuracy scheme RK4 are $\nu_1=1$, $\nu_2 = 1/2$, $\nu_3=1/6$ and $\nu_4 = 1/24$.\\ 

We can now compress the problem proceeding as in the time continuous case. In  particular, using Equations \eqref{eq_ad_fourier0} 
one easily shows that the problem can be written in terms of the local $d\times d$ matrices $\widetilde{A}:= \widetilde{\mass}^{-1}\left( a_x\widetilde{\mathcal{K}_x}+a_y\widetilde{\mathcal{K}_y}+\delta \widetilde{\mathbb{S}} \right)$ 
and in particular that
%
\begin{equation}
\widetilde{\mathbf{U}}^{n+1} = G \widetilde{\mathbf{U}}^{n}\quad\text{with}\quad
G:=  \left(\widetilde{\mathcal{I}} + \sum_{j=1}^{S} \nu_{j} \Delta t^j \widetilde{A}^j \right)  = e^{\epsilon \Delta t } e^{-i\omega \Delta t} ,
\label{eq:ampliG}
\end{equation}
where $G\in \R^{d\times d}$ is the amplification matrix depending on $\theta,\,\delta,\, \Delta t,\, \Delta x$ and $ \Delta y$. Considering each eigenvalue $\lambda_i$ of $G$, we can write the following formulae for the corresponding phase $\omega_i$ and damping coefficient $\epsilon_i$
\begin{align*}
\begin{cases}
   	e^{\epsilon_i \Delta t } \cos(\omega_i \Delta t)  = \text{Re}(\lambda_i) ,\\
   	- e^{\epsilon_i \Delta t } \sin(\omega_i \Delta t)  = \text{Im}(\lambda_i),
\end{cases} 
\Leftrightarrow \, \begin{cases} 
\omega_i\Delta t = \arctan \left( \frac{-\text{Im}(\lambda_i)}{\text{Re}(\lambda_i)} \right)  ,\\
(e^{\epsilon_i \Delta t })^2 = \text{Re}(\lambda)^2 + \text{Im}(\lambda)^2,
\end{cases} \Leftrightarrow \, \begin{cases} 
\dfrac{\omega_i}{k} = \arctan \left( \frac{-\text{Im}(\lambda_i)}{\text{Re}(\lambda_i)} \right)  \frac{1}{k \Delta t},\\
 \epsilon_i = \log\left( | \lambda_i | \right) \frac{1}{\Delta t}.
\end{cases}
\end{align*}
For the DeC method we can proceed with the same analysis transforming also the other involved matrices into their Fourier equivalent ones. 
Using \cref{eq:DeCasSSPRK} these terms would contribute to the construction of $G$ not only in the $\widetilde{A}$ matrix, but also in the coefficients $\nu_j$, which become matrices as well.
At the end we just study the final matrix $G$ and its eigenstructure, whatever process was needed to build it up.

The matrix $G$ describes one timestep evolution of the Fourier modes for all the $d$ different types of degrees of freedom. The damping coefficients $\epsilon_i$ tell if the modes are increasing or decreasing in amplitude and the phase coefficients $\omega_i$ describe the phases of such modes.

We remark that a necessary condition for stability of the scheme is that $ |\lambda_i | \leq 1$ or, equivalently, $\epsilon_i \leq 0 $ for all the eigenvalues. The goal of our study is to find the largest CFL number for which the stability condition is fulfilled and such that the dispersion error is \textit{not too large}.\\

For our analysis, we focus on the \textit{X} type triangular mesh in \cref{fig_Xmesh_P2} 
 with elements of degree $1$, $2$ and $3$. This \textit{X} type triangular mesh is also used in \cite{articleTaoLiu} for Fourier analysis of the acoustic wave propagation system.

\subsection{Methodology}
The methodology we explain in the following, will be applied to all the combination of schemes we presented above (in time: RK, SSPRK and DeC, discretisation in space: \textit{Basic}, \textit{Cubature} and \textit{Bernstein}, stabilization techniques: CIP, OSS and SUPG), in order to find the best coefficients (CFL, $\delta$), as in \cite{michel2021spectral}. \\

It must be remarked that the dispersion analysis must satisfy the Nyquist stability criterion, i.e., $\Delta x_{max} \leq \frac{L}{2}$ with $\Delta x_{max}$ the maximal distance between two nodes on edges. In other words, $k_{max} = \frac{2\pi}{L_{min}} = \frac{2\pi}{2 \Delta x_{max}}=\frac{\pi}{\Delta x_{max}}$. This tells us where $k$ should vary, i.e., $k \in ] 0,\pi / \Delta x_{max} ]$. \\

What we aim to do is an optimization process also on the stabilization parameter and the CFL number. 
With the notation of \cite{michel2021spectral}, we will set for the different stabilizations 
\begin{equation*}
\begin{split}
\quad \text{OSS :} \;\;& \tau_K =\delta  \Delta x  | a |,\\[5pt]
\quad\text{CIP :}   \;\;& \tau_f = \delta\Delta x^2  | a |,\\[5pt]
\quad \text{SUPG :}\;\;& \tau_K =\delta  \Delta x/|a|.
\end{split} 
\end{equation*}

One of our objectives is to explore the space of parameters (CFL,$\delta$), and to propose criteria allowing to set these parameters to provide the most stable, least dispersive and least expensive methods.
A clear and natural criterion is to exclude all parameter values for which there exists at least a wavenumber $\theta$ or an angle $\Phi \in [0,2\pi]$ such that we obtain an amplification of the mode, i.e., $\epsilon(\theta)>10^{-12}$ (taking into account the machine precision errors that might occur). Doing so, we obtain what we will denote as \textit{stable area} in $(\CFL,\theta)$ space. 
For all the other points we propose 3 strategies to minimize a combination of dispersion error and computational cost. \\

In the following we describe the strategy we adopt to find the best parameters couple (\CFL,$\delta$) that minimizes a global solution error, denoted by $\eta_u$,  while maximizing the {\CFL} in the stable area. In particular, we start from  the relative square error of $u$ 
\begin{align}
	\left| \frac{u(t)-u_{ex}(t)}{u_{ex}(t)}\right|^2= &\left|e^{\epsilon t - i t(\omega-\omega_{ex})}-1\right|^2\\
	=&\left[e^{\epsilon t}\cos(t(\omega-\omega_{ex}))-1\right]^2+\left[e^{\epsilon t}\sin(t(\omega-\omega_{ex}))\right]^2\\
	=&e^{2\epsilon t} - 2 e^{\epsilon t} \cos (t(\omega-\omega_{ex})) +1.
\end{align}
Here, we denote with $\epsilon$ and $\omega$ the damping and phase of the \textit{principal }mode and with $\omega_{ex}=\mathbf{k} \cdot \mathbf{a}$ the exact phase.
For a small enough  dispersion error $|\omega-\omega_{ex} |\ll 1$, we can expand the cosine in the previous formula in a truncated Taylor series as 
\begin{align}
	\left|\frac{u(t)-u_{ex}(t)}{u_{ex}(t)}\right|^2\approx&\underbrace{\left[e^{\epsilon t} -1\right]^2}_{\text{Damping error}} + \underbrace{e^{\epsilon t}t^2 \left[\omega- \omega_{ex}\right]^2}_{\text{Dispersion error}}.
\end{align}
We then compute an   error at the final time $T=1$, over the whole phase domain, using   at least 3 points per wave $0\leq k \Delta x_p \leq \frac{2\pi}{3}$, with $\Delta x_p=\frac{\Delta x}{p}$, and   $p$ the degree of the polynomials. 
We  obtain the following $\mathbb{L}_2$ error definition, 
\begin{equation}\label{eq:dispersionError}
	\eta_u(\omega,\epsilon)^2:= \frac{3}{2\pi} \left[\int_{0}^{\frac{2\pi}{3}} (e^{\epsilon}-1 )^2 dk + \int_{0}^{\frac{2\pi}{3}} e^\epsilon(\omega-\omega_{ex})^2 dk \right].
\end{equation}
Recalling that $\epsilon=\epsilon(k\Delta x,\CFL,\delta, \Phi)$ and $\omega=\omega(k,\Delta x,\CFL,\delta, \Phi)$, we need to further set the parameter $\Delta x_p$. We choose it to be large $\Delta x_p=1$, with the hope that for finer grids the error will be smaller. 
Moreover, we need to check that the stability condition holds for all the possible angles $\Phi\in [0,2\pi]$.

Finally,  we seek for the couple $(\text{CFL}^*,\delta^*)$ such that
\begin{equation}
\label{cfl_d}
	(\text{CFL}^*,\delta^*)=\arg \max_{\text{CFL}} \left\lbrace \eta (\omega, \epsilon, \Phi')< \mu  \min_{\text{stable } (\text{CFL},\delta)}  \max_{\Phi} \eta(\omega,\epsilon, \Phi), \quad \forall\,  \Phi' \in [0,2\pi] \,\right\rbrace,
\end{equation}
where the dependence on $\Phi$ of $\eta$ is highlighted with an abuse of notation.
For this strategy, the parameter $\mu$ must be chosen in order to balance the  requirements on stability and accuracy. 
After having tried different values,  we  have  set $\mu$ to $10$ providing a sufficient flexibility to obtain results of practical usefulness. Indeed, the found values will be tested in the numerical section.

\begin{figure}
	\centering
	\subfigure[$\Phi =0$]{\includegraphics[width=0.49\linewidth,trim={100 0 100 0}, clip]{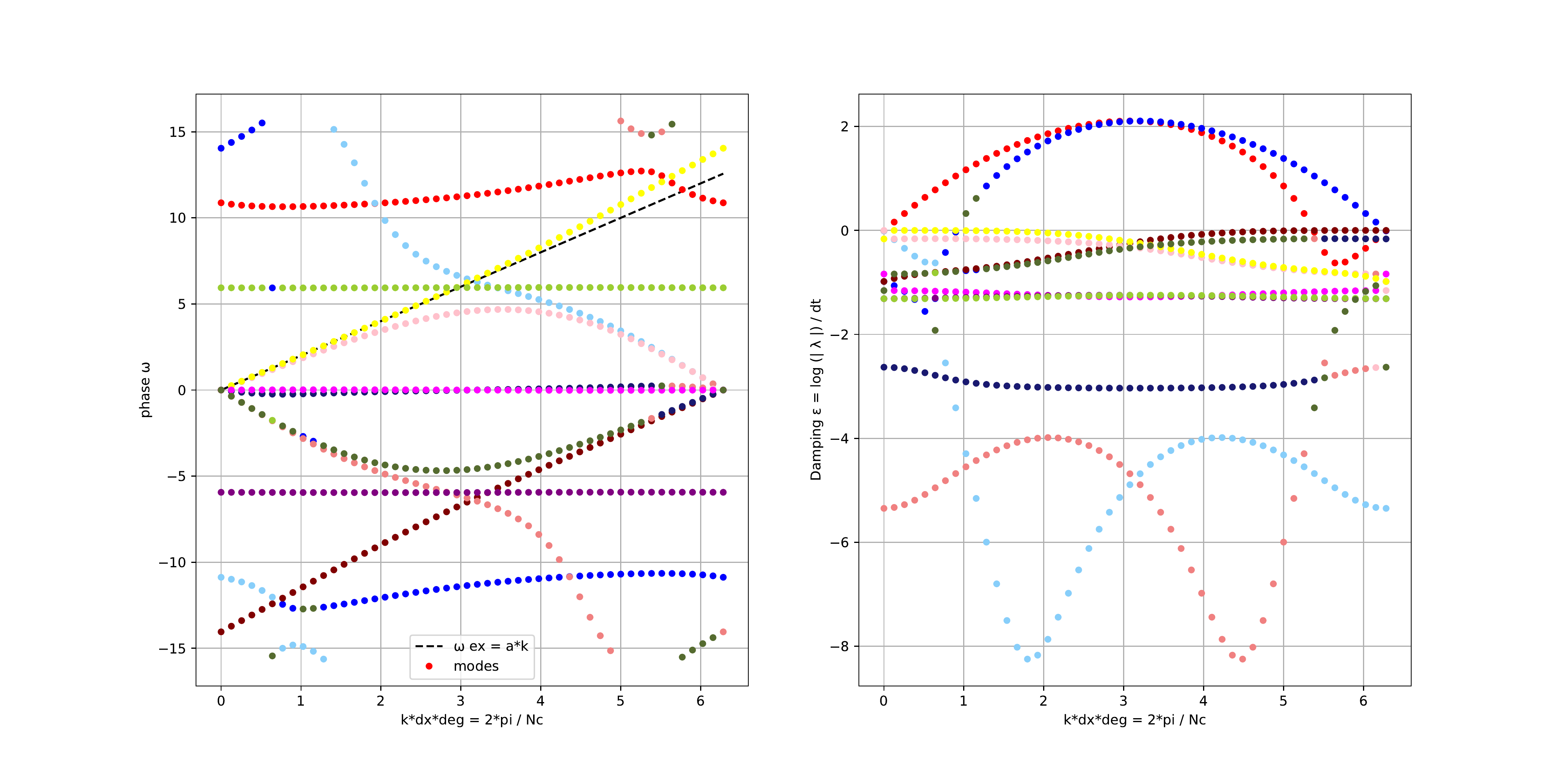}}
	\subfigure[$\Phi=3\pi / 16$]{\includegraphics[width=0.49\linewidth,trim={100 0 100 0}, clip]{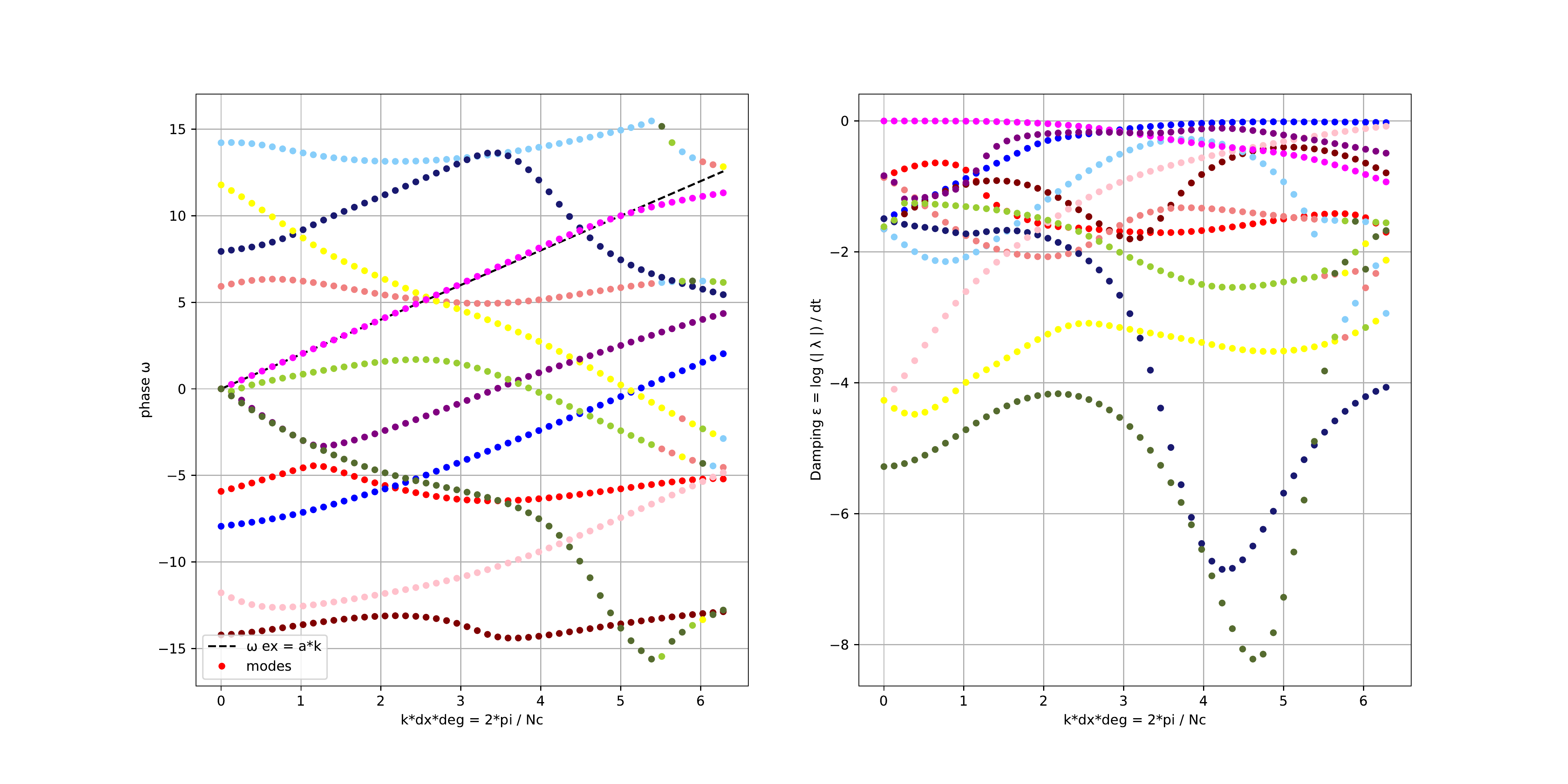}}
	\caption{Comparison of dispersion curves $\omega_i$ and damping coefficients $\epsilon_i$,  for \textit{Cubature} $\TP_2$ elements, with SSPRK time discretization and OSS stabilization. $\Phi =0$ at the left and $\Phi=3\pi / 16$ at the right.} 		\label{fig:fourier_two_angles}
\end{figure}

To show the influence of the angle $\Phi$ on the optimization problem we show an example for the \textit{X} mesh. For a given couple of parameters (\CFL,$\delta$) = $(0.4,0.01)$ we compare the results for $\Phi =0$ and $\Phi=3\pi/16$. In \cref{fig:fourier_two_angles} we compare the phases $\omega_i$ and the damping coefficients $\epsilon_i$ for the two angles. It is clear that for the angle $\Phi=0$, on the left, there are some modes which are not stable $\epsilon_i>0$, while for $\Phi=3\pi/16$ all modes are stable.

\begin{figure}
	\centering
	\subfigure[$\Phi =0$]{\includegraphics[width=0.45\linewidth,trim={500 0 460 0}, clip]{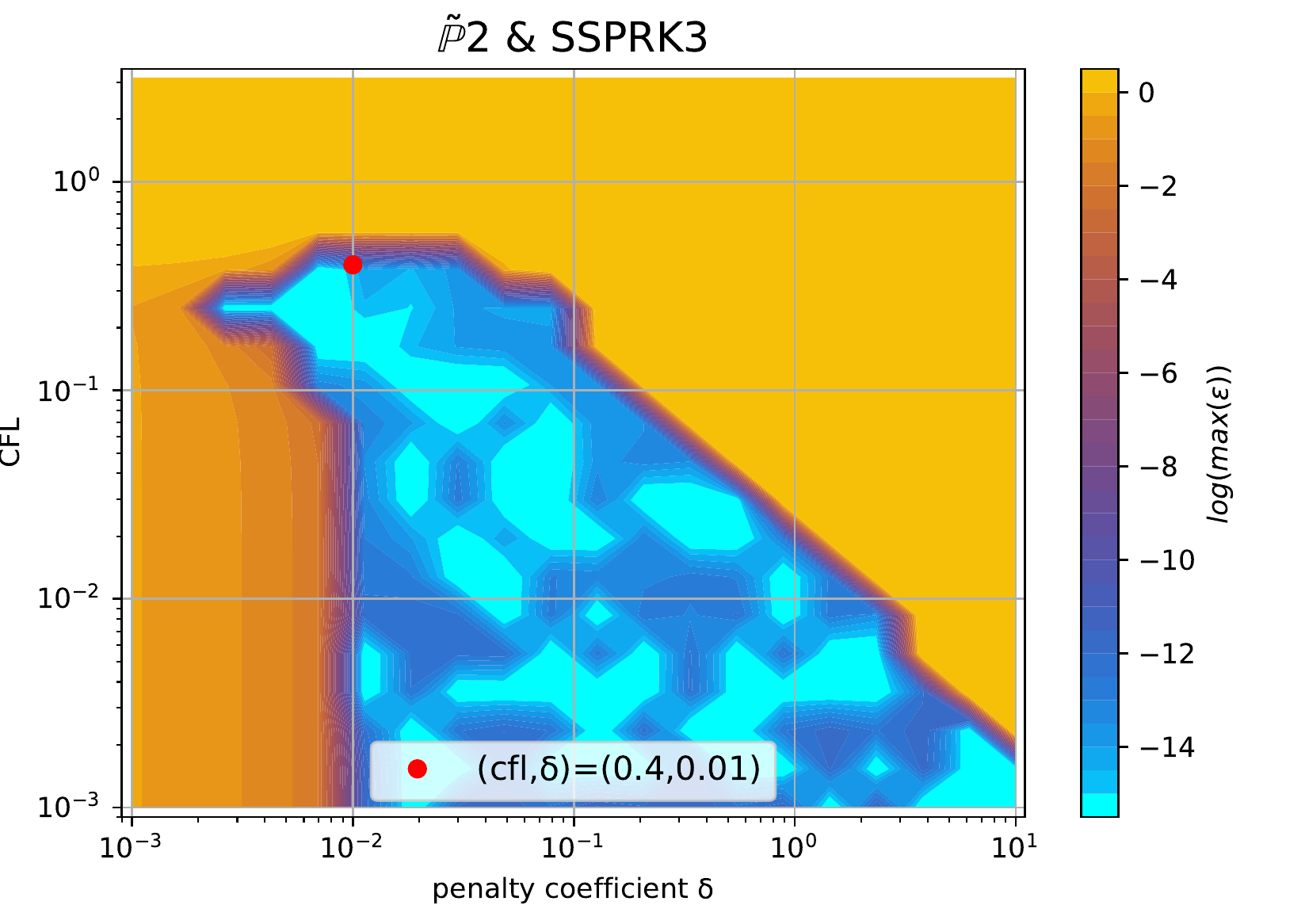}}
	\subfigure[$\forall \Phi \in {[} 0,2 \pi {]} $]{\includegraphics[width=0.45\linewidth,trim={500 0 460 0}, clip]{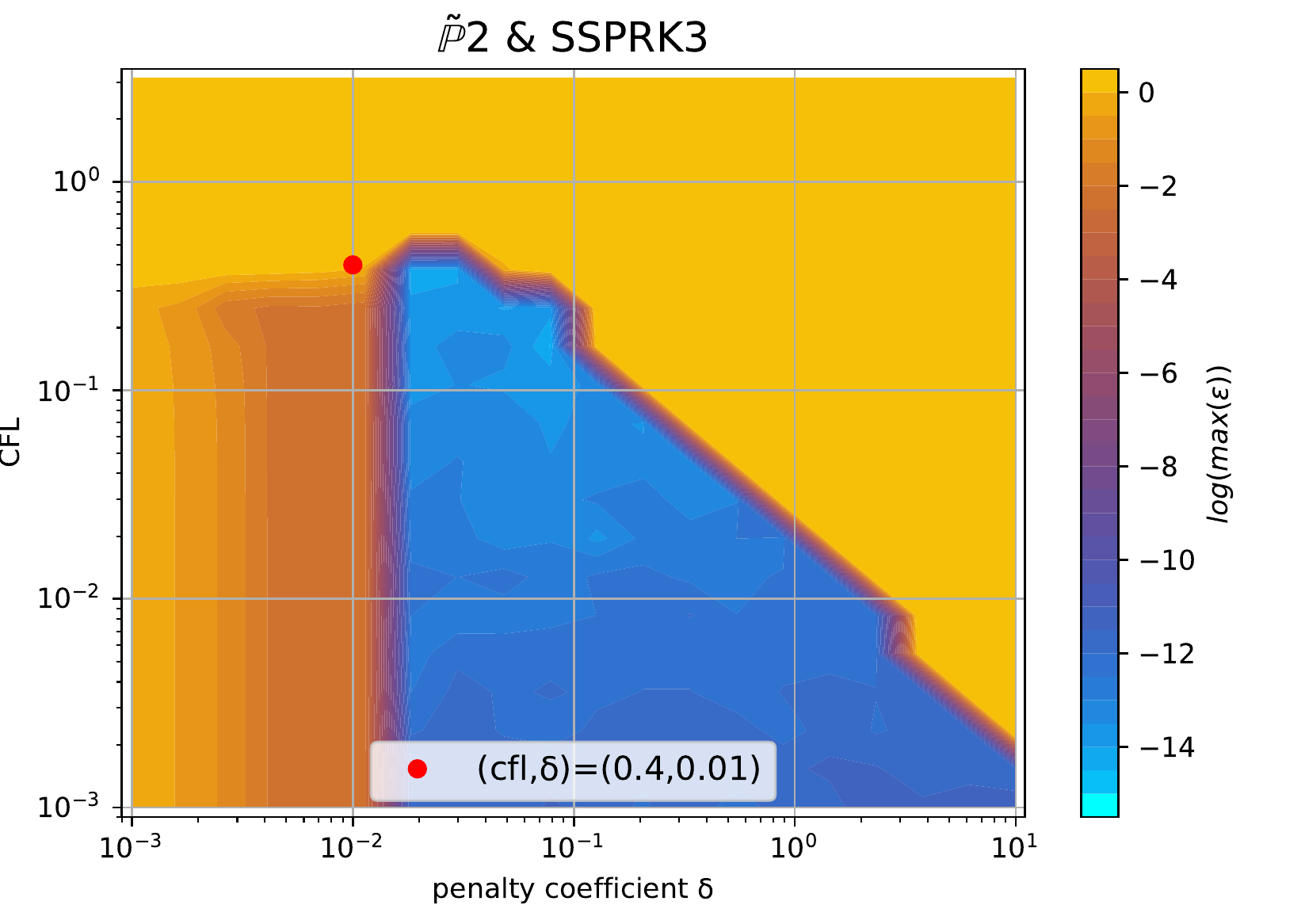}}
	\caption{Plot of $\log(\max_i \epsilon_i)$ for \textit{Cubature} $\TP_2$ elements, SSPRK time discretization and OSS stabilization. The blue and light blue region is the stable one. At the left only for $\Phi =3\pi/16$, at the right we plot the maximum over all $\Phi$. }\label{fig:fourier_ex_comparison_phi}
\end{figure}

The angle can widely influence the whole analysis as  one can observe in the plot of $\max_i  \epsilon_i$ in \cref{fig:fourier_ex_comparison_phi}, where we observe that for the only angle $\Phi = 3\pi/16$ we would obtain an optimal parameter in (\CFL,$\delta$) = $(0.4,0.01)$, while, using all angles, this value is not stable anymore.
%

\begin{remark}
To define the stable region, we should only consider configurations for which the damping is below machine accuracy.
In practice, this cannot be done due to the fact that the eigenvalue problem  arising from  \eqref{eq:ampliG}  is only solved
approximately using  the linear algebra package of \texttt{numpy}. This introduces some uncertainty in the definition
of the stability region as   machine accuracy needs  to be replaced
by some other finite threshold.
\end{remark} 

\subsection{Results of the fourier analysis using the \textit{X} type mesh} \label{sec:fourier_globalresults}
In this section, we illustrate the result obtained with the methodology explained above. 
For clarity not all the results are reported in this work, however we place all the plots for all possible combination of schemes in an online repository \cite{michel2021spectral}. We will provide some examples here and a summary of the main results that we obtained.

 \begin{figure}
	\centering
	\subfigure[SUPG]{\includegraphics[width=0.3\linewidth,trim={1000 0 0 0}, clip]{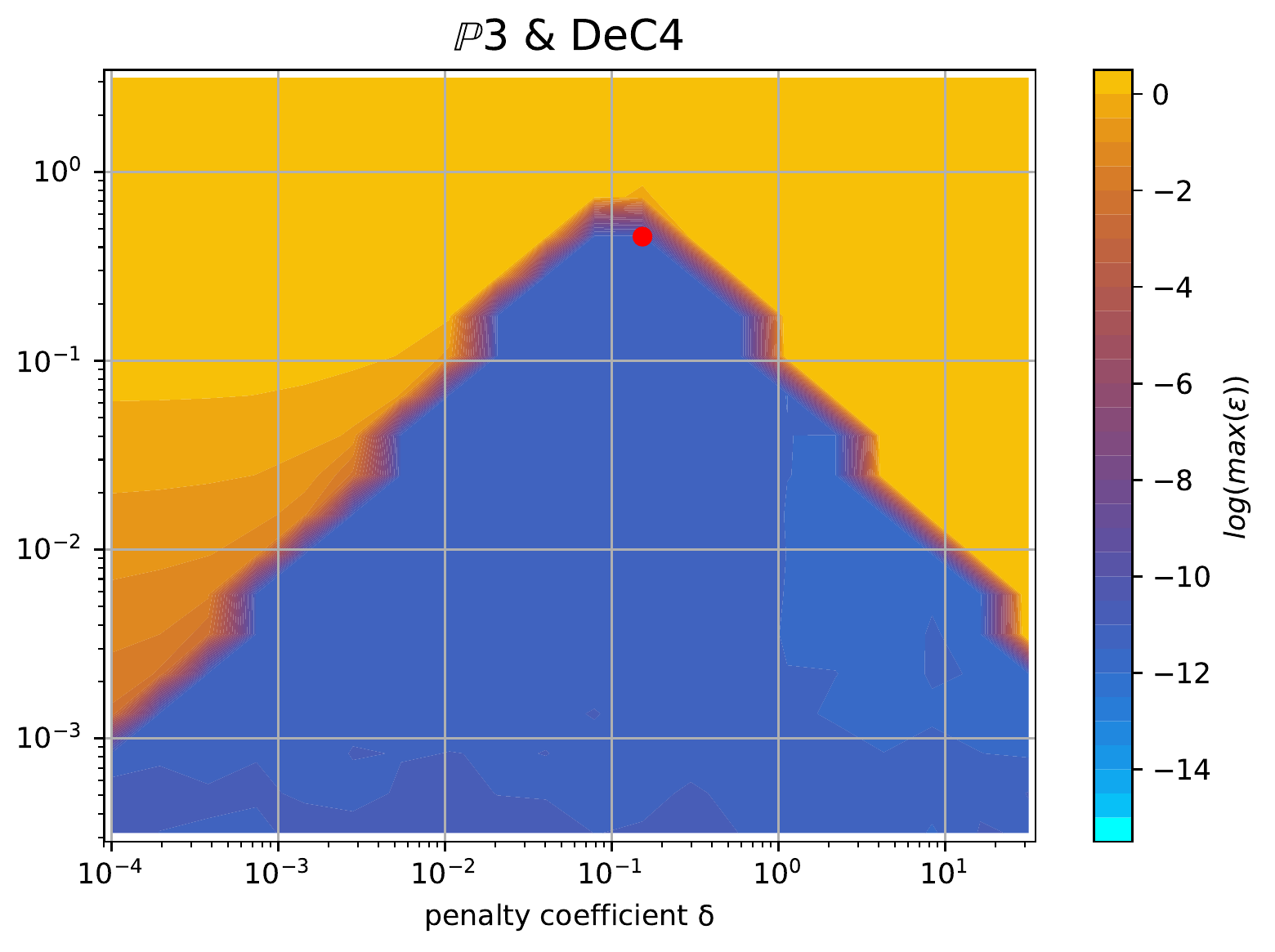}}
	\subfigure[OSS]{\includegraphics[width=0.3\linewidth,trim={1000 0 0 0}, clip]{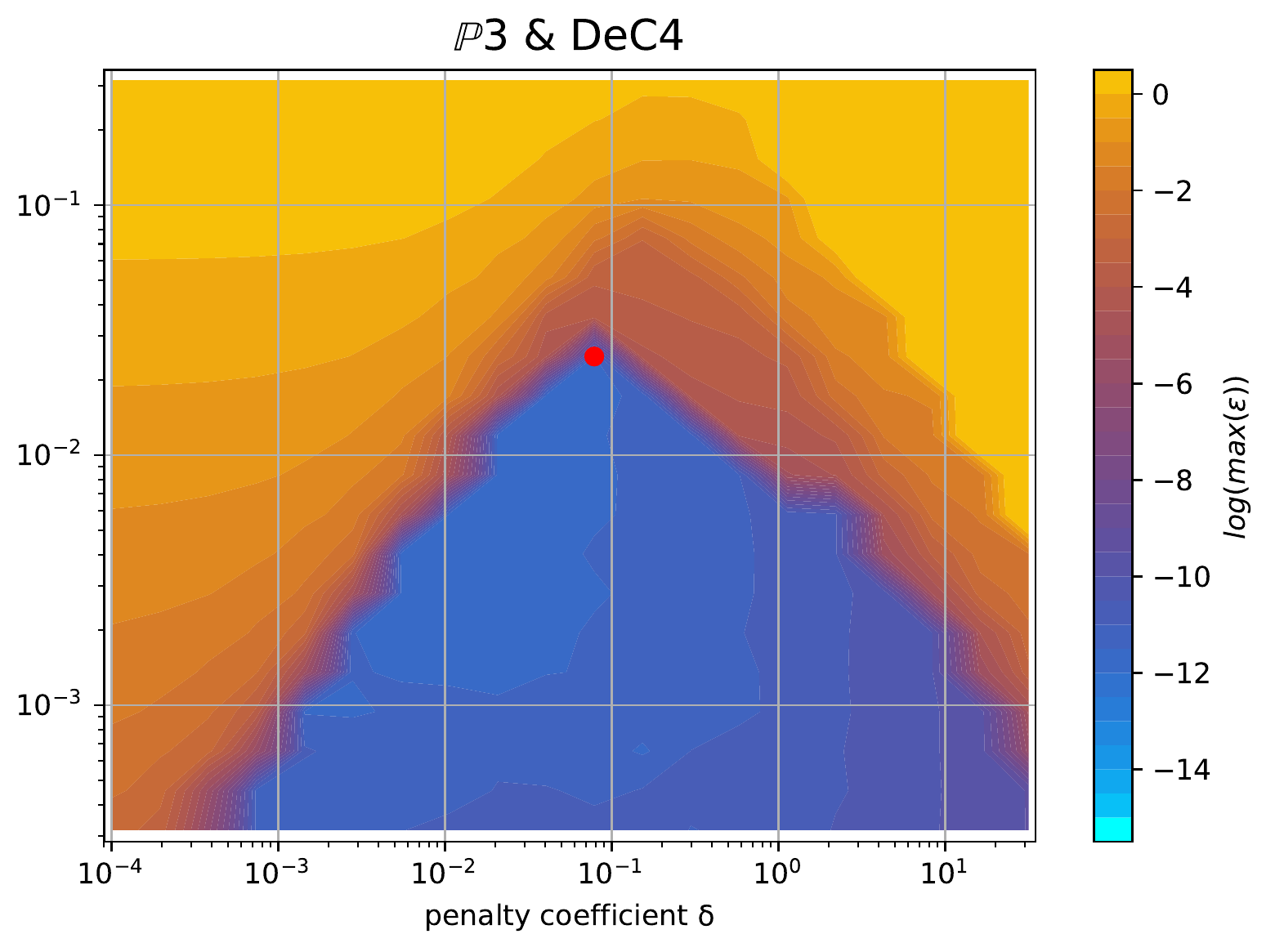}}
	\subfigure[CIP]{\includegraphics[width=0.3\linewidth,trim={1000 0 0 0}, clip]{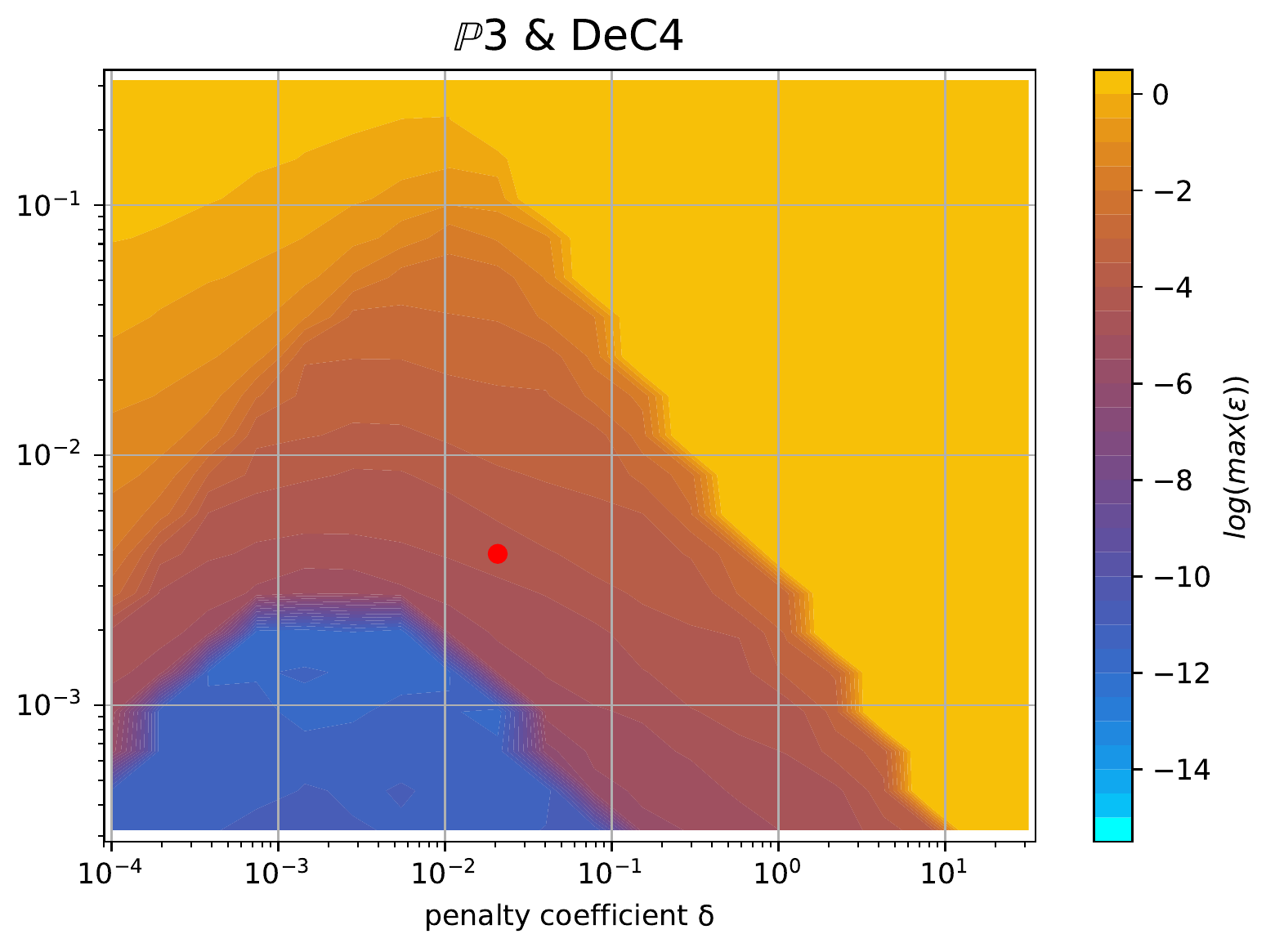}}
	\caption{Damping coefficients $\log(\max_i \epsilon_i)$ for $\B_3$ \textit{Bernstein} elements and the DeC method with, from left to right, SUPG, OSS and CIP stabilization. The red dot is the optimum according to \eqref{cfl_d}.} \label{fig:particularfourier_all}
\end{figure}

The first type of plot we introduce is helping us in understanding how we can define the stability region in the $(\CFL,\delta)$ plane. 
So, for every $(\CFL,\delta)$ we plot the maximum of $\log(\epsilon_i)$ over all modes and angles $\Phi \in [0,2\pi]$ (thanks to the symmetry of the mesh we can reduce this interval). An example is given in the right plot of \cref{fig:fourier_ex_comparison_phi}, 
it is clear that the whole blue area is stable and the yellow/orange area is unstable. In other cases, this boundary is not so clear and setting a threshold to determine the stable area can be challenging. 
In \cref{fig:particularfourier_all} we compare different stabilizations for DeC with $\B_3$ elements. 
In the CIP stabilization case, we clearly see that there is no clear discontinuity between unstable values and stable ones, as in SUPG, because there is a transient region where $\max_i \epsilon_i$ varies between $10^{-7}$ and $10^{-4}$.

\begin{figure}
	\centering
	\includegraphics[width=\linewidth]{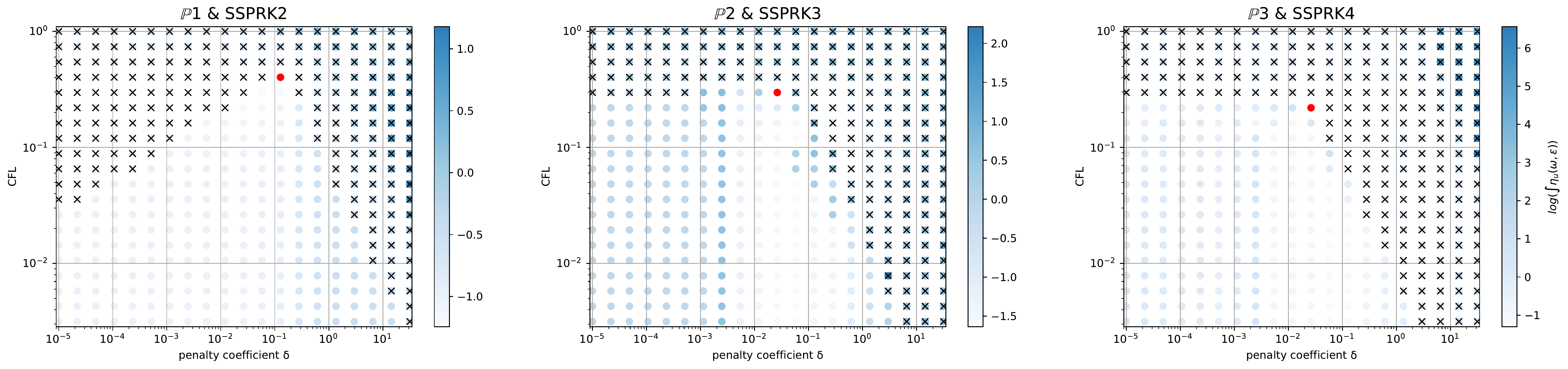}
	\caption{$\log(\eta_u)$ values (blue scale) and stable area (unstable with black crosses), on $(\CFL,\delta)$  plane. The red dot denotes the optimal value. From left to right $\P_1$, $\P_2$, $\P_3$ \textit{Basic} elements with SSPRK scheme and OSS stabilization}\label{fig:fourier_global_basic_SSPRK_OSSX}
	\centering
	\includegraphics[width=\linewidth]{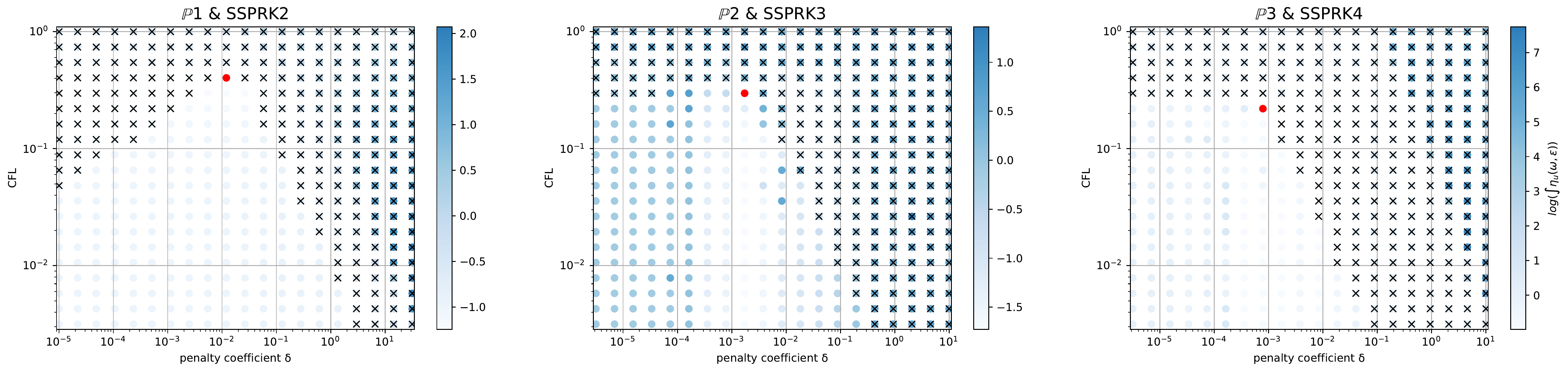}
	\caption{$\log(\eta_u)$ values (blue scale) and stable area (unstable with black crosses), on $(\CFL,\delta)$  plane. The red dot denotes the optimal value. From left to right $\P_1$, $\P_2$, $\P_3$ \textit{Basic} elements with SSPRK scheme and CIP stabilization }\label{fig:fourier_global_basic_SSPRK_CIPX}
\end{figure}

\begin{figure}
	\centering
	\includegraphics[width=\linewidth]{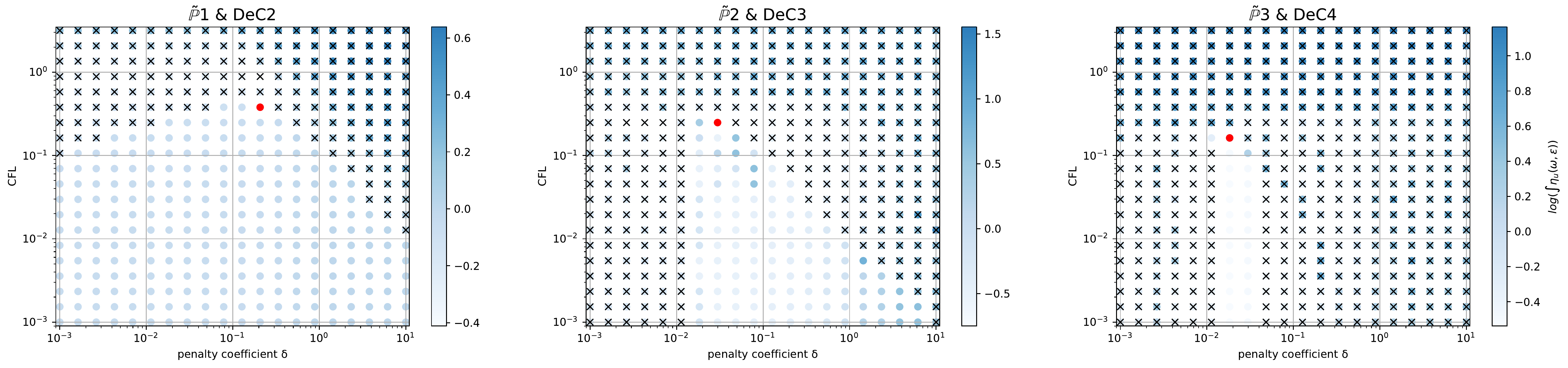}
	\caption{$\log(\eta_u)$ values (blue scale) and stable area (unstable with black crosses), on $(\CFL,\delta)$  plane. The red dot denotes the optimal value. From left to right $\TP_1$, $\TP_2$, $\TP_3$ \textit{Cubature} elements with DeC scheme and OSS stabilization}\label{fig:fourier_global_cub_DeC_OSSX}
	\centering
	\includegraphics[width=\linewidth]{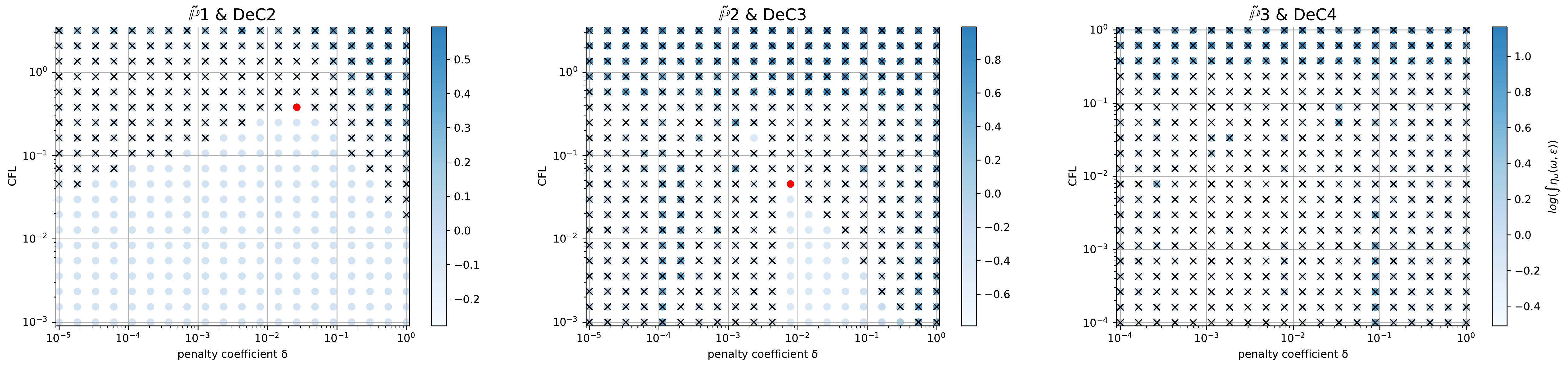}
	\caption{$\log(\eta_u)$ values (blue scale) and stable area (unstable with black crosses), on $(\CFL,\delta)$  plane. The red dot denotes the optimal value. From left to right $\TP_1$, $\TP_2$, $\TP_3$ \textit{Cubature} elements with DeC scheme and CIP stabilization}\label{fig:fourier_global_cub_DeC_CIPX}
\end{figure}

The second type of plot combines the chosen stability region with the error $\eta_u$. We plot on the $(\CFL,\delta)$ plane some black crosses on the unstable region, where there exists an $i$ and $\Phi$ such that $\epsilon_i > 10^{-7}$. The color represents $\log(\eta_u)$ and the best value according to the previously described method is marked with a red dot.
In	\cref{fig:fourier_global_basic_SSPRK_OSSX,fig:fourier_global_basic_SSPRK_CIPX,fig:fourier_global_cub_DeC_OSSX,fig:fourier_global_cub_DeC_CIPX} we show some examples of these plots for some schemes, for different $p=1,2,3$. In \cref{fig:fourier_global_basic_SSPRK_OSSX,fig:fourier_global_basic_SSPRK_CIPX} we test the \textit{Basic} elements with the SSPRK time discretization, while in \cref{fig:fourier_global_cub_DeC_OSSX,fig:fourier_global_cub_DeC_CIPX} we use the \textit{Cubature} elements with DeC time discretization. We compare also different stabilization technique: in \cref{fig:fourier_global_basic_SSPRK_OSSX,fig:fourier_global_cub_DeC_OSSX} we use the OSS, while in \cref{fig:fourier_global_basic_SSPRK_CIPX,fig:fourier_global_cub_DeC_CIPX} the CIP.
One can observe many differences among the schemes. For instance, for $p=3$ we see a much wider stable area for SSPRK than with DeC and, in the \textit{Cubature} DeC case, we see that the CIP requires a reduction in the CFL number with respect to the OSS stabilization.

%

We summarize the results obtained by the optimization strategy in \cref{tab:dispersion_cfl_LinearAdvection_meshX-RES} for all the combinations of spatial, time and stabilization discretization. The CFL and $\delta$ presented there are optimal values obtained by the process above described, which we aim to use in simulations to obtain stable and efficient schemes. Unfortunately, as already mentioned above, for some schemes the stability area is not so well defined for several reasons. One of these reasons is the \textit{"shape"} of the stability area as for one-dimensional problems, see  \cite{michel2021spectral}. Other issues that affect this analysis are the numerical precision, see Section \ref{sec:fourier_rocket}, and the mesh configuration, see Section \ref{sec:fourier_mesh_pattern}. In the following we study more in details these cases and how one can find better values.


\begin{table}[H] 
\small  
 \begin{center} 
		\begin{tabular}{| c | c || c | c | c | }  
	     \hline 
	     \multicolumn{2}{|c||}{Element $\&$ }  & \multicolumn{3}{|c|}{SUPG}  \\ \hline 
	     \multicolumn{2}{|c||}{ Time scheme }  & $\mathbb{P}_1$ & $\mathbb{P}_2$ & $\mathbb{P}_3$  \\ \hline \hline 
        \parbox[t]{2mm}{\multirow{2}{*}{\rotatebox[origin=c]{90}{Basic}}}              &  SSPRK  &  0.739 (0.127)   &  0.298 (0.058)   &  0.22 (0.026)  \\ 
               &  RK  &  0.403 (0.127)   &  0.298 (0.026)   &  0.22 (5.46e-03)  \\ 
         \hline 
       \parbox[t]{2mm}{\multirow{3}{*}{\rotatebox[origin=c]{90}{\centering  Cub.}}}              &  DeC  &  0.616 (0.28)   &  0.234 (0.04)$^*$   &  0.144 (0.04)  \\ 
               &  SSPRK  &  1.062 (0.28)   &  0.379 (0.021)$^*$   &  0.234 (0.011)$^*$ \\ 
               &  RK  &  0.616 (0.28)   &  0.234 (0.04)   &  0.144 (0.04)  \\ 
         \hline 
       \parbox[t]{2mm}{\multirow{3}{*}{\rotatebox[origin=c]{90}{\centering  Bern.}}}              &  DeC  &  0.739 (0.298)   &  0.455 (0.298)$^*$   &  0.455 (0.153)$^*$  \\ 
               &  SSPRK  &  0.739 (0.127)   &  0.298 (0.058)   &  0.22 (0.026)  \\ 
               &  RK  &  0.403 (0.127)   &  0.298 (0.026)   &  0.22 (5.46e-03)  \\ 
         \hline 
        \end{tabular} 
    \end{center} 
 \begin{center} 
		\begin{tabular}{| c | c || c | c | c || c | c | c | }  
	     \hline 
	     \multicolumn{2}{|c||}{Element $\&$ }  & \multicolumn{3}{|c||}{OSS}  & \multicolumn{3}{|c|}{CIP}  \\ \hline 
	     \multicolumn{2}{|c||}{ Time scheme }  & $\mathbb{P}_1$ & $\mathbb{P}_2$ & $\mathbb{P}_3$   & $\mathbb{P}_1$ & $\mathbb{P}_2$ & $\mathbb{P}_3$  \\ \hline \hline 
        \parbox[t]{2mm}{\multirow{2}{*}{\rotatebox[origin=c]{90}{Basic}}}              &  SSPRK  &  0.403 (0.127)   &  0.298 (0.026)   &  0.22 (0.026)   &  0.403 (0.012)   &  0.298 (1.73e-03)   &  0.22 (7.85e-04)$^*$  \\ 
               &  RK  &  0.22 (0.058)   &  0.22 (0.026)   &  0.22 (0.012)   &  0.298 (0.012)   &  0.22 (1.73e-03)   &  0.22 (3.57e-04)   \\ 
         \hline 
       \parbox[t]{2mm}{\multirow{3}{*}{\rotatebox[origin=c]{90}{\centering  Cub.}}}              &  DeC  &  0.379 (0.207)   &  0.248 (0.03)   &  0.162 (0.018)   &  0.379 (0.026)   &  0.045 (7.85e-03)$^*$   &   /   \\ 
               &  SSPRK  &  0.58 (0.336)   &  0.379 (0.03)   &  0.248 (0.018)   &  0.58 (0.048)   &  0.07 (7.85e-03)$^*$   &   /   \\ 
               &  RK  &  0.379 (0.207)   &  0.248 (0.03)   &  0.162 (0.018)   &  0.379 (0.026)   &  0.045 (7.85e-03)   &   /   \\ 
         \hline 
       \parbox[t]{2mm}{\multirow{3}{*}{\rotatebox[origin=c]{90}{\centering  Bern.}}}              &  DeC  &  0.173 (0.58)   &  0.036 (0.298)  &  0.025 (0.078)$^*$   &  0.173 (0.153)   &  0.012 (0.021)   &  0.004 (0.021)$^*$    \\ 
               &  SSPRK  &  0.403 (0.127)   &  0.298 (0.026)   &  0.22 (0.026)   &  0.403 (0.012)   &  0.298 (1.73e-03)   &  0.22 (7.85e-04)   \\ 
               &  RK  &  0.22 (0.058)   &  0.22 (0.026)   &  0.22 (0.012)   &  0.298 (0.012)   &  0.22 (1.73e-03)   &  0.22 (3.57e-04)   \\ 
         \hline 
        \end{tabular} 
    \end{center} 
     \caption{\textit{X} mesh: Optimized CFL and penalty coefficient $\delta $ in parenthesis, minimizing $\eta_u$. \newline
     \textit{"/"} means that the fourier analysis shown that the scheme is unstable. \newline
     $^*$ These values are not reliable, see \cref{sec:fourier_rocket}.} 
     \label{tab:dispersion_cfl_LinearAdvection_meshX-RES}
 \end{table}%

\subsection{Comparison with a space-time split stability analysis} \label{sec:fourier_rocket}
\begin{figure}
	\centering
	\includegraphics[width=0.5\linewidth,trim={1000 0 0 0}, clip]{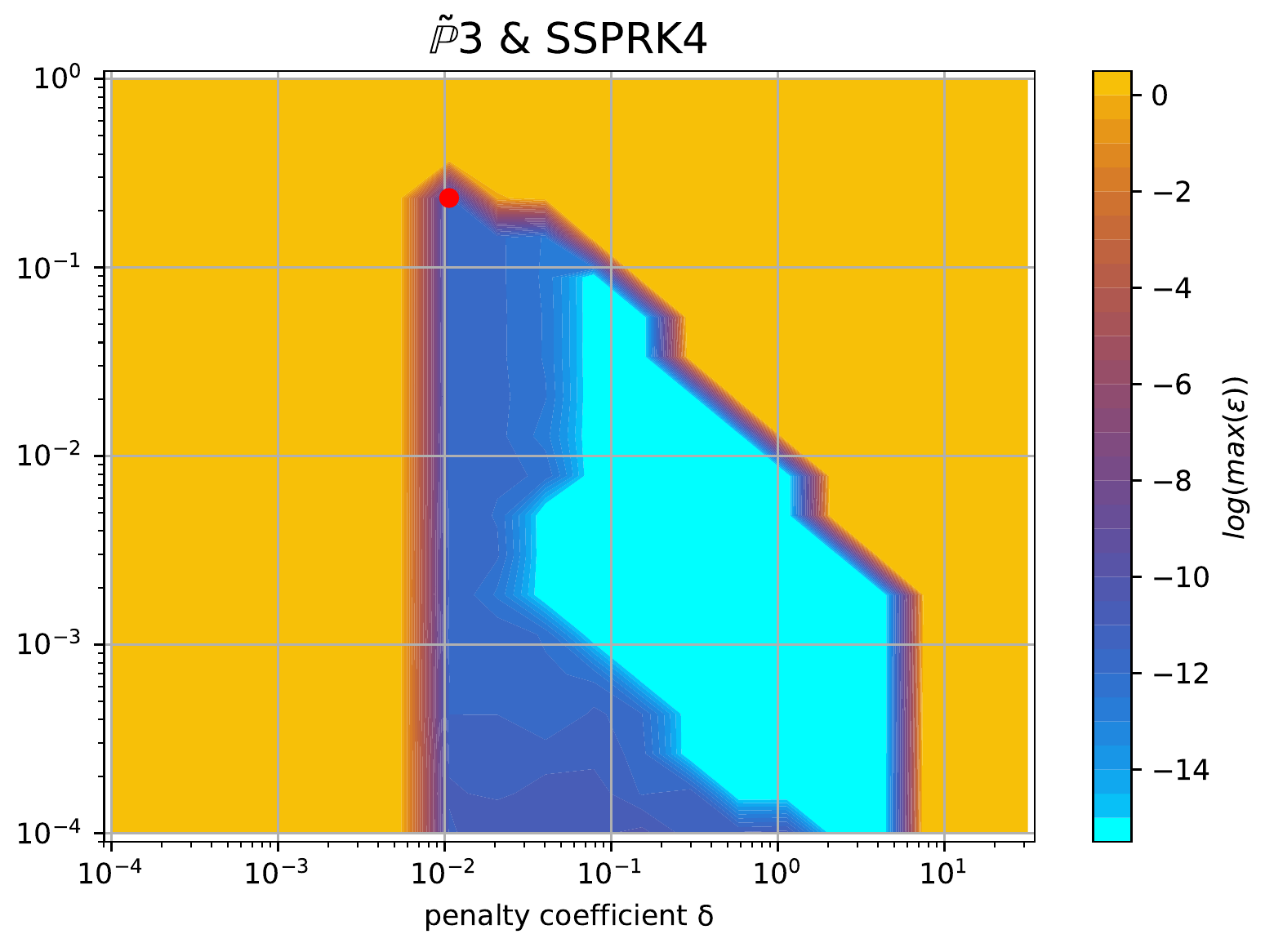}
	\caption{Logarithm of the amplification coefficient $\log(\max_i (\varepsilon_i))$ for SUPG stabilization with $\TP_3$ \textit{Cubature} elements and the SSPRK method. Unstable region in yellow, the red dot is the optimal parameter according to \eqref{cfl_d}} \label{fig:particularfourier_cub_supg_ssp4}
\end{figure}
In this section, we show another stability analysis to slightly improve the results obtained above. Indeed, the solution of the eigenvalue problem \eqref{eq:ampliG} is only obtained within some   approximation from the \texttt{numpy} numerical library. In some cases,
the threshold used to define the stability region is defined in a somewhat heuristic manner.  So to confirm the results, we use
independently another criterion. To this end we 
 treat independently  the temporal and spatial 
 discretizations as in the
  method of lines. We then study only the spectral properties of the spatial discretization alone, computing the eigenvalues of the corresponding matrix $A$ (cf. \eqref{eq:discreteSSPRK}). With this information, we then
    check whether they belong to the stability area of the time discretization.
%

\begin{figure}
\centering
\subfigure[(\CFL,$\delta$)=$(0.234 ,0.011)$ \label{fig:unstableSSPRKSUPG_rocket}]{\includegraphics[width=0.4\linewidth,trim={0 0 0 40}, clip]{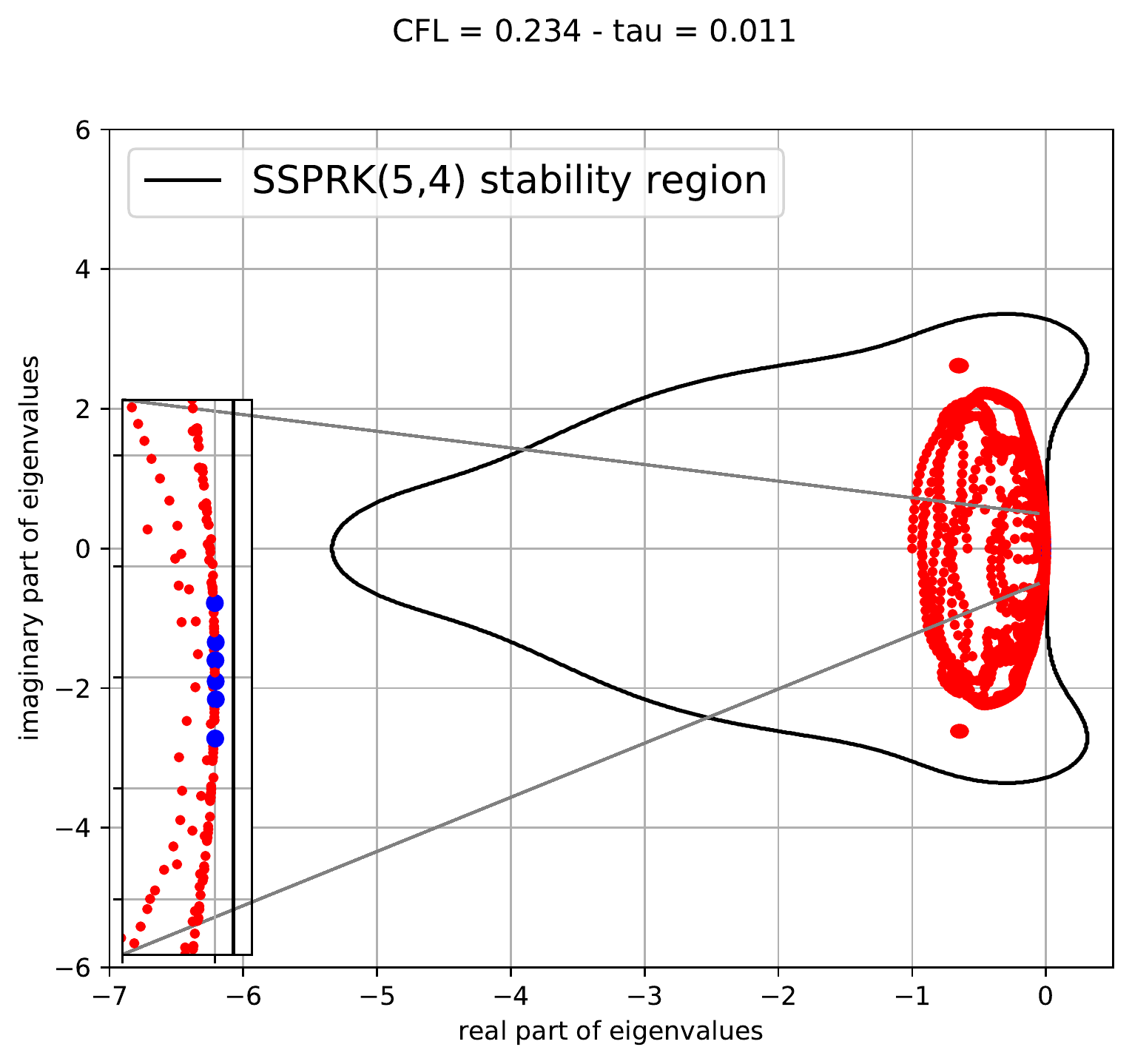}
}\;\;\;\;\;
\subfigure[(\CFL,$\delta$)=$(0.18 ,0.04)$ \label{fig:stableSSPRKSUPG_rocket}]{\includegraphics[width=0.4\linewidth,trim={0 0 0 40}, clip]{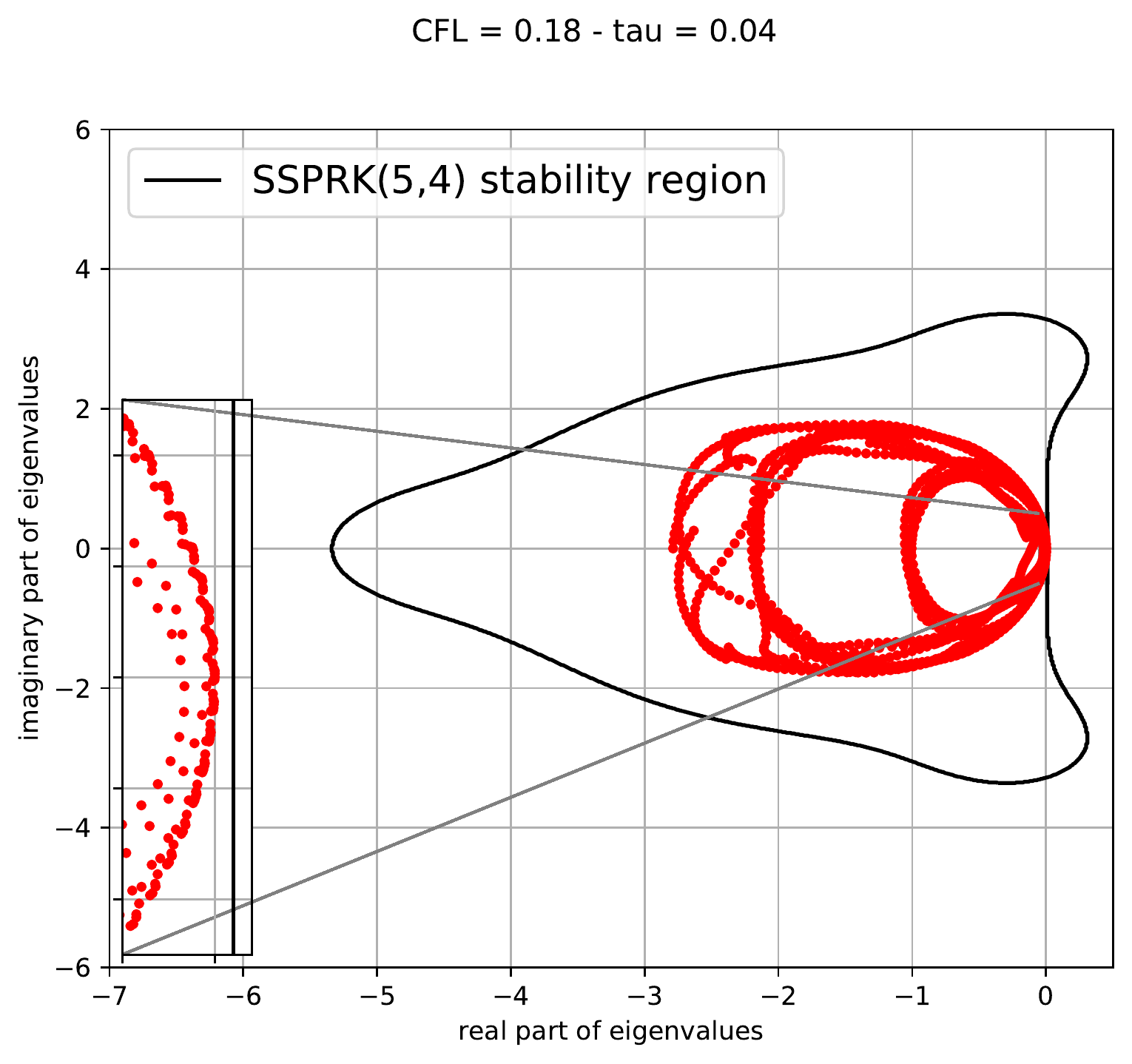}}
\caption{Eigenvalues of $\widetilde{A}$ using cubature discretization and the SUPG stabilization (varying $k$) and stability area of the SSPRK method. In red the stable eigenvalues, in blue the unstable ones.}
\label{fig:rocket_meshX_cub_SSPRK_SUPG_P3_unstable}
\end{figure}

In particular, following \cite{book_stab_analysis_ch2}, we write the time discretization for Dahlquist's equation
\begin{equation}
\partial_t u - \lambda u = 0 \label{eq:stab_rocket_eq},
\end{equation}
in this example, we consider the SSPRK discretization \eqref{eq:discreteSSPRK}.
From \cref{RK_gamma0} we can write the amplification coefficient $\Gamma(\lambda)$, i.e.,
\begin{equation}
\mathbf{U}^{n+1} = \mathbf{U}^{(0)} + \sum_{j=1}^{S} \nu_{j} \Delta t^j\lambda^j  \mathbf{U}^{(0)} = \underbrace{ \left(\mathcal{I} + \sum_{j=1}^{S} \nu_{j}\Delta t^j \lambda^j  \right)}_{\Gamma(\lambda)} \mathbf{U}^n. \label{RK_gamma02}
\end{equation}
The stability condition for this SSPRK scheme is given by $\Gamma(\lambda) \leq 1$.
Now, when we substitute the Fourier transform of the spatial semidiscretization $\widetilde{A}$ to the coefficient $\lambda$ and we diagonalize the system (or we put it in Jordan's form), we obtain a condition on the eigenvalues of $\widetilde{A}$.
Then, using the parameters provided by the previous analysis (\CFL,$\delta$)=$(0.234 ,0.011)$, in Table~\ref{tab:dispersion_cfl_LinearAdvection_meshX-RES}, we plot the eigenvalues of $\widetilde{A}$ and the stability region of the SSPRK scheme for different $\theta \in [ 0, \pi ]$. We notice that for some values of $\theta$ some of the eigenvalues
fall slightly outside the
stable area, see Figure~\ref{fig:unstableSSPRKSUPG_rocket}. There are, indeed, few eigenvalues  dangerously close to the imaginary axis and some of them have actually positive real part (blue dots).
As suggested before, if we decrease the CFL and increase $\delta$, we move towards a safer region, so considering (\CFL,$\delta$)=$(0.18 ,0.04)$ with the same $\theta$, we obtain all stable eigenvalues, as shown in Figure~\ref{fig:stableSSPRKSUPG_rocket}.


The summary of the optimal parameters of Table~\ref{tab:dispersion_cfl_LinearAdvection_meshX-RES} updated taking into account also a larger safety region in the (CFL, $\delta$) plane (as explained in this section) can be found in Table~\ref{tab:new_param_meshX_LinearAdvection-2D-RES} in Appendix~\ref{sec:app_fourier_param}.

\subsection{Different mesh patterns} \label{sec:fourier_mesh_pattern}
Another important aspect about this stability analysis is the influence of the mesh structure on the results. As an example, we use the \textit{T}-mesh, another regular and structured mesh type depicted in Figure~\ref{fig_Tmesh_P2}. In Figure~\ref{fig_Tmesh_P2} we plot also the degrees of freedom for elements of degree 2 and the periodic elementary unit that we take into consideration for the Fourier analysis. The number of modes in the periodic unit for this mesh type are summarized in Table~\ref{tab:number_of_modes_2D-meshT}.  The elements of degree 3 can be found in Figure~\ref{fig:meshDofsP3} in Appendix~\ref{sec:app_mesh4fourier}.
\begin{table}
	\small  
	\begin{center} 
		\begin{tabular}{| c || c | c | c | }  
			\hline 
			Element  & $\mathbb{P}_1$ & $\mathbb{P}_2$ & $\mathbb{P}_3$    \\ \hline \hline 
			Cub.              &  1 & 6 & 13  \\   \hline 
			Basic.             &  1 & 4 & 9  \\   \hline 
			Bern.               &  1 & 4 & 9  \\   \hline 
		\end{tabular} 
	\end{center} 
	\caption{Number of modes in the periodic unit for different elements in the \textit{T} mesh.} \label{tab:number_of_modes_2D-meshT}
\end{table}%

\begin{figure}
	\centering
	\subfigure[\textit{X} mesh \label{fig:fourier_global_basic_SSPRK_CIPXbis}]{
	\includegraphics[width=\linewidth]{image_dispersion_2D/meshX/cfl_vs_tau_angleLoop/error_disp_cfl_tau_2D_lagrange_SSPRK_CIP.pdf}}
	\subfigure[\textit{T} mesh\label{fig:fourier_global_basic_SSPRK_CIPT}]{
		\includegraphics[width=\linewidth]{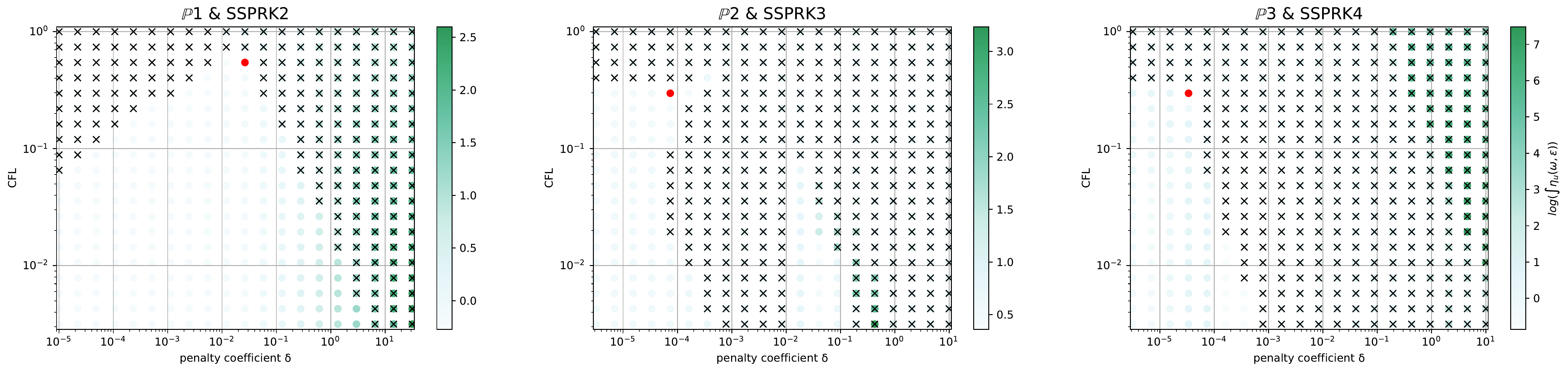}}
	\caption{$\log(\eta_u)$ values (blue scale) and stable area (unstable with black crosses), on $(\CFL,\delta)$  plane. The red dot denotes the optimal value. From left to right $\P_1$, $\P_2$, $\P_3$ \textit{Basic} elements with SSPRK scheme and CIP stabilization} 
\end{figure}

Even if for several methods we observe comparable results for the two mesh types, for some of them the analyses are quite different. An example is given by the \textit{Basic} elements with SSPRK schemes and CIP stabilization. For this method, we plot the dispersion error \eqref{eq:dispersionError} and the stability area in Figure~\ref{fig:fourier_global_basic_SSPRK_CIPXbis} for the \textit{X }mesh and in Figure~\ref{fig:fourier_global_basic_SSPRK_CIPT} for the \textit{T} mesh. We see huge differences in $\P_2$ and $\P_3$ where in the former a wide region becomes unstable for $\delta_L\leq \delta\leq \delta_R$ and for the latter we have to decrease a lot the value of $\delta$ to obtain stable schemes.

\begin{figure}
\subfigure[\textit{X} mesh \label{fig:fourier_global_cubature_SSPRK_OSSXbis}]{	\includegraphics[width=\linewidth]{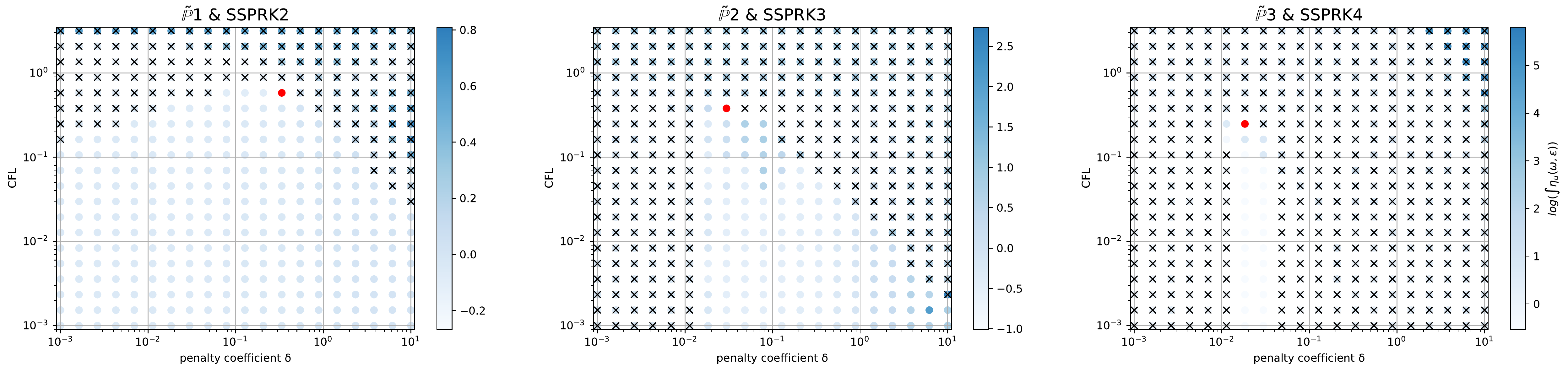}}
	\subfigure[\textit{T} mesh\label{fig:fourier_global_cubature_SSPRK_OSST}]{
	\includegraphics[width=\linewidth]{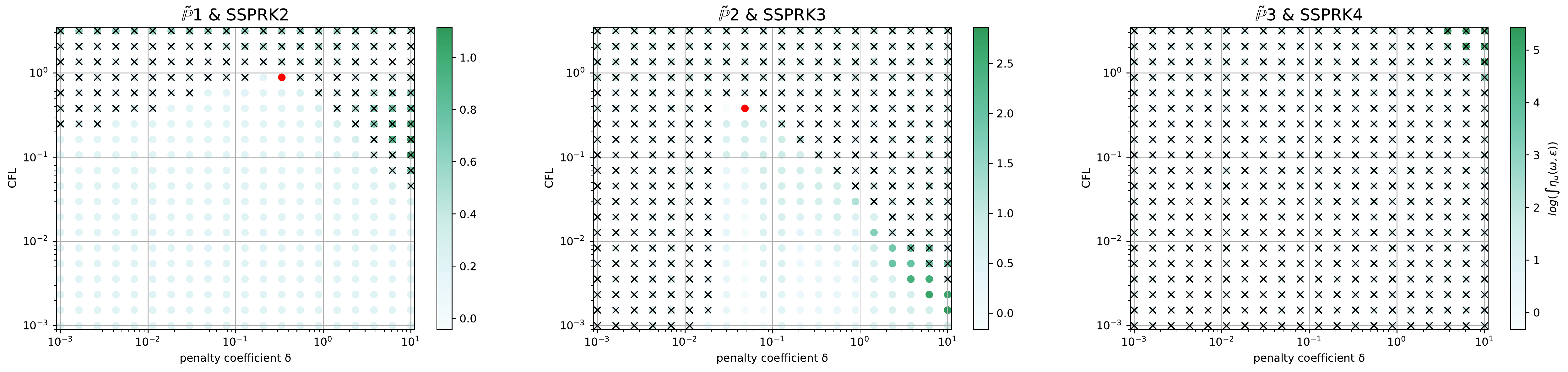}}
	\caption{$\log(\eta_u)$ values (blue scale) and stable area (unstable with black crosses), on $(\CFL,\delta)$  plane. The red dot denotes the optimal value. From left to right $\TP_1$, $\TP_2$, $\TP_3$ \textit{Cubature} elements with SSPRK scheme and OSS stabilization }	
\end{figure}

In the case of \textit{Cubature} elements with the OSS stabilization and SSPRK time integration, we have already seen in the previous section that the optimal parameters found were in a dangerous area. Repeating the stability analysis for the \textit{T} mesh we see that the situation is even more complicated. In Figure~\ref{fig:fourier_global_cubature_SSPRK_OSSXbis}  we plot the analysis for the \textit{X} mesh and in Figure~\ref{fig:fourier_global_cubature_SSPRK_OSST} the one for the \textit{T} mesh. $\TP_3$ elements, though being stable for some parameters for the \textit{X} mesh, are never stable on the \textit{T} mesh. This means, that, when searching general parameters for the schemes, we have to keep in mind that different meshes leads to different results.
%

For completeness, we present the optimal parameters also for the \textit{T} mesh in Table~\ref{tab:new_param_meshT_LinearAdvection-2D-RES} in Appendix~\ref{sec:app_fourier_param}.

In general, it is important to consider more mesh types when doing this analysis. In practice, we will use the two presented above (\textit{X} and \textit{T} meshes). In the following, we will consider the stability region as the intersection of stability regions of both meshes.
%


\subsection{Final results of the stability analysis} \label{sec:fourier2D_final_results}
\begin{figure}
	\centering
	\subfigure[SSPRK with OSS]{\includegraphics[height=0.2\linewidth,trim={1000 0 80 0}, clip]{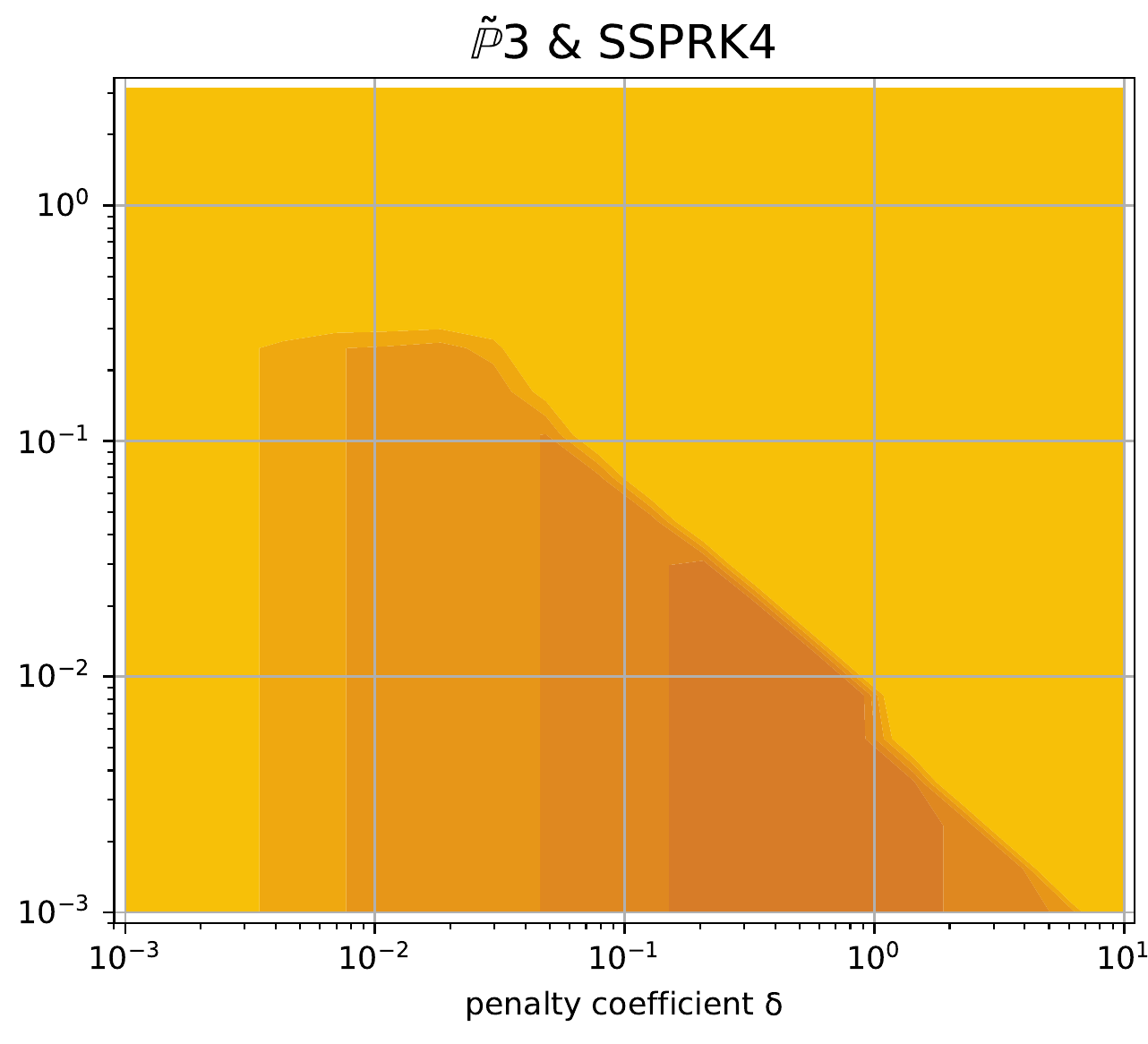}}
	\subfigure[SSPRK with CIP]{\includegraphics[height=0.2\linewidth,trim={1000 0 80 0}, clip]{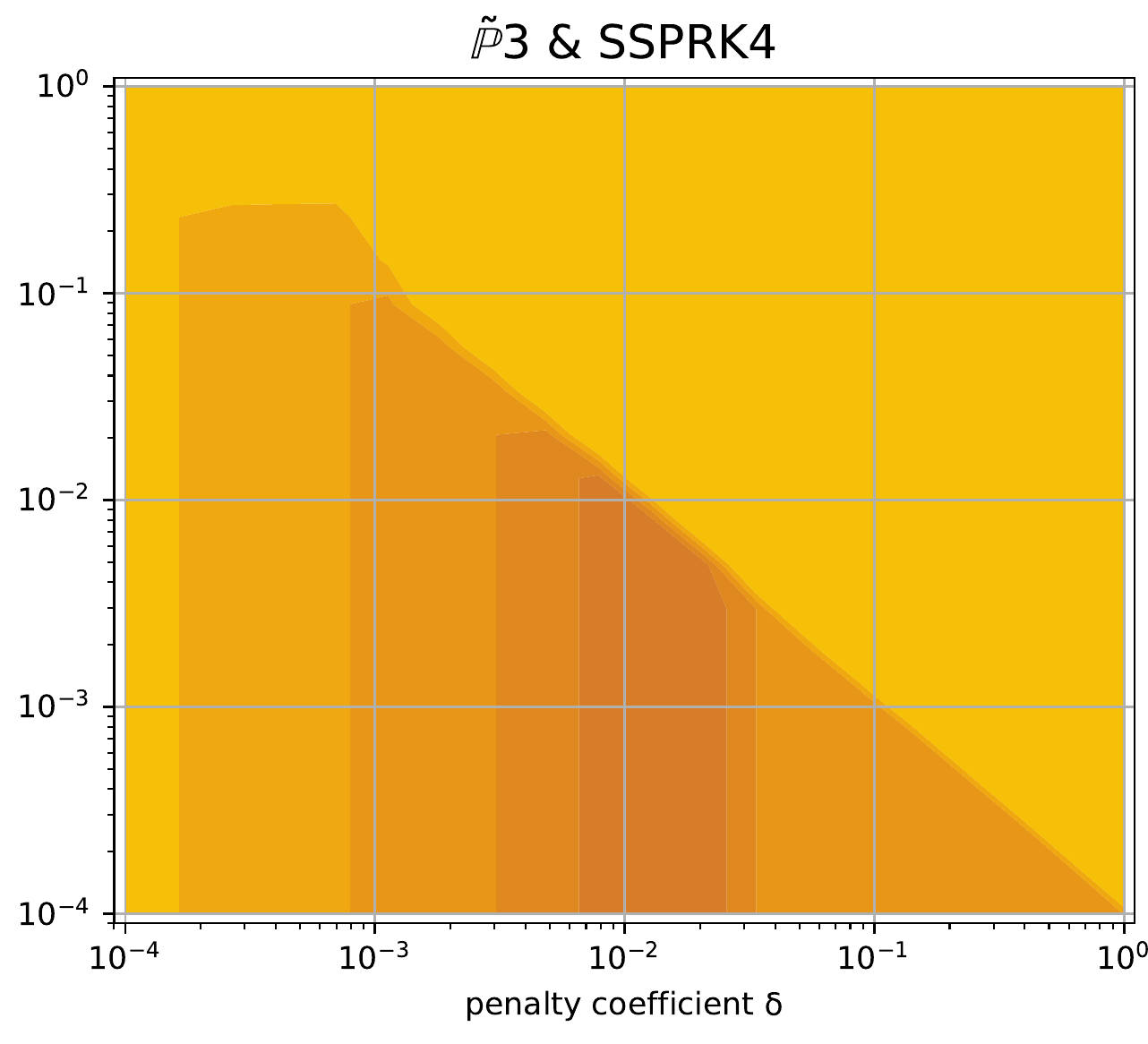}}
	\subfigure[DeC with OSS]{\includegraphics[height=0.2\linewidth,trim={1000 0 80 0}, clip]{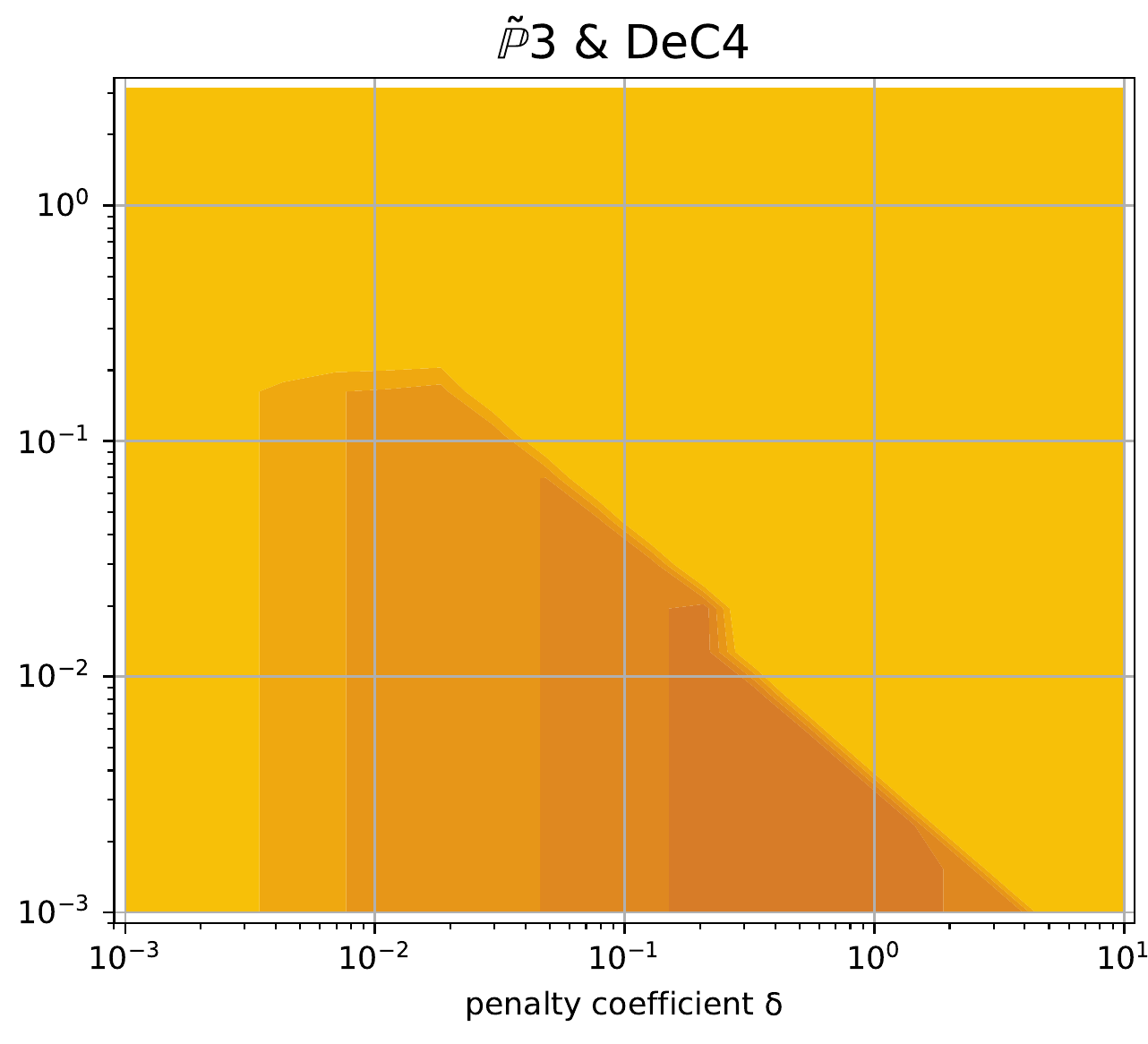}}
	\subfigure[DeC with CIP]{\includegraphics[height=0.2\linewidth,trim={1000 0 0 0}, clip]{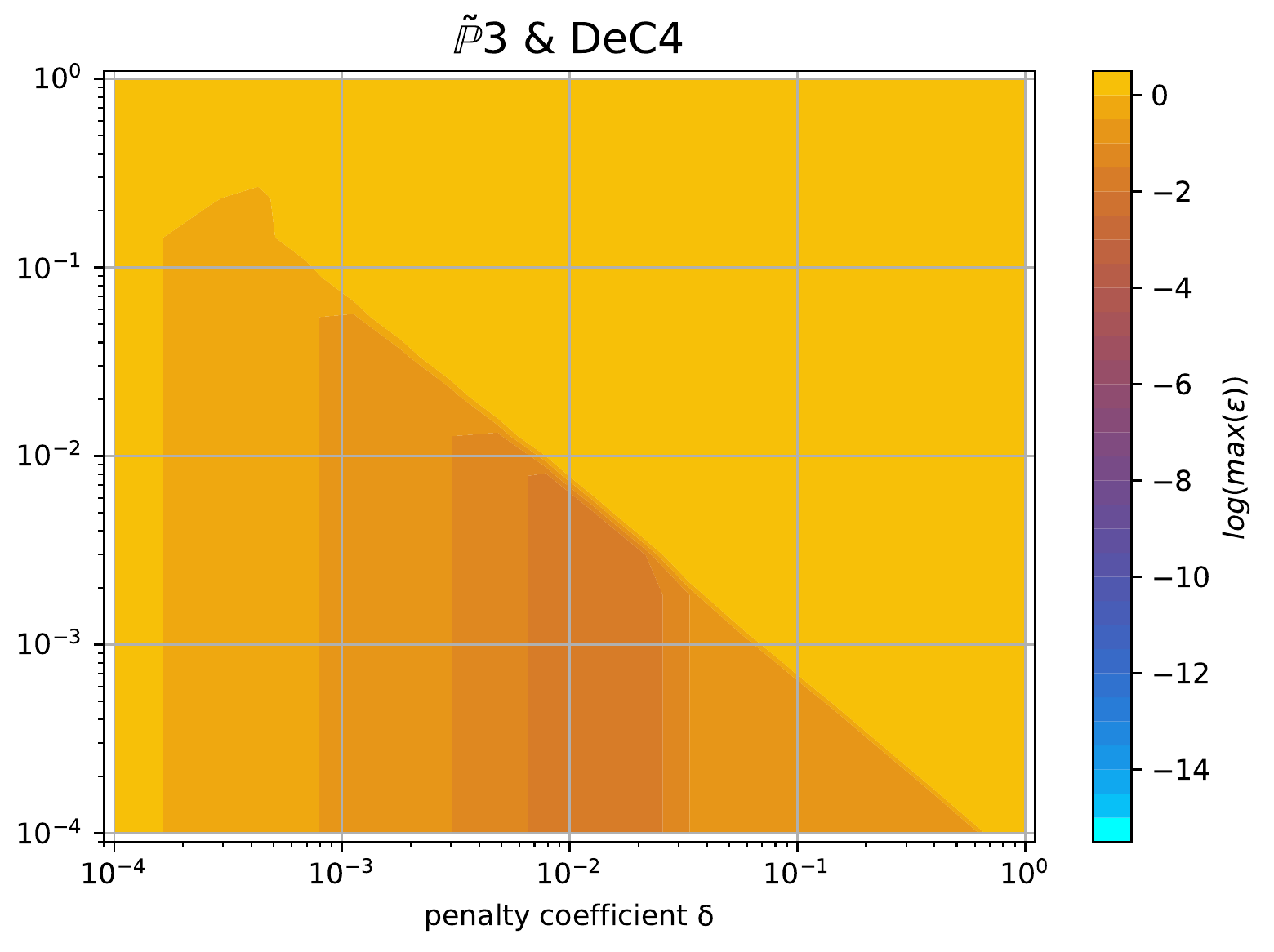}}
	\caption{Maximum logarithm of the amplification coefficient $\log(\max_i (\varepsilon_i))$ for  $\TP_3$ \textit{Cubature} elements on the \textit{X} and \textit{T} meshes} \label{fig:fourier_combined_cubP3_problem}
\end{figure}

\begin{figure}
	\centering
	\subfigure[SSPRK with OSS]{\includegraphics[height=0.2\linewidth,trim={1000 0 80 0}, clip]{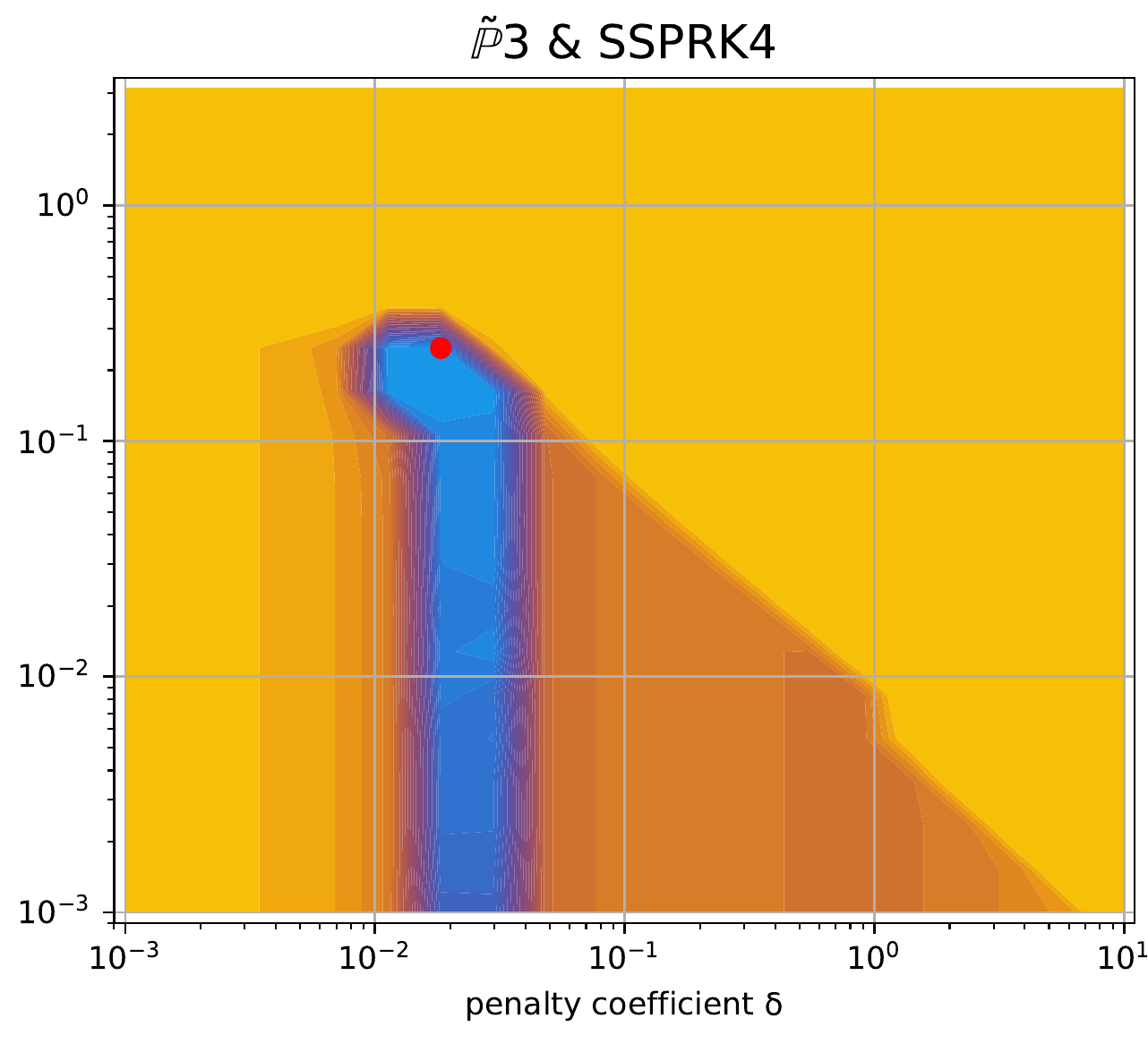}}
	\subfigure[SSPRK with CIP]{\includegraphics[height=0.2\linewidth,trim={1000 0 80 0}, clip]{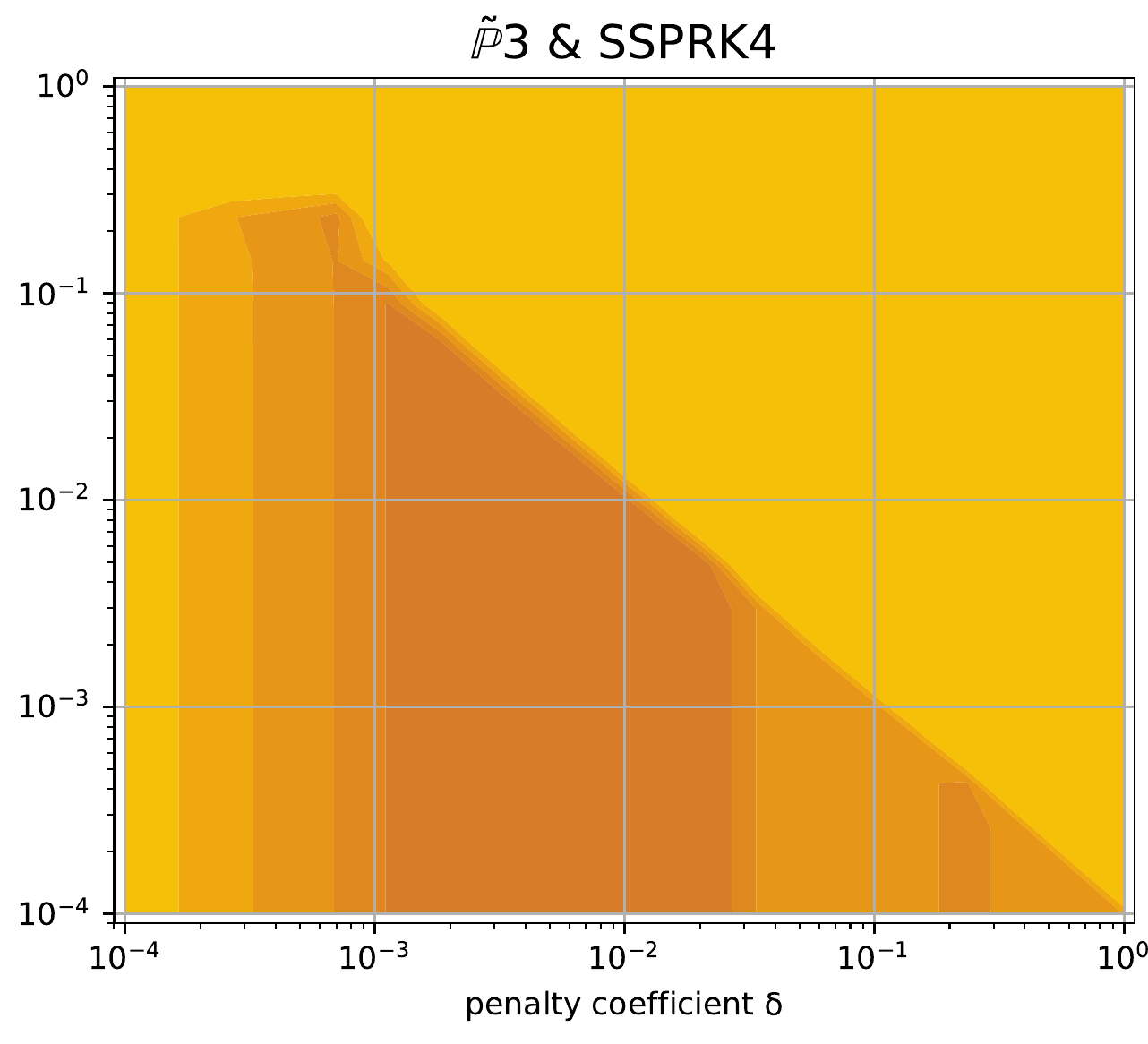}}
	\subfigure[DeC with OSS]{\includegraphics[height=0.2\linewidth,trim={1000 0 80 0}, clip]{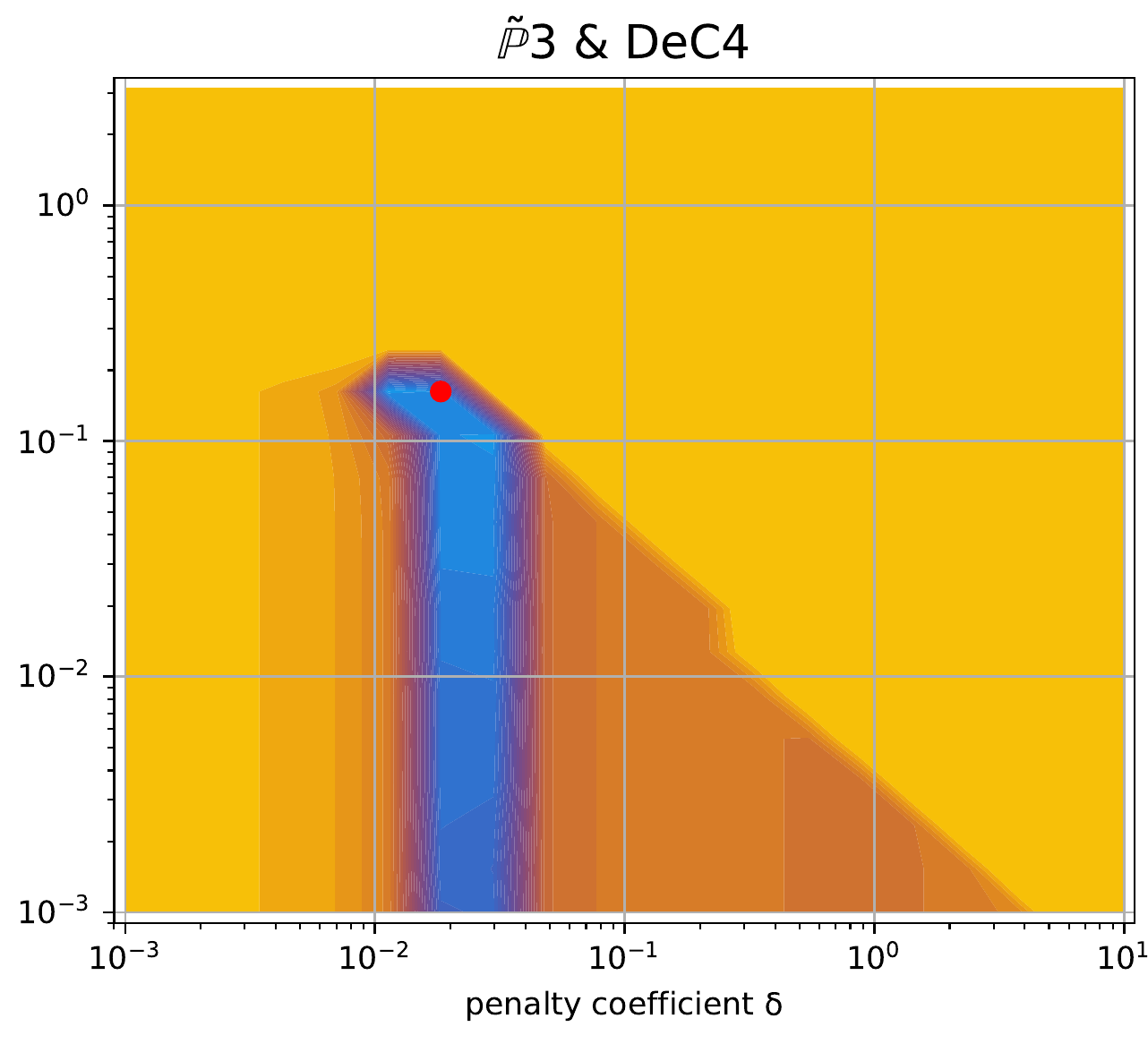}}
	\subfigure[DeC with CIP]{\includegraphics[height=0.2\linewidth,trim={1000 0 0 0}, clip]{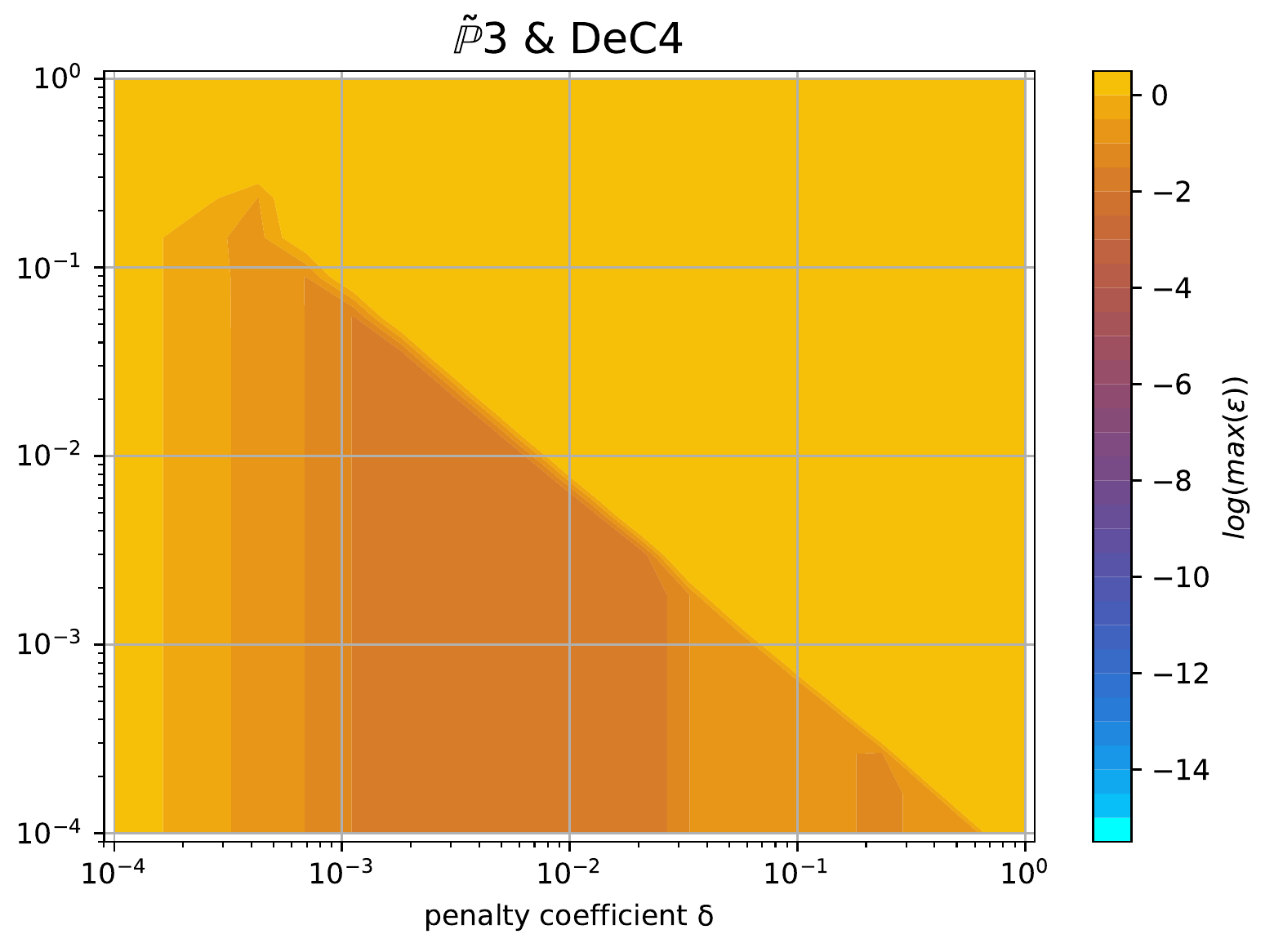}}
	\caption{Logarithm of the amplification coefficient $\log(\max_i (\varepsilon_i))$ for  $\TP_3$ \textit{Cubature} elements on the \textit{X} mesh} \label{fig:fourier_Xmesh_cubP3_problem}
\end{figure}
Taking into consideration all the aspects seen in the previous sections, it is important to have a comprehensive result, which tells which parameters can be used in the majority of the situations. A summary of the parameters obtained for the \textit{X} and \textit{T} mesh is available in Appendix~\ref{sec:app_fourier_param}.
In Table~\ref{tab:restrictive_param_LinearAdvection-2D-RES}, instead, we present parameters obtained using the most restrictive case among different meshes and that insure an enough big area of stability around them, as explained in Section~\ref{sec:fourier_rocket}. These parameters can be safely used in many cases and we will validate them in the numerical sections, where, first, we validate the results of the \textit{X} mesh on a linear problem on an \textit{X} mesh, then we used the more general parameters in Table~\ref{tab:restrictive_param_LinearAdvection-2D-RES} for nonlinear problems on unstructured meshes. \\
\begin{table} 
\small  
 \begin{center} 
		\begin{tabular}{| c | c || c | c | c | }  
	     \hline 
	     \multicolumn{2}{|c||}{Element $\&$ }  & \multicolumn{3}{|c|}{SUPG}  \\ \hline 
	     \multicolumn{2}{|c||}{ Time scheme }  & $\mathbb{P}_1$ & $\mathbb{P}_2$ & $\mathbb{P}_3$  \\ \hline \hline 
Basic              &  SSPRK  &  0.739 (0.127)   &  0.2 (0.1)$^*$   &  0.22 (0.026)  \\ 
         \hline 
       \multirow{2}{*}{ Cub.}              &  SSPRK  &  1.062 (0.28)   &  0.12 (0.13)$^*$   &  0.09 (0.05)$^*$  \\ 
               &  DeC  &  0.616 (0.28)   &  0.144 (0.078)   &  0.05 (0.05)$^*$  \\ 
         \hline 
Bern.              &  DeC  &  0.739 (0.298)   &  0.12 (0.45)$^*$   &  0.2 (0.153)$^*$  \\ 
         \hline 
        \end{tabular} 
    \end{center} 
 \begin{center} 
		\begin{tabular}{| c | c || c | c | c || c | c | c | }  
	     \hline 
	     \multicolumn{2}{|c||}{Element $\&$ }  & \multicolumn{3}{|c||}{OSS}  & \multicolumn{3}{|c|}{CIP}  \\ \hline 
	     \multicolumn{2}{|c||}{ Time scheme }  & $\mathbb{P}_1$ & $\mathbb{P}_2$ & $\mathbb{P}_3$   & $\mathbb{P}_1$ & $\mathbb{P}_2$ & $\mathbb{P}_3$  \\ \hline \hline 
Basic              &  SSPRK  &  0.403 (0.127)   &  0.2 (0.05)$^*$   &  0.22 (0.026)   &  0.403 (0.012)   &  0.1 (1.00e-03)$^*$   &  0.1 (5.00e-04)$^*$   \\ 
         \hline 
       \multirow{2}{*}{ Cub.}               &  SSPRK  &  0.58 (0.336)   &  0.2 (0.08)$^*$   &   0.28 (0.018)$^{**}$  &  0.58 (0.048)   &  0.06 (0.01)$^*$   &   /     \\ 
               &  DeC  &  0.379 (0.207)   &  0.12 (0.07)$^*$   &   0.162 (0.018)$^{**}$    &  0.379 (0.026)   &  0.025 (0.01)$^*$   &  /     \\ 
         \hline 
Bern.              &  DeC  &  0.173 (0.58)   &  0.02 (0.2)$^*$   &  0.015 (0.078)$^*$   &  0.173 (0.153)   &  0.012 (0.01)$^*$   &  0.001 (0.01)$^*$   \\ 
         \hline 
        \end{tabular} 
     \caption{Optimized CFL and penalty coefficient $\delta $ in parenthesis, combining the two mesh configurations.
     The values denoted by $^*$ are not the optimal one, but they lay in a safer region, see Section~\ref{sec:fourier_rocket}. 
     The values marked by $^{**}$ cannot be used on the \textit{T} mesh. ``/" means that it is unstable for every parameter.} \label{tab:restrictive_param_LinearAdvection-2D-RES}        
    \end{center} 
\end{table}%

A special remark must be done for \textit{Cubature} $\TP_3$ elements combined with the OSS and the CIP stabilizations. In Figure~\ref{fig:fourier_combined_cubP3_problem} we see how the amplification coefficient $\max_i \varepsilon_i$ has always values far away from zero. 
For the CIP stabilization this is always true and even for the $\TP_2$ elements the stability region is very thin. 
As suggested in \cite{Burman2020ACutFEmethodForAModelPressure,larson2019stabilizationHighOrderCut} higher order derivatives jump stabilization terms might fix this problem, but it introduces more parameters. This has not been considered here. Another remark is that the \textit{T} configuration is very peculiar and, as we will see, on classical Delauney triangulations the issue seem to not affect the results. Moreover, the use of additional discontinuity capturing operators may alleviate this issue as some additional, albeit small, dissipation is explicitly introduced in smooth regions.

In \cref{sec:app_fourier_visco}, we propose to add an additional stabilization term for these unstable schemes, i.e., \textit{Cubature} $\TP_3$ elements and OSS or CIP stabilization techniques. This term is based on viscous term \cite{abgrallEntropyConservative,2011JCoPh.230.4248G,KUZMIN2020104742,inbookLlobell2020} and allows to stabilize numerical schemes for any mesh configuration.\\


For the OSS stabilization we observe a similar behavior in Figure~\ref{fig:fourier_combined_cubP3_problem}. The stability that we see in that plot are only due to the the $T$ mesh. Indeed, for the OSS stabilization on the \textit{X} mesh there exists a corridor of stable values, which turn out to be unstable for the \textit{T} mesh, see Figure~\ref{fig:fourier_Xmesh_cubP3_problem}. 
In practice, also on unstructured grids we have not noticed instabilities when running with the parameters found with the \textit{X} mesh. 
 Hence, we suggest anyway some values of CFL and $\delta$ for these schemes, which are valid for the \textit{X} mesh, noting that they might be dangerous for very simple structured meshes.
The validation on unstructured meshes also for more complicated problems will be done in the next sections.  \\


Overall, Table~\ref{tab:restrictive_param_LinearAdvection-2D-RES} gives some insight on the efficiency of the schemes. We remind that, in general, we prefer matrix free schemes, so this aspect must be kept in mind while evaluating the efficiency of the schemes. All the SUPG schemes, except when with DeC, and all the \textit{Basic} element schemes have a mass matrix that must be inverted. Among the others we see that for first degree polynomials schemes the DeC with \textit{Bernstein} polynomials and SUPG stabilization gives one of the largest CFL result,  while for second degree polynomials the OSS \textit{Cubature} SSPRK scheme seems the one with best performance and, for fourth order schemes, again the \textit{Bernstein} DeC SUPG is one of the best.

In conclusion of this section, there are important points to highlight:
\begin{itemize}
	\item The extension of the Fourier analysis to the two-dimensional space leads to significantly different results with respect to the one-dimensional one. Both in terms of global stability of the schemes, and in terms of optimal parameters. Moreover, in opposition to \cite{michel2021spectral}, \textit{Bernstein} elements with SUPG stabilization technique lead to stable  and efficient schemes. \textit{Cubature} elements, which were the most efficient in one-dimensional problems, have stability issues on the two-dimensional mesh topologies studied. 
	\item The complexity of the analysis in two-dimensional space is increased. This not only implies a larger number of degrees of freedom, but also more parameters to keep into account, including the angle of the advection term and the possible different configuration of the mesh. 
	The visualization of the stability region of the time scheme as shown in Figure~\ref{fig:rocket_meshX_cub_SSPRK_SUPG_P3_unstable} with the eigenvalues of the semi-discretization operators helps in understanding the effect of CFL and penalty coefficient on the stability of the scheme, only for methods of lines. This helps in choosing and optimizing the couple of parameters. 
\end{itemize}

\begin{remark} 
Another possibility to characterize the linear stability of numerical method is proposed by J. Miller \cite{articleMiller1971}. This method is based on the study of the characteristic polynomial of the amplification matrix $G$. 
However, this method does not provide information about the phase $\omega$, since it does not compute eigenvalues of $G$. For this reason, we choose the eigenanalysis.
\end{remark} 
%
\subsection{Accounting for discontinuity capturing corrections }\label{sec:app_fourier_visco}

The stabilization terms accounted for so far are linear stabilization operators. For more challenging simulations, 
additional non-linear stabilization techniques might be added to control the numerical solution in vicinity of strong non-linear fronts
and/or discontinuities. We consider here the effect of adding an extra viscosity term, as in 
the entropy stabilization formulations proposed e.g.  in 
\cite{abgrallEntropyConservative,kuzmin2020entropyDG,2011JCoPh.230.4248G,KUZMIN2020104742,inbookLlobell2020}.
We in particular look at the approach proposed in   \cite{2011JCoPh.230.4248G}, and used for shallow water waves 
in \cite{10.1007/978-3-319-19800-2_36,inbookLlobell2020} and in \cite{ARPAIA2022101915,filippini:hal-01824108}.
In this approach the viscosity is designed to provide a first order correction $\mu_K=\mathcal{O}(h)$ close to
discontinuities, while for smooth enough solutions $\mu_K =  c h^{p+1}$.

Our idea is to embed this high order correction explicitly in the analysis of the previous section to provide a heuristic
characterization of the fully discrete stability of the resulting stabilized formulation: 
find $u_h\in V_h^p$ that satisfies  for any $v_h\in W_h$
\begin{equation}
\int_{\Omega} v_h ( \partial_t u_h   + \nabla \cdot f(u_h)) dx + \underbrace{ S(v_h,u_h)}_{\mbox{Diffusive term}} 
+ \underbrace{\sum_K \int_K \mu_K(u_h) \nabla v_h \cdot \nabla u_h}_{\mbox{Viscosity term}} =0. \label{eq:FV-visco_app}
\end{equation}

%
\subsubsection{Note on the stability of the method} 
As it is done for previous stabilization terms in \cref{sec:stabilization}, we can characterize the accuracy of this method estimating the truncation error for a polynomial approximation of degree $p$. Considering the smooth exact solution $u^e(t,x)$ of \eqref{eq:FV-visco_app}, for all functions $\psi$ of class at least $\mathcal{C}^1(\Omega)$ of which  $\psi_h$  denotes the finite element projection, we obtain
\begin{equation}
\begin{split}
\epsilon(\psi_h) := \Big|
\int_{\Omega_h} &\psi_h \partial_t (u_h^e - u^e) \; dx - \int_{\Omega_h} \nabla \psi_h \cdot (f(u_h^e)-f(u^e))\; dx \\ 
              +& \sum\limits_{K\in \Omega_h}\mu_K \int\limits_{K} \nabla \psi_h \cdot \nabla(  u^e_h - u^e  ) dx
           	\Big| \le C h^{p+1},
\end{split}
\label{eq_entrop2}
\end{equation}
with $C$ a constant independent of $h$. The estimate can be derived from standard approximation results applied to $u_h^e-u^e$ and to its derivatives, knowing that $\mu_K = \mathcal{O}(h^{p+1})$. 

Then, for a linear flux, periodic boundaries and taking $\mu_K=\mu$ constant along the mesh, we can test with  $v_h=u_h$ in \eqref{eq:FV-visco_app}, we get 
\begin{equation}
\begin{split}
\int\limits_{\Omega_h} d_t \frac{ u^2_h}{2} = - \sum\limits_{K} \int\limits_{K} \mu ( \nabla  u_h )^2 \leq 0,
\end{split}
\label{eq_entrop3}
\end{equation}
which can be integrated in time to obtain a bound on the $\mathbb{L}_2$ norm of the solution. 

\subsubsection{The von Neumann analysis} 
\begin{figure}[h]
		\centering
	\subfigure[$\mu = 0$ \label{fig:fourier_cubP3_SSPRK_OSS_viso0}]{\includegraphics[height=0.25\linewidth,trim={950 0 80 0}, clip]{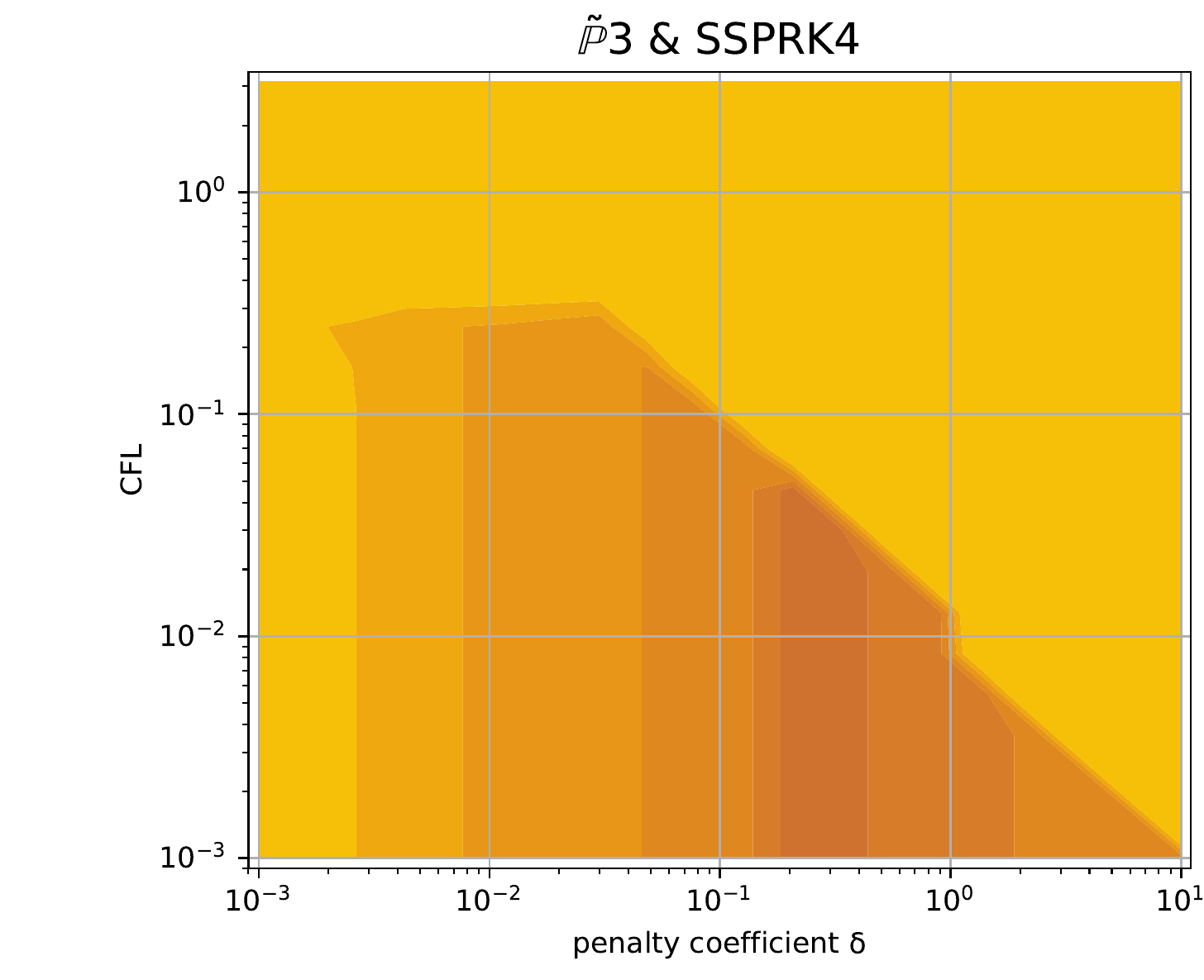}} \hfill
	\subfigure[$\mu   = 0.005 h_K^{p+1}$ \label{fig:fourier_cubP3_SSPRK_OSS_viso0005}]{\includegraphics[height=0.25\linewidth,trim={0 0 80 0}, clip]{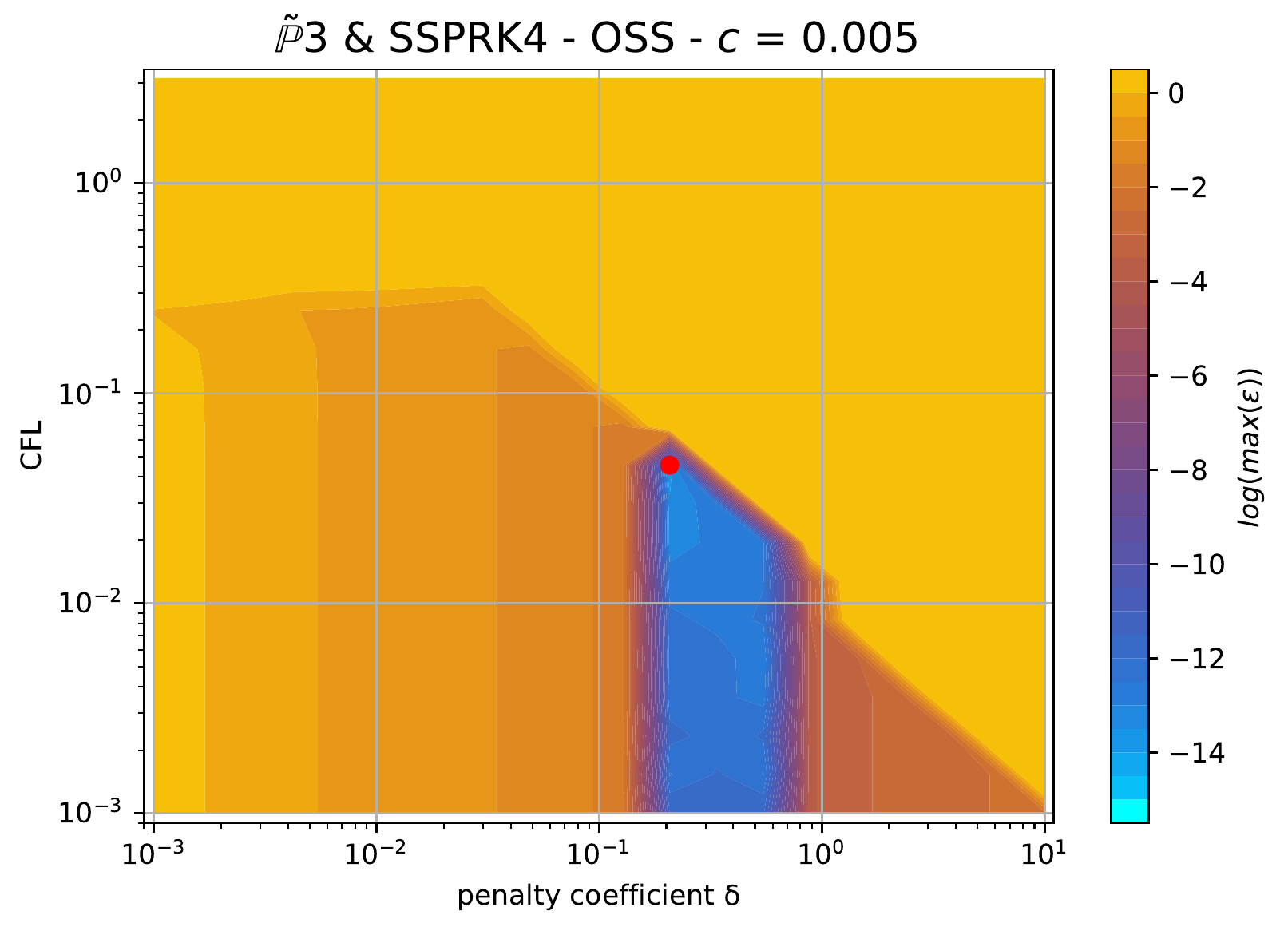}} \hfill
	\subfigure[$\mu   = 0.05 h_K^{p+1}$ \label{fig:fourier_cubP3_SSPRK_OSS_viso005}]{\includegraphics[height=0.25\linewidth,trim={0 0 80 0}, clip]{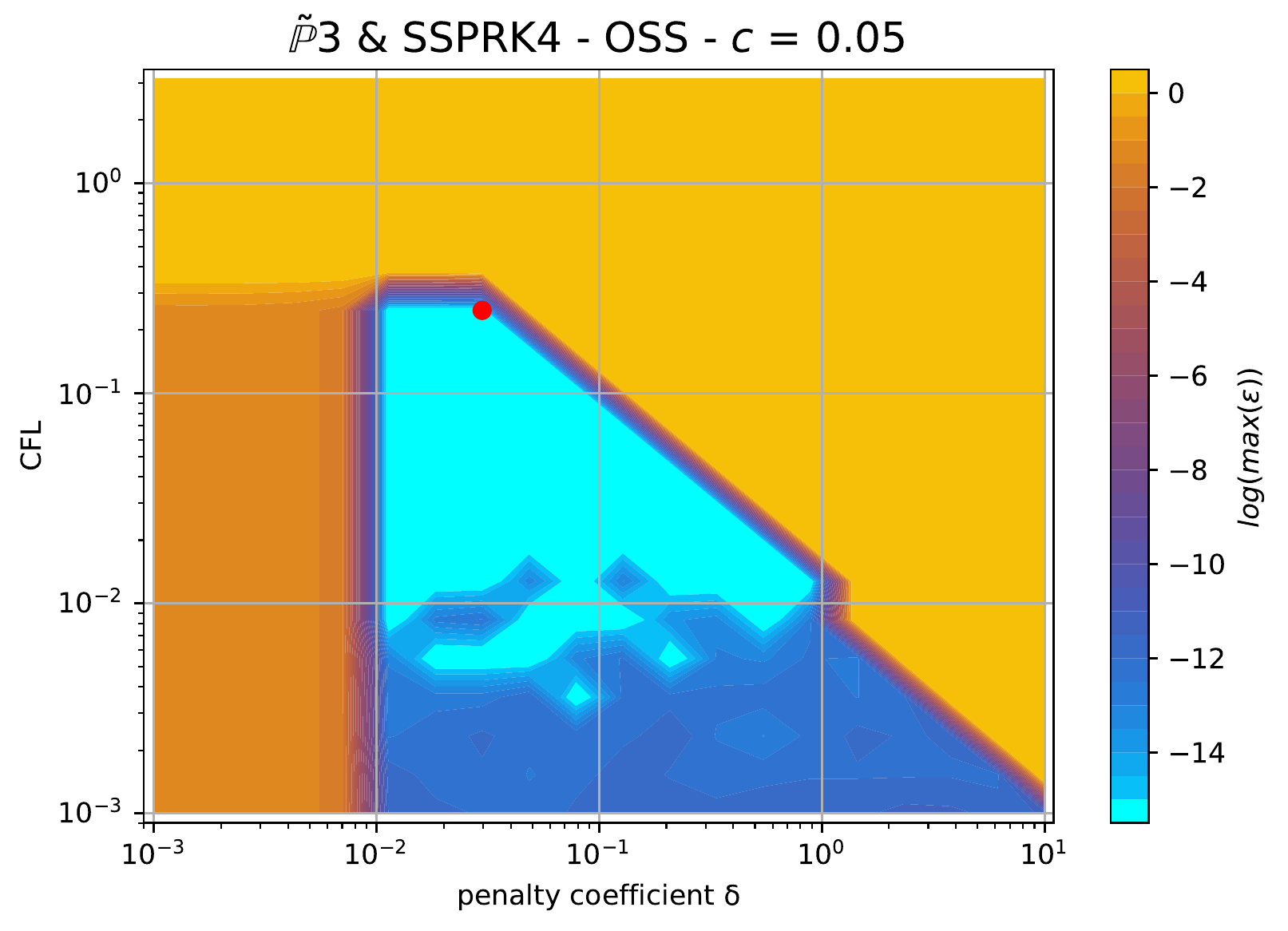}} \hfill \\
	\subfigure[$\mu  = 0.5 h_K^{p+1}$ \label{fig:fourier_cubP3_SSPRK_OSS_viso05}]{\includegraphics[height=0.25\linewidth,trim={0 0 80 0}, clip]{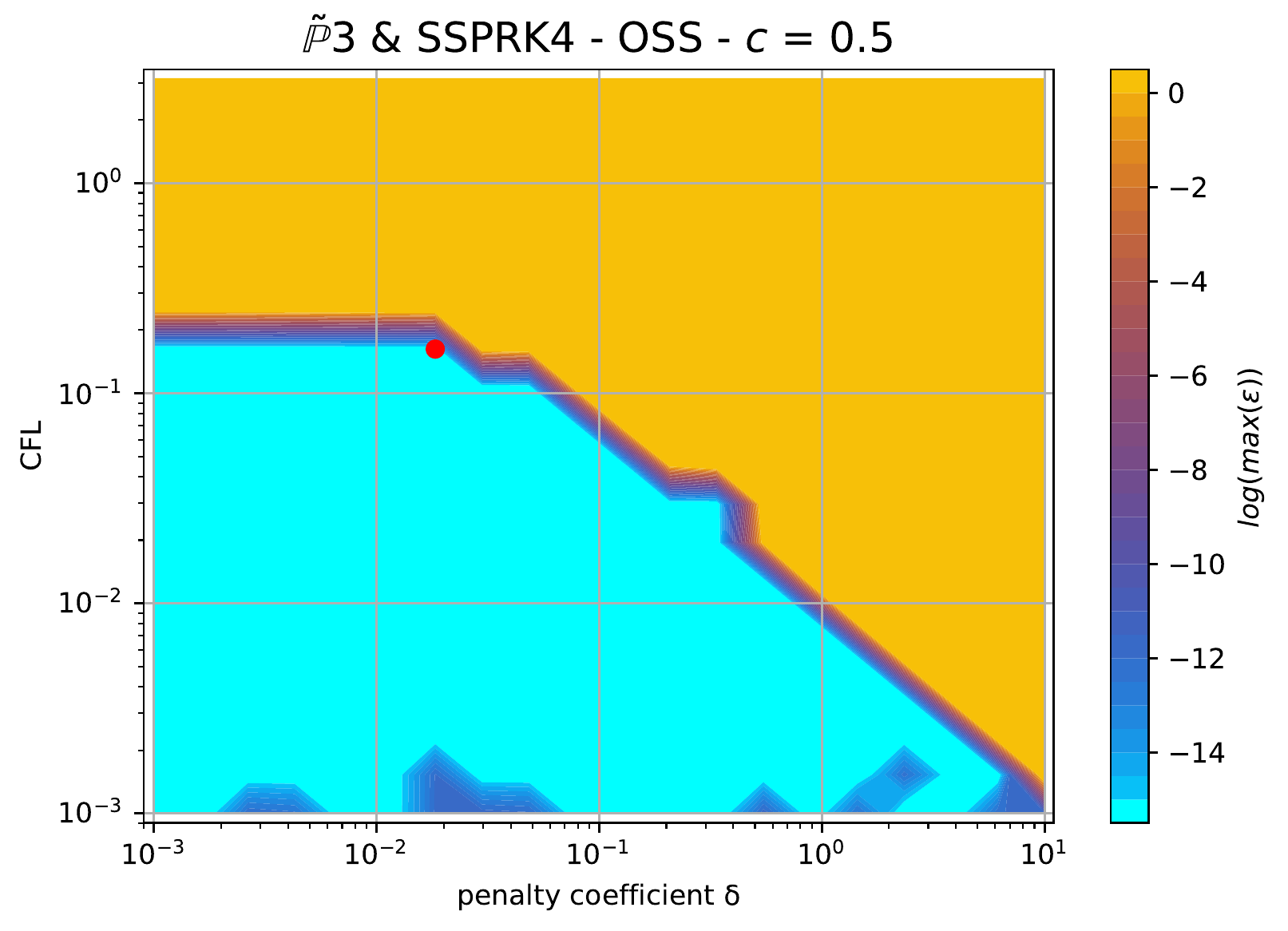}} 
	\subfigure[$\mu  = 5 h_K^{p+1}$ \label{fig:fourier_cubP3_SSPRK_OSS_viso5}]{\includegraphics[height=0.25\linewidth,trim={0 0 0 0}, clip]{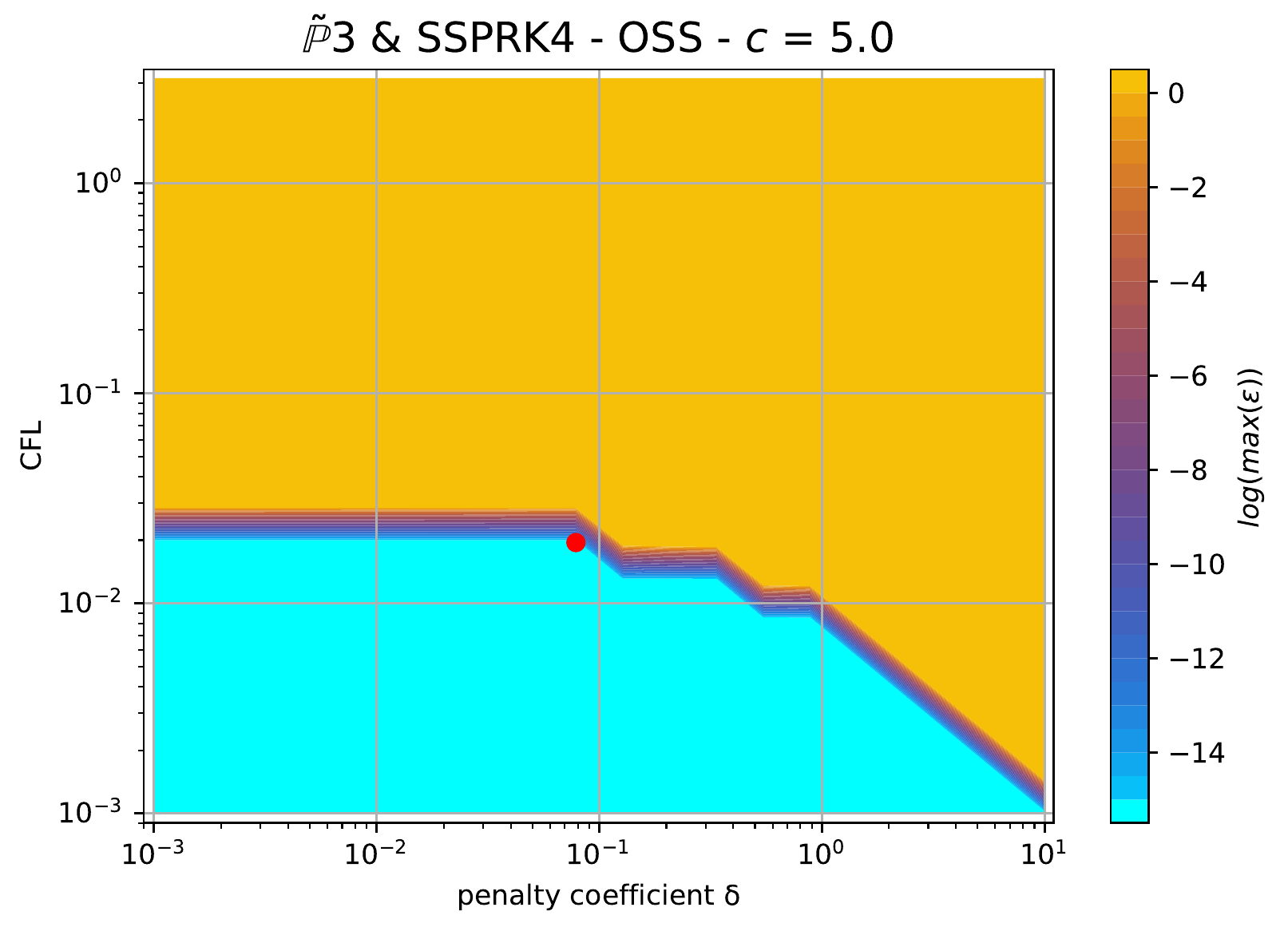}} 
	\caption{\textit{T} mesh - Von Neumann analysis using an additional viscosity term (see \eqref{eq:FV-visco_app}). \textit{Cubature} $\TP_3$ elements with SSPRK and OSS. Comparison of different $\mu$.} \label{fig:fourier_cubP3_SSPRK_OSS_viso}
\end{figure}

As we saw in \cref{sec:fourier2D_final_results}, the \textit{T} mesh configuration has stability issues. In particular, the numerical schemes using \textit{Cubature} $\TP_3$ elements, SSPRK and DeC time integration methods, and the OSS and the CIP stabilization techniques are unstable. We propose to evaluate these schemes adding the viscosity term in \eqref{eq:FV-visco_app}. For the von Neumann analysis, we use $\mu_K(u) = c h_K^{p+1}$ in \eqref{eq:FV-visco_app}, with $c \in \R^+$, $h_K$ the cell diameter and $p$ the degree of polynomial approximation. We show the plot of $\max_i  \epsilon_i$ to understand how the stability region behaves with respect to $c$ using \textit{Cubature} $\TP_3$ elements. In \cref{fig:fourier_cubP3_SSPRK_OSS_viso} the maximum amplification factor $\epsilon$ is represented for varying $c$, using the OSS stabilization technique and the SSPRK time integration method. We note that the same behaviour is observed with CIP and DeC. Plots are available online \cite{TorloMichel2021git}. 

We can observe two main results. First, increasing the parameter $c$ up to around 0.1 allows to expand the stability region. Second, when the viscosity coefficients reaches too high values, it is necessary to decrease the CFL (see \cref{fig:fourier_cubP3_SSPRK_OSS_viso005} with $\mu=0.05$ and \cref{fig:fourier_cubP3_SSPRK_OSS_viso05} with $\mu=0.5$ as an example).
\section{Numerical verification} \label{sec:validation_num-meshX}
We now perform numerical tests to check the validity of our theoretical findings. 
We initially focus on the structured grids, and in particular on the  \textit{X} mesh configuration, although similar verifications
have been performed on the  \textit{T} mesh.
We will use elements of degree $p$, with $p$ up to 3, with time integration schemes of the corresponding order of accuracy to ensure an overall error of $\mathcal{O}(\Delta x ^{p+1})$, under the CFL conditions discussed earlier  (see also  \cref{tab:new_param_meshX_LinearAdvection-2D-RES} in \cref{sec:app_fourier_param}). As already stressed, numerical integration is performed with Gauss--Legendre  \emph{formulae} of the appropriate order
to exactly integrate the variational form for  \textit{Basic} and \textit{Bernstein} elements, while
for \textit{Cubature} elements we use those 
associated to the interpolation points.

The mesh used in the Fourier analysis is the basis of the one we will use in the numerical simulations. We will extend it periodically for the whole domain, see an example in \cref{fig:2d_meshX}.

\subsection{Linear advection equation test} \label{sec:LA2D_results-meshX}
We start with the linear advection  equation \ref{eq:conservation_law2D} on the domain  $\Omega = [0,2]\times[0,1]$ using Dirichlet inlet boundary conditions:
\begin{equation}\label{lin_adv_2D}
	\begin{cases}
    	\partial_t u (t,\mathbf{x}) + \ba \cdot \nabla u (t,\mathbf{x})  = 0,  \qquad & \quad (t,\mathbf{x}) \in [t_0,t_f] \times \Omega, \quad \ba = (a_x,a_y)^T \in \mathbb{R}^2,  \\
    	 u (0,\mathbf{x}) = u_0(\mathbf{x}), & \\
    	 u (t,\mathbf{x_D}) = u_{ex} (t,\mathbf{x_D}), & \quad \mathbf{x_D} \in \Gamma_D = \{ (x,y)\in \R^2 ,x \in \{0,2\} \mbox{ or } y \in \{0,1\} \},
	\end{cases}
\end{equation}
where $u_0((x,y)^T) = 0.1 \cos(2\pi\, r(x,y) )$, with $r(x,y)=\cos(\theta)x+\sin(\theta)y$ the rotation by an angle $\theta$ around $(0,0)$, $\mathbf{a}=(a_x,a_y)^T=(\cos(\theta), \sin(\theta) )^T$ and $\theta=3\pi/16$. The final time of the simulation is $t_f=2s$. 

\begin{figure}
		\centering
\subfigure[\textit{X} mesh on  $\Omega=\left(0,2\right) \times \left( 0,1\right)$  \label{fig:2d_meshX}]{\includegraphics[width=0.46\textwidth]{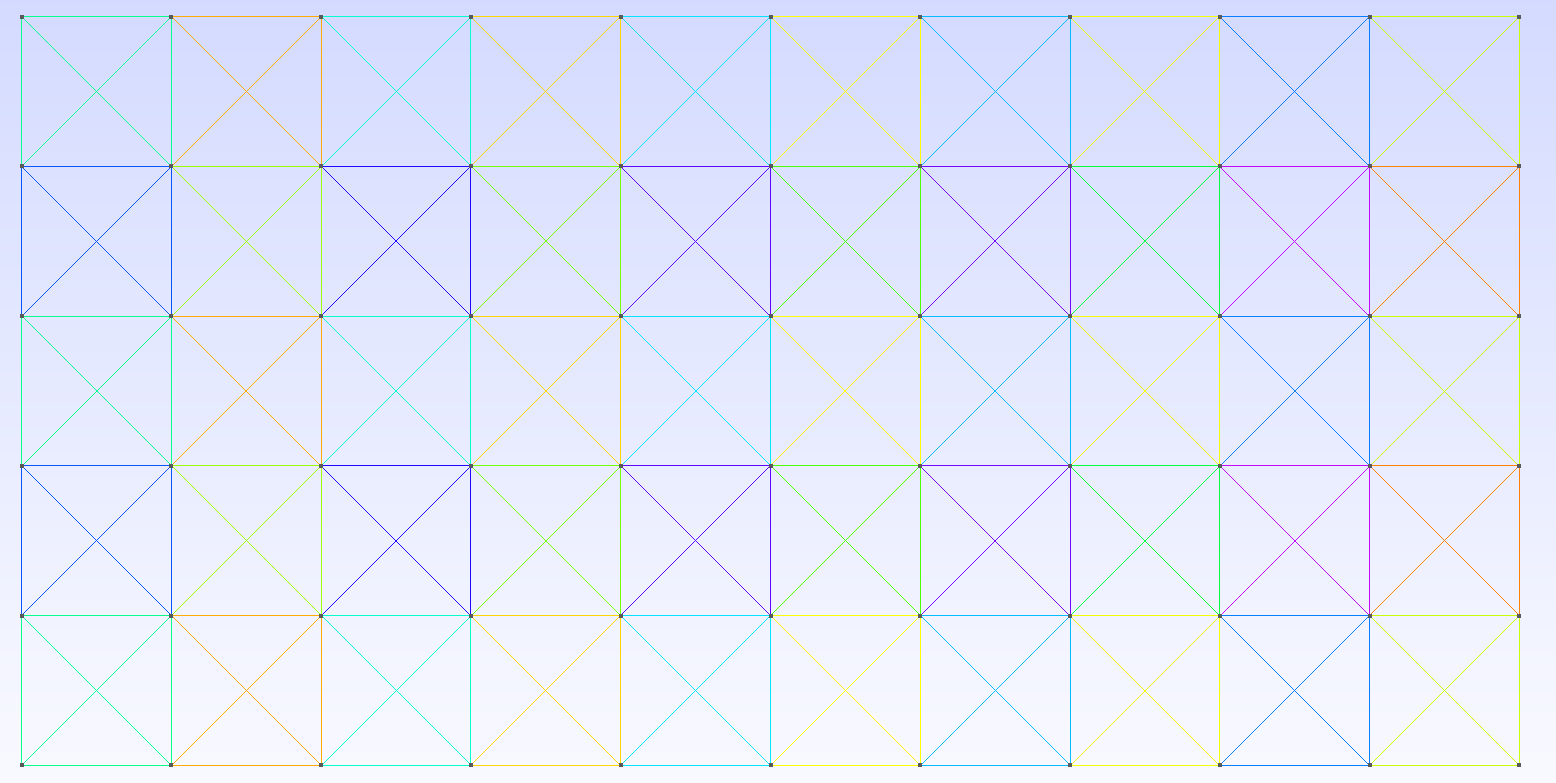}}
\hfill \subfigure[Cosinus test case with $\theta=3\pi/16$ \label{fig:ICcosLA}]{ \includegraphics[width=0.46\textwidth]{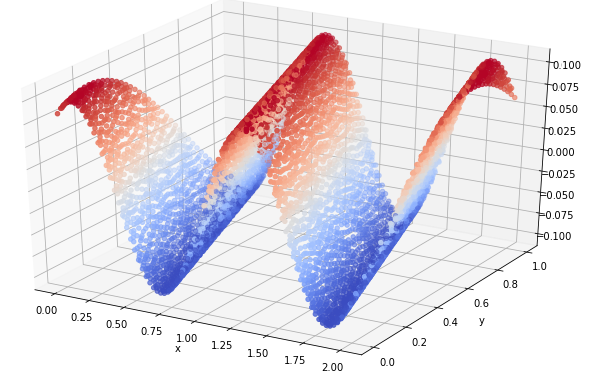}}
\caption{Linear advection simulation on the \textit{X} mesh}
\end{figure}

The exact solution is $u_{ex}(\mathbf{x},t)=u_0(x-a_x\, t,y-a_y\, t )$ for all $\mathbf{x}=(x,y)\in\Omega$ and $t\in\R^+$. The initial conditions are displayed in \cref{fig:ICcosLA}.
We discretize the domain with 
 the \textit{X} mesh pattern, see \cref{fig:2d_meshX}.
To have approximately the same number of degrees of freedom for different degrees $p$, 
we  use different mesh sizes for each order of accuracy: $\Delta x_1 = \{ 0.1, 0.05, 0.025 \}$ for $\mathbb{P}_1$, $\Delta x_2 = 2\Delta x_1 $ for $\mathbb{P}_2$, and $\Delta x_3 = 3 \Delta x_1 $ for $\mathbb{P}_3$ elements. 



\begin{figure}
	\centering
	\subfigure[\textit{Basic} elements \label{fig_err_2d_EDP-SSPRK_lag-OSS} ]{		\includegraphics[width=0.43\textwidth]{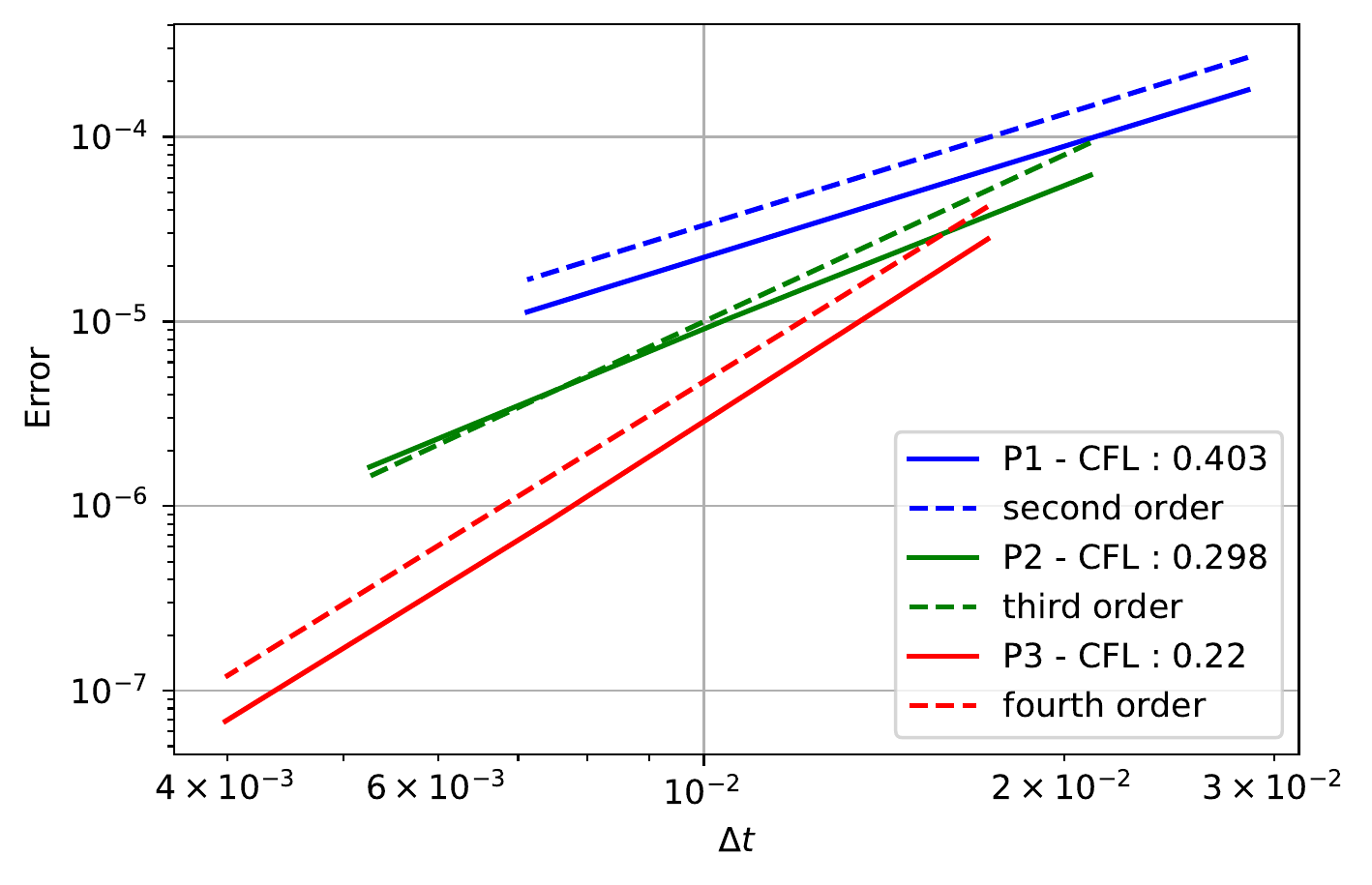}}
	\hfill
	\subfigure[\textit{Cubature} elements \label{fig_err_2d_EDP-SSPRK_cohen-OSS} ]{\includegraphics[width=0.41\textwidth]{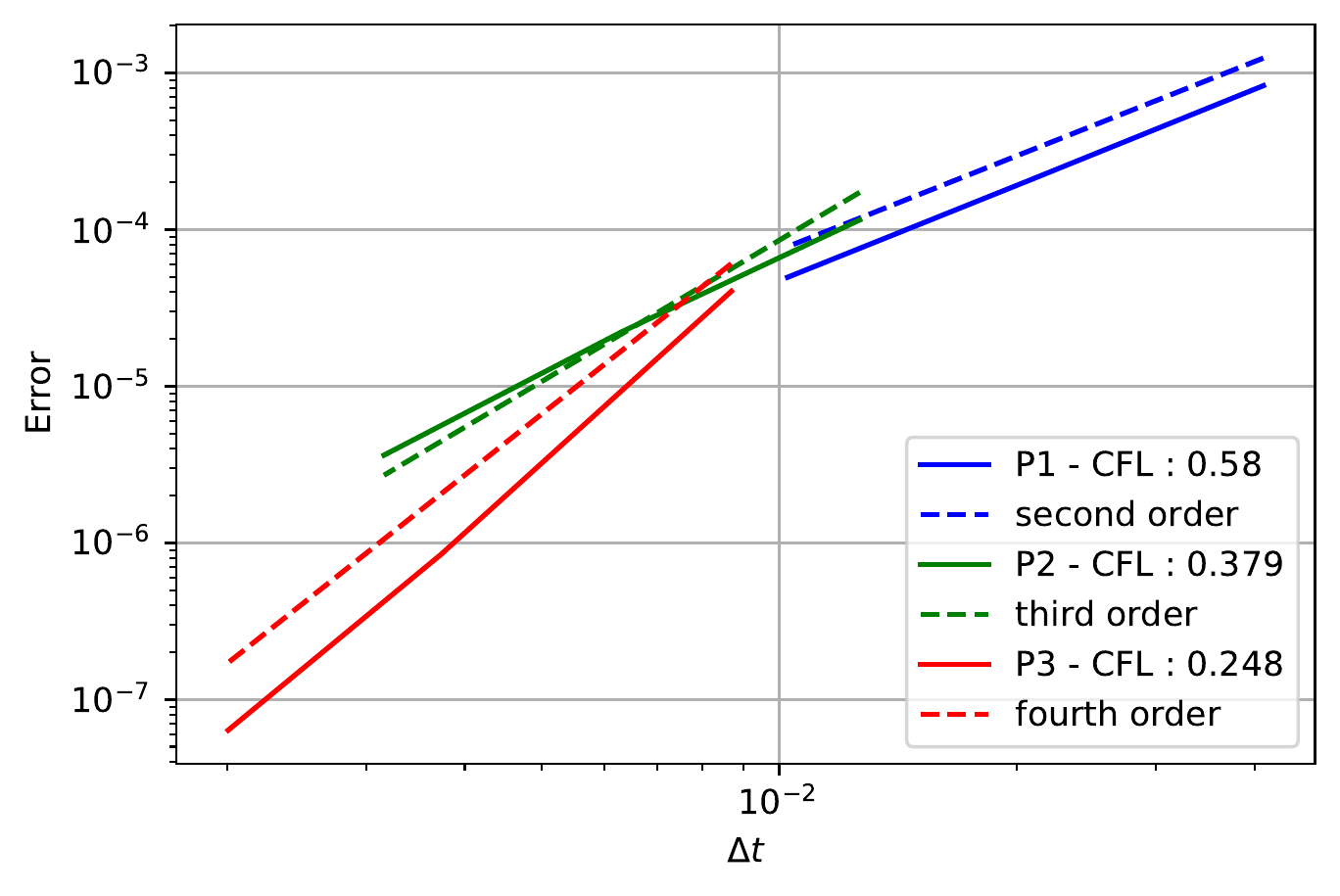}}
	\caption{Error decay for linear advection problem with different elements and OSS stabilization and SSPRK time discretization: $\mathbb{P}_1$ in blue, $\mathbb{P}_2$ in green and $\mathbb{P}_3$ in red}
\end{figure}
A representative result is provided in  \cref{fig_err_2d_EDP-SSPRK_lag-OSS,fig_err_2d_EDP-SSPRK_cohen-OSS}: it shows a comparison between \textit{Cubature} and \textit{Basic} elements with OSS stabilization and SSPRK time integration. As we can see, the two schemes have very similar errors except for $\P^1$ where the larger CFL increases the error. The \textit{Basic} elements require stricter CFL conditions, see \cref{tab:new_param_meshX_LinearAdvection-2D-RES}, and have larger computational costs because of the inversion of the mass matrix.

\begin{figure}
	\begin{center}
		\subfigure[\textit{Cubature} elements]{
		\includegraphics[height=0.19\textheight,trim={0 0 64mm 0}, clip]{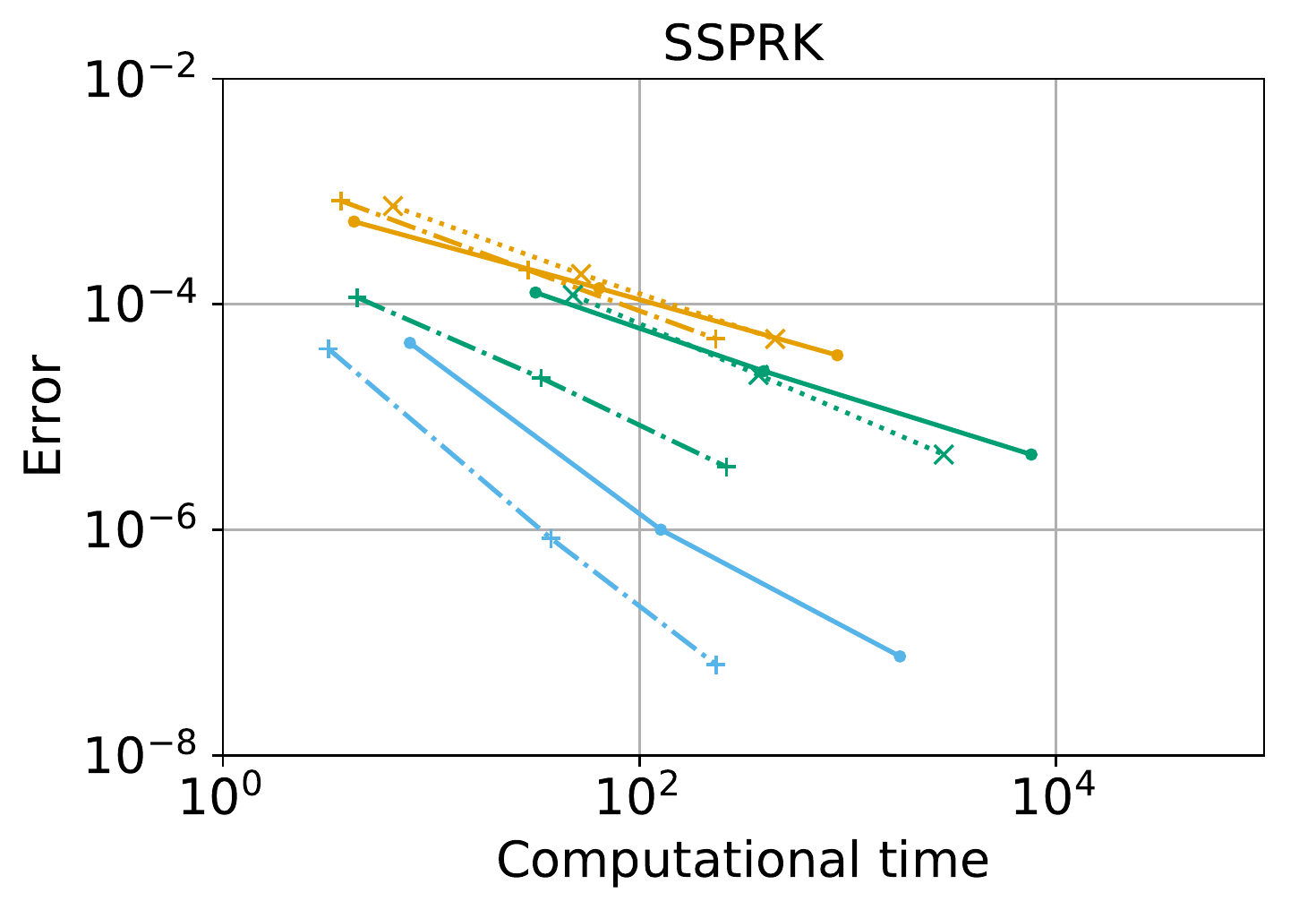}}
	\subfigure[\textit{Basic }elements]{
		\includegraphics[height=0.19\textheight,trim={8mm 0 0 0}, clip]{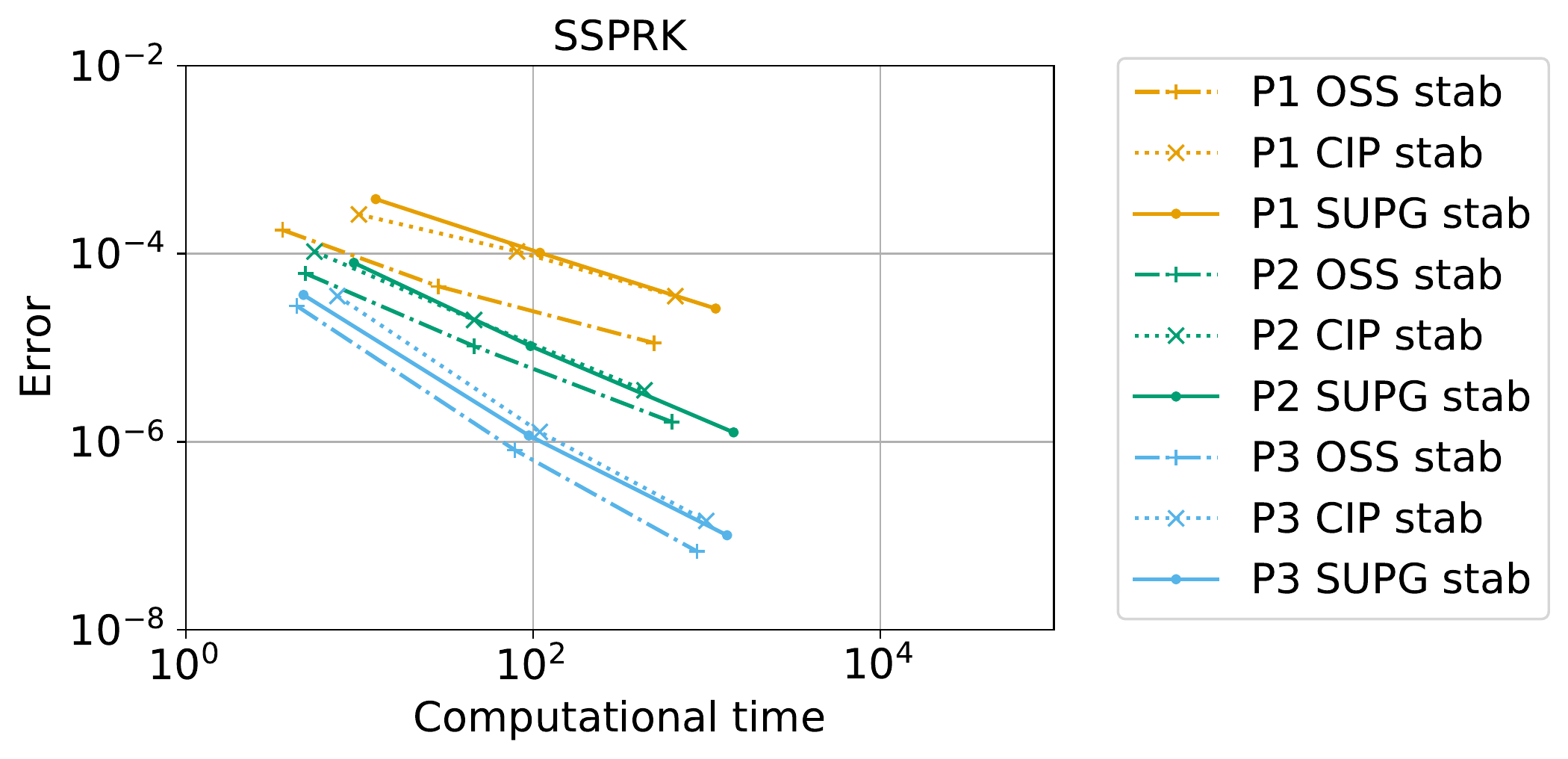}}
	\caption{Error for linear advection problem \eqref{lin_adv_2D} with respect to computational time for SSPRK time discretization, comparing \textit{Basic }and \textit{Cubature }elements and all stabilization techniques}
	\label{fig:2D_timeVsErrorLinAdv_meshX}
	\end{center}
\end{figure}

To show the main benefit of using the \textit{Cubature} elements (diagonal mass matrix), we plot  in \cref{fig:2D_timeVsErrorLinAdv_meshX} the computational time of \textit{Basic} and \textit{Cubature} elements for the SSPRK time scheme and all stabilization techniques.
As a first interesting result of numerical test, looking at the \cref{fig:2D_timeVsErrorLinAdv_meshX}, we can clearly see that, for a fixed accuracy, \textit{Cubature} elements obtain better computational times with respect to \textit{Basic} elements. Moreover, as expected, the SUPG stabilization technique requires more computational time as it requires the inversion of a mass matrix, even in the case where the CFL used in is larger than the ones for OSS or CIP stabilization, see \cref{tab:new_param_meshX_LinearAdvection-2D-RES}.

The order of accuracy reached by each simulations is shown in \cref{tab:conv_order_meshX_LinearAdvection-2D-RES}. The plots and all the errors are available at the repository \cite{TorloMichel2021git}.
\begin{table}[H] 
\small  
 \begin{center} 
		\begin{tabular}{| c | c || c | c | c || c | c | c || c | c | c | }  
	     \hline 
	     \multicolumn{2}{|c||}{Element $\&$ }  & \multicolumn{3}{|c||}{SUPG}  & \multicolumn{3}{|c||}{OSS}  & \multicolumn{3}{|c|}{CIP}  \\ \hline 
	     \multicolumn{2}{|c||}{ Time scheme }  & $\mathbb{P}_1$ & $\mathbb{P}_2$ & $\mathbb{P}_3$   & $\mathbb{P}_1$ & $\mathbb{P}_2$ & $\mathbb{P}_3$   & $\mathbb{P}_1$ & $\mathbb{P}_2$ & $\mathbb{P}_3$  \\ \hline \hline 
Basic              &  SSPRK & 1.93 & 2.96 & 4.02 & 2.0 & 2.62 & 4.1 & 1.44 & 2.45 & 3.77 \\ 
         \hline 
       \multirow{2}{*}{ Cub.}                &  SSPRK & 1.97 & 2.39 & 4.38 & 2.03 & 2.49 & 4.41 & 1.96 & 2.35 &  /  \\ 
               &  DeC & 1.97 & 2.27 & 4.34 & 2.02 & 2.49 & 4.41 & 2.01 & 2.35 &  /  \\ 
         \hline 
Bern.              &  DeC & 1.97 & 2.61 & 1.8 & 2.29 & 2.52 & 2.27 & 1.97 & 2.7 & 2.06 \\ 
         \hline 
        \end{tabular} 
    \end{center} 
     \caption{Convergence order for all schemes on linear advection test, using coefficients obtained in \cref{tab:new_param_meshX_LinearAdvection-2D-RES}. \newline 
     ``/" means that the Fourier analysis showed that the scheme is unstable.} \label{tab:conv_order_meshX_LinearAdvection-2D-RES}
\end{table}%

Looking at the table \ref{tab:conv_order_meshX_LinearAdvection-2D-RES}, we observe that almost all the stabilized schemes provide the expected order of accuracy. Exception to this rule are several $\P_2$ discretization which reach an order of accuracy of $\approx 2.5$, and all \textit{Bernstein} $\B_3$ polynomials with the \textit{DeC} which reach an order of accuracy of $2$. This result is very disappointing and it does not improve even adding more corrections, as suggested in \cite{paola_svetlana,DeC_2017}. Moreover, it has been independently verified that also in Fourier space the accuracy of DeC with Bernstein polynomials of degree 3 is only of order 2. This problem do not show up for steady problems, as there only the spatial discretization determines the order of accuracy. We will show it in \cref{sec:numerical_nlsw2D-stedyVSunsteady}, where we study also some steady vortexes. The authors still do not understand why the optimal order of accuracy is not reached. This opens doors to further research on this family of schemes. \\
Note that we do not show results for \textit{Bernstein} elements with \textit{SSPRK} technique because they are identical to \textit{Basic} elements, but are more expensive because of the projection in the \textit{Bernstein} element space and the interpolation in the quadrature points. \\
More comparisons on different grids (unstructured) will be done in \cref{sec:numericalSimulations}.


\subsection{Shallow water equations} \label{sec:nlsw2D_results-meshX}
We consider the non linear shallow water equations (no friction and constant topography):
\begin{equation}
    \left  \{
    \begin{array}{lll}
    	\partial_t h + \partial_x (hu) + \partial_y (hv) & = 0,  \qquad \quad & x\in \Omega = [0,2]\times[0,1], \\
    	\partial_t (hu) + \partial_x (hu^2 +g\frac{h^2}{2} ) + \partial_y (huv) & =0, & t \in [0,t_f]\\
    	\partial_t (hv) + \partial_x (huv) + \partial_y (hv^2 +g\frac{h^2}{2} )   & =0, & t_f =1s. \label{eq:num_test_SW2}
	\end{array}	
    \right .
\end{equation}
An analytical solution of this system is given by travelling vortexes \cite{ricchiuto2021analytical}. We use here a vortex with compact support and in $\mathcal{C}^6(\Omega)$ described by 
%
%
\begin{align}\label{eq_sol2_SW2}
	\begin{pmatrix}
		h(x,t)\\u(x,t)\\v(x,t)
	\end{pmatrix}=
\begin{cases}
\begin{pmatrix}
	h_c + \frac{1}{g} \frac{\Gamma^2}{\omega^2} \cdot \left( \lambda(\omega  \mathcal{R}( \mathbf{x},t) ) - \lambda (\pi) \right) ,  \\
	u_c + \Gamma(1+\cos (\omega \mathcal{R}( \mathbf{x},t)))^2 \cdot (- \mathcal{I}(\mathbf{x},t)_y),  \\
	v_c + \Gamma(1+\cos (\omega \mathcal{R}( \mathbf{x},t)))^2 \cdot ( \mathcal{I}(\mathbf{x},t)_x),
\end{pmatrix}, &\mbox{if } \omega  \mathcal{R}( \mathbf{x},t) \leq \pi,\\[10pt]
\begin{pmatrix}
	h_c & u_c & v_c
\end{pmatrix}^T, &\mbox{else,}
\end{cases}
\end{align}    
with 
\begin{align*}
\lambda(r) = &\frac{20\cos(r)}{3} + \frac{27\cos(r)^2}{16} + \frac{4\cos(r)^3}{9} + \frac{\cos(r)^4}{16} + \frac{20r\sin(r)}{3} \\
              &+ \frac{35r^2}{16} + \frac{27r\cos(r)\sin(r)}{8} + \frac{4r\cos(r)^2 \sin(r)}{3} + \frac{r\cos(r)^3 \sin(r)}{4}.
\end{align*}
where $\mathbf{X_c} = (0.5,0.5)$ is the initial vortex center,  $(h_c,\, u_c,\, v_c)=(1.,\,0.6,\,0)$ is the far field state, $ r_0 = 0.45$ is the vortex radius, 
$\Delta h = 0.1$ is the vortex amplitude, and the remaining paramters are defined as 
\begin{equation}
    \left  \{
    \begin{array}{ll}
        \omega = \pi / r_0  \qquad & \mbox{ angular wave frequency}, \\
        \Gamma = \frac{12 \pi \sqrt{g \Delta h } }{r_0 \sqrt{315 \pi^2-2048}}  \qquad & \mbox{ vortex intensity parameter}, \\
		\mathcal{I}(\mathbf{x},t) = \mathbf{x} - \mathbf{X_c} - (u_c t,v_c t)^T \qquad & \mbox{ coordinates with respect to the vortex center}, \\
		\mathcal{R}( \mathbf{x},t) = \|   \mathcal{I}(\mathbf{x},t) \|  	\qquad & \mbox{ distance from the vortex center}.  \label{eq_param_SW2}
	\end{array}
    \right .
\end{equation}

\begin{figure}[h!]
	\begin{center}
		\includegraphics[height=0.18\textheight,trim={0 0 64mm 0}, clip]{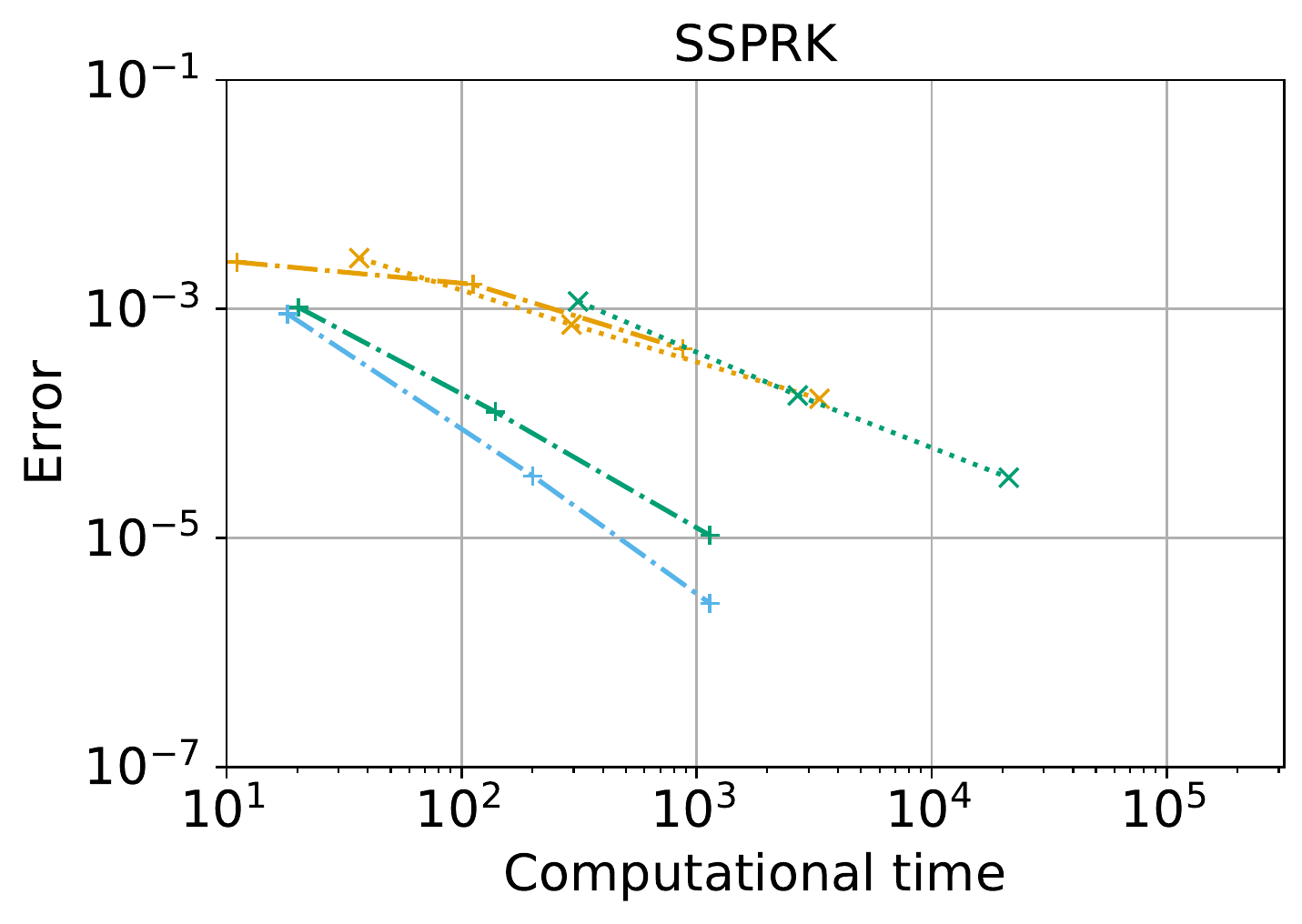} \;\;
		\includegraphics[height=0.18\textheight,trim={8mm 0 0 0}, clip]{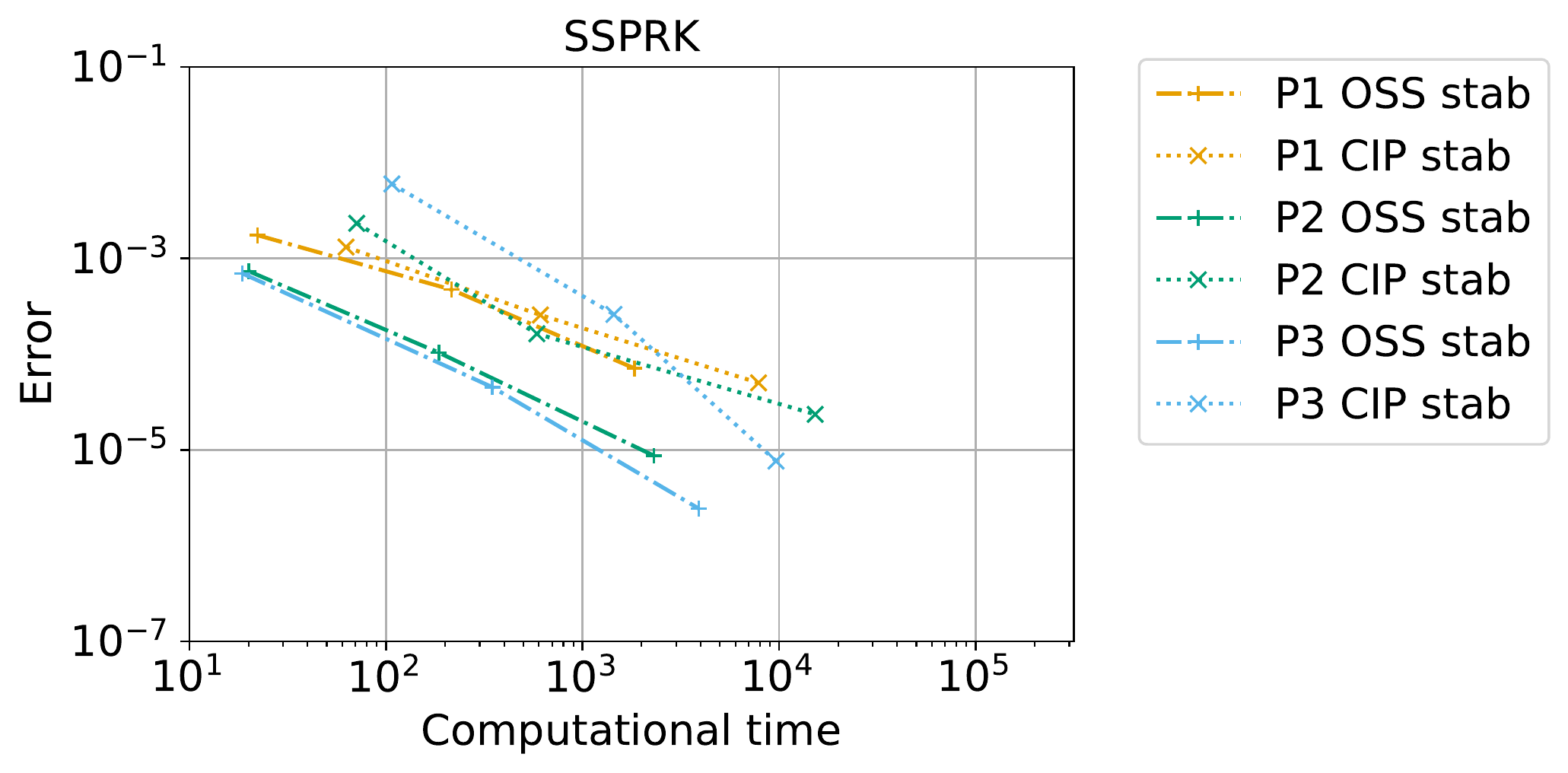}\\\vspace{2mm}
	\end{center}
	\caption{Error for shallow water system \eqref{eq:num_test_SW2} with respect to computational time for SSPRK method with \textit{Cubature} (left) and \textit{Basic} (right) elements and CIP and OSS stabilizations.}
	\label{fig:2D_timeVsErrorSW_meshX}
\end{figure}

We discretize the mesh with uniform square intervals of length $\Delta x$  (see figure \ref{fig:2d_meshX}), and as before we perform a grid convergence by respecting the constraint $\Delta x_2 = 2\Delta x_1 $ for $\mathbb{P}_2$ elements and $\Delta x_3 = 3 \Delta x_1 $ for $\mathbb{P}_3$ elements. Because of the high cost of the SUPG technique, we only compare the OSS and the CIP stabilization techniques. 
As an example of results, we again show the benefit of using \textit{Cubature} elements in \ref{fig:2D_timeVsErrorSW_meshX}. We can see that since the dimension of the discretized system is even larger than before (three times larger), the differences between \textit{Cubature} and \textit{Basic} elements are even more highlighted in the error-computational time plot.

In \cref{tab:conv_order_meshX_nlsw-2D-RES} we show the convergence orders for this shallow water problem with the CFL and $\delta$ coefficients found in \cref{tab:new_param_meshX_LinearAdvection-2D-RES}.
\begin{table}[H] 
\small  
 \begin{center} 
		\begin{tabular}{| c | c || c | c | c || c | c | c | }  
	     \hline 
	     \multicolumn{2}{|c||}{Element $\&$ }  & \multicolumn{3}{|c||}{OSS}  & \multicolumn{3}{|c|}{CIP}  \\ \hline 
	     \multicolumn{2}{|c||}{ Time scheme }  & $\mathbb{P}_1$ & $\mathbb{P}_2$ & $\mathbb{P}_3$   & $\mathbb{P}_1$ & $\mathbb{P}_2$ & $\mathbb{P}_3$  \\ \hline \hline 
Basic              &  SSPRK & 2.3 & 3.18 & 3.8 & 2.34 & 3.3 & 4.47 \\ 
         \hline 
       \multirow{2}{*}{ Cub.}               &  SSPRK & 1.25 & 3.31 & 3.94 & 2.03 & 2.56 &  /  \\ 
               &  DeC & 1.45 & 3.31 & 3.94 & 1.98 & 2.56 &  /  \\ 
         \hline 
Bern.              &  DeC & 1.52 & 2.93 & 2.97 & 2.92 & 2.12 & 2.91 \\ 
         \hline 
        \end{tabular} 
    \end{center} 
     \caption{Convergence order on shallow water, using coefficients obtained in \cref{tab:new_param_meshX_LinearAdvection-2D-RES}. \newline
     \textit{"/"} means that the fourier analysis shown that the scheme is unstable.} \label{tab:conv_order_meshX_nlsw-2D-RES}
\end{table}%

The results obtained are similar to those of the \textit{linear advection} case. 
We can also notice the $\P_2$ discretization reaching the proper convergence order, i.e., $3$, and \textit{Bernstein} $\B_3$ elements reaching an order of accuracy of $\approx 3$ which is more satisfying than the results obtained for the linear advection test, but still disappointing knowing that we were expecting $4$.

\section{Simulations on unstructured meshes}\label{sec:numericalSimulations}
\begin{figure}
	\centering
	\includegraphics[width=0.49\linewidth , height=4cm]{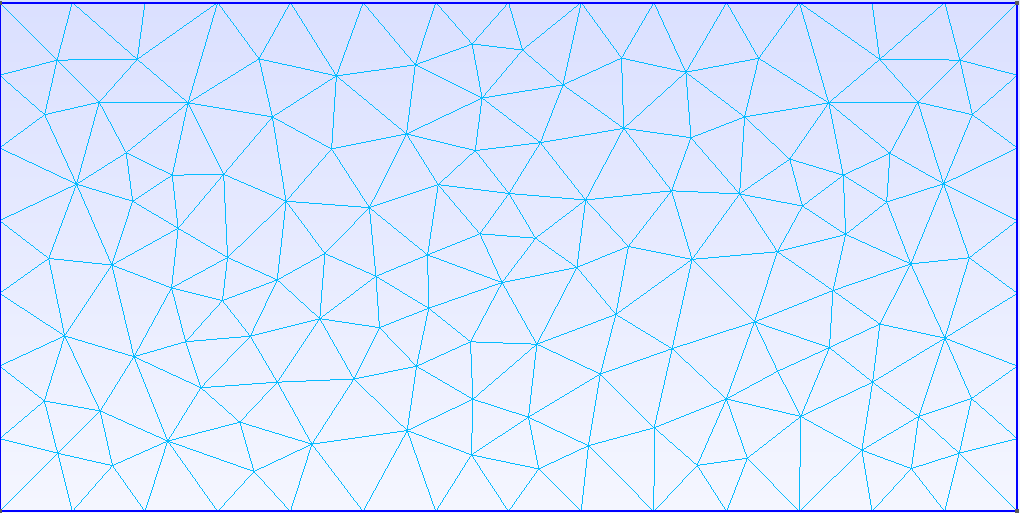}
		\caption{Unstructured mesh on $\Omega=[0,2]\times [0,1]$.}
	\label{fig:2d_mesh}
\end{figure}

We now perform numerical tests to check the validity of our theoretical findings using an unstructured mesh, and the most restrictive parameters in \cref{tab:restrictive_param_LinearAdvection-2D-RES}. These parameters make sure that we are stable for both \textit{T} and \textit{X} mesh configurations. The results have similar convergence rate to the tests on the structured meshes of the previous section. \\
The unstructured mesh used in this section is shown in \cref{fig:2d_mesh}, and it was created by the mesh generator \textit{gmsh}\footnote{https://gmsh.info/}.

\subsection{Linear advection test} \label{sec:LA2D_results-unstructure_mesh}
\begin{table}[H] 
\small  
 \begin{center} 
		\begin{tabular}{| c | c || c | c | c || c | c | c || c | c | c | }  
	     \hline 
	     \multicolumn{2}{|c||}{Element $\&$ }  & \multicolumn{3}{|c||}{SUPG}  & \multicolumn{3}{|c||}{OSS}  & \multicolumn{3}{|c|}{CIP}  \\ \hline 
	     \multicolumn{2}{|c||}{ Time scheme }  & $\mathbb{P}_1$ & $\mathbb{P}_2$ & $\mathbb{P}_3$   & $\mathbb{P}_1$ & $\mathbb{P}_2$ & $\mathbb{P}_3$   & $\mathbb{P}_1$ & $\mathbb{P}_2$ & $\mathbb{P}_3$  \\ \hline \hline 
Basic              &  SSPRK & 1.9 & 2.57 & 3.76 & 1.99 & 2.5 & 3.76 & 1.57 & 2.14 & 3.66 \\ 
         \hline 
       \multirow{2}{*}{\centering  Cub.}               &  SSPRK & 1.73 & 2.4 & 3.83 & 1.81 & 2.53 &  3.98$^{**}$  & 1.8 & 2.17 &  /  \\ 
               &  DeC & 1.81 & 2.21 & 2.56 & 1.82 & 2.48 &   3.98$^{**}$  & 1.83 & 2.17 &  /  \\ 
         \hline 
Bern.              &  DeC & 1.78 & 2.12 & 1.94 & 2.31 & 2.48 & 2.12 & 1.56 & 2.03 & 2.24 \\ 
         \hline 
        \end{tabular} 
     \caption{Convergence order for linear advection on unstructured mesh, using coefficients obtained in \cref{tab:restrictive_param_LinearAdvection-2D-RES}.\newline
     $^{**}$ These values are found using only the $X$ mesh  (see \cref{fig:fourier_combined_cubP3_problem}). \newline
     \textit{"/"} means that the scheme is clearly unstable.} \label{tab:conv_order_mesh_combined_LinearAdvection-2D-RES}    
    \end{center} 
\end{table}%

We use the same test case of \cref{sec:LA2D_results-meshX}. 
Convergence orders for all schemes are summarized in \cref{tab:conv_order_mesh_combined_LinearAdvection-2D-RES}. 
We observe that all $\P_1$ discretizations provide the proper convergence order. For $\P_2$ discretization we spot a slight reduction of the order of accuracy, which lays for most of the schemes between $2$ and $\approx 2.5$ instead of being $3$. For polynomials of degree $3$, we observe an order reduction to 2 for the same schemes that lost the right order of accuracy also for \textit{X} mesh in the previous section. 
In particular, we have that  \textit{Bernstein} $\B_3$ polynomials with the \textit{DeC} result in an order of accuracy of $\approx 2$ instead of $4$, as well as the $\TP_3$ discretization with the combination DeC and SUPG stabilization. 
As for the X mesh, the \textit{Basic} $\P_3$ discretization reach order of accuracy $\approx 4$ for all stabilization techniques, as well as \textit{Cubature} $\TP_3$ with SUPG and OSS stabilizations. \\
Also in this case, the results obtained with $\TP_3$ \textit{Cubature} elements and OSS stabilization are stable as we can see from the convergence analysis. This might mean that just few unfortunate mesh configurations, as the $T$ one, result in an unstable scheme and that, most of the time, the parameters found in \cref{tab:restrictive_param_LinearAdvection-2D-RES}  are reliable for this scheme.
On the other hand, the combination $\TP_3$ and CIP gives an unstable scheme. 

\begin{figure}[h!]
	\begin{center}
		\subfigure[\textit{Cubature} elements with DeC]{
			\includegraphics[height=0.18\textheight,trim={0 0 64mm 0}, clip]{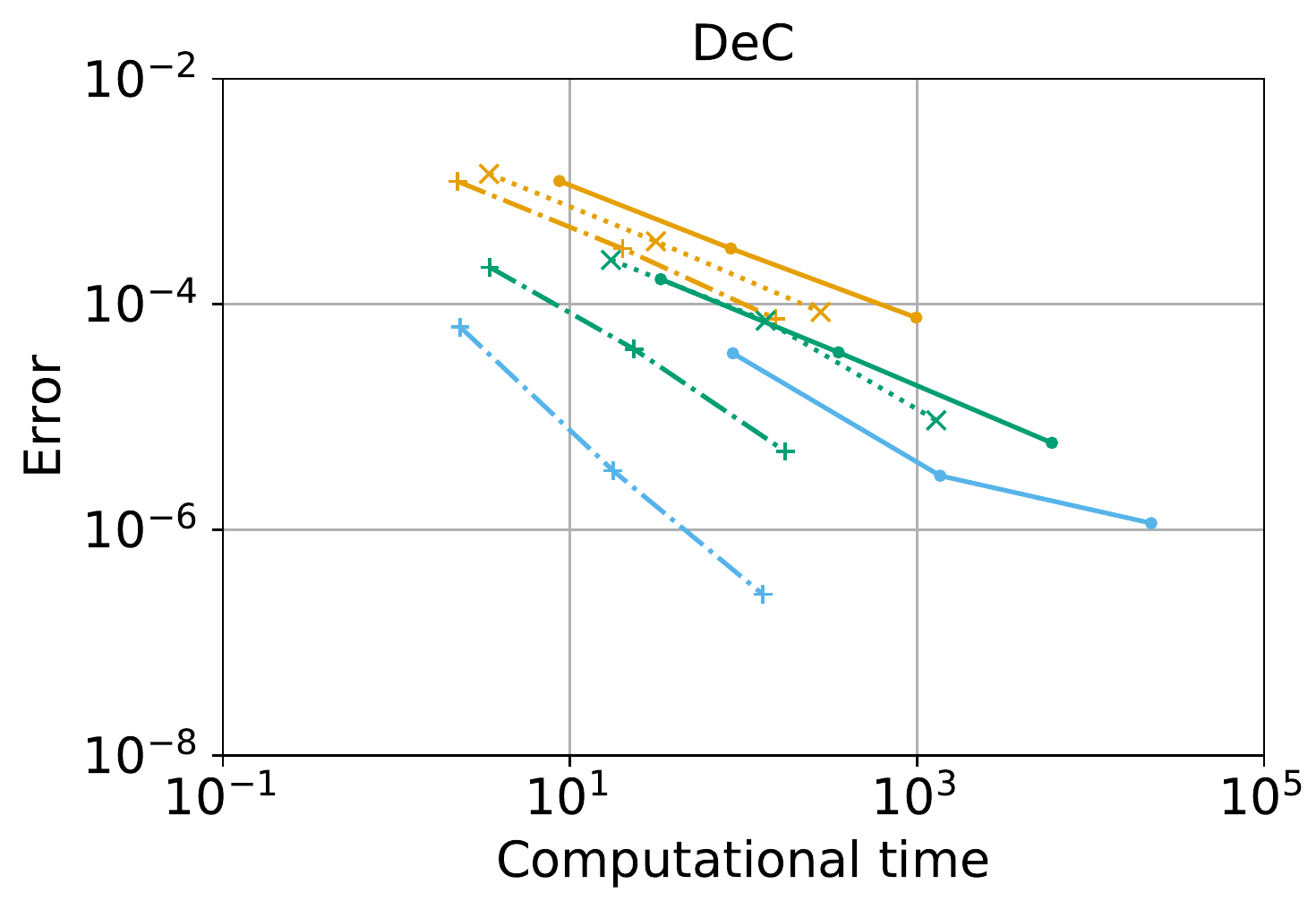}} \qquad
		\subfigure[\textit{Cubature} elements with SSPRK]{
			\includegraphics[height=0.18\textheight,trim={10mm 0 0 0}, clip]{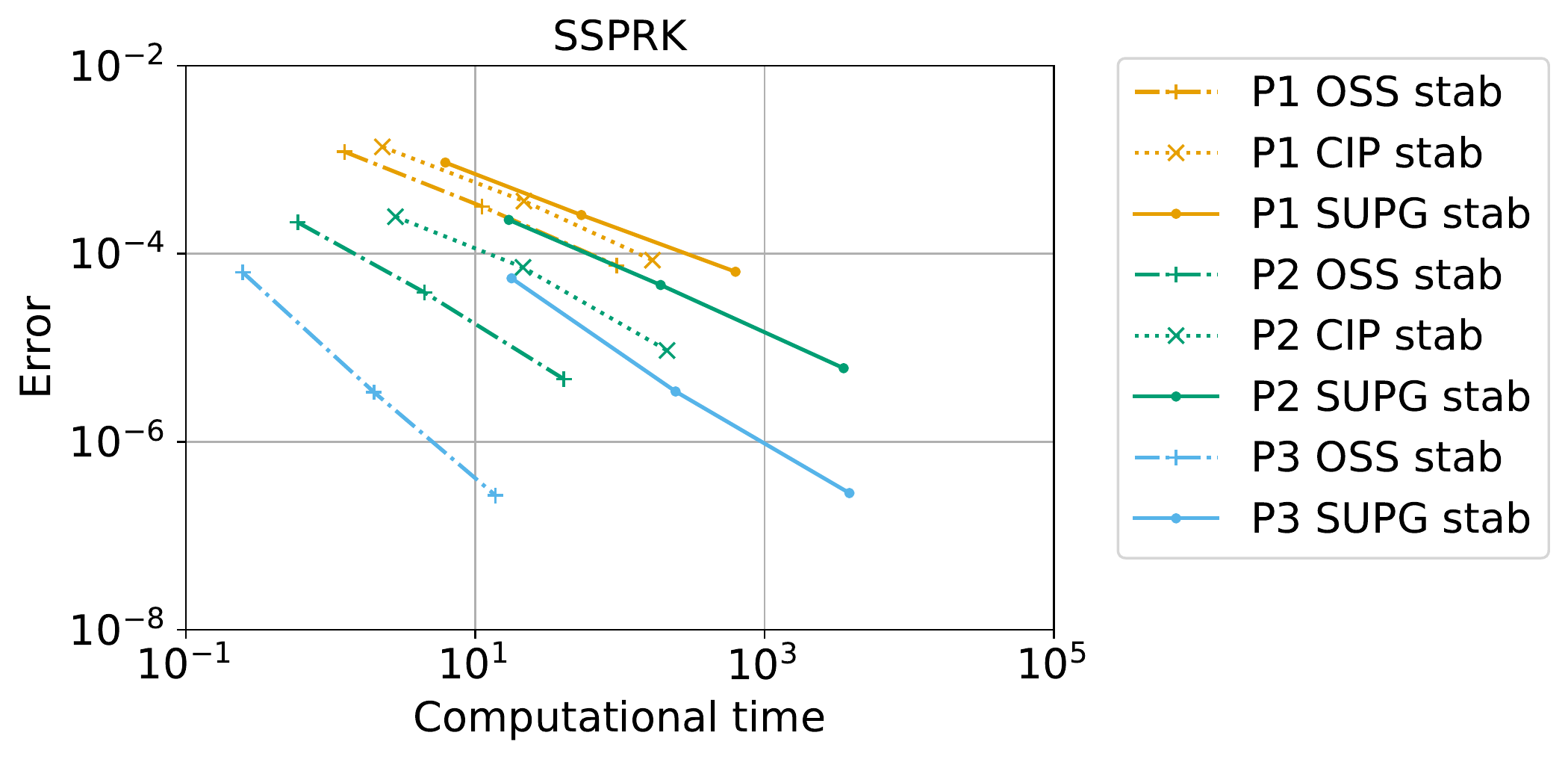}}\\
		\subfigure[\textit{Bernstein} elements with DeC]{		\includegraphics[height=0.18\textheight,trim={0 0 64mm 0}, clip]{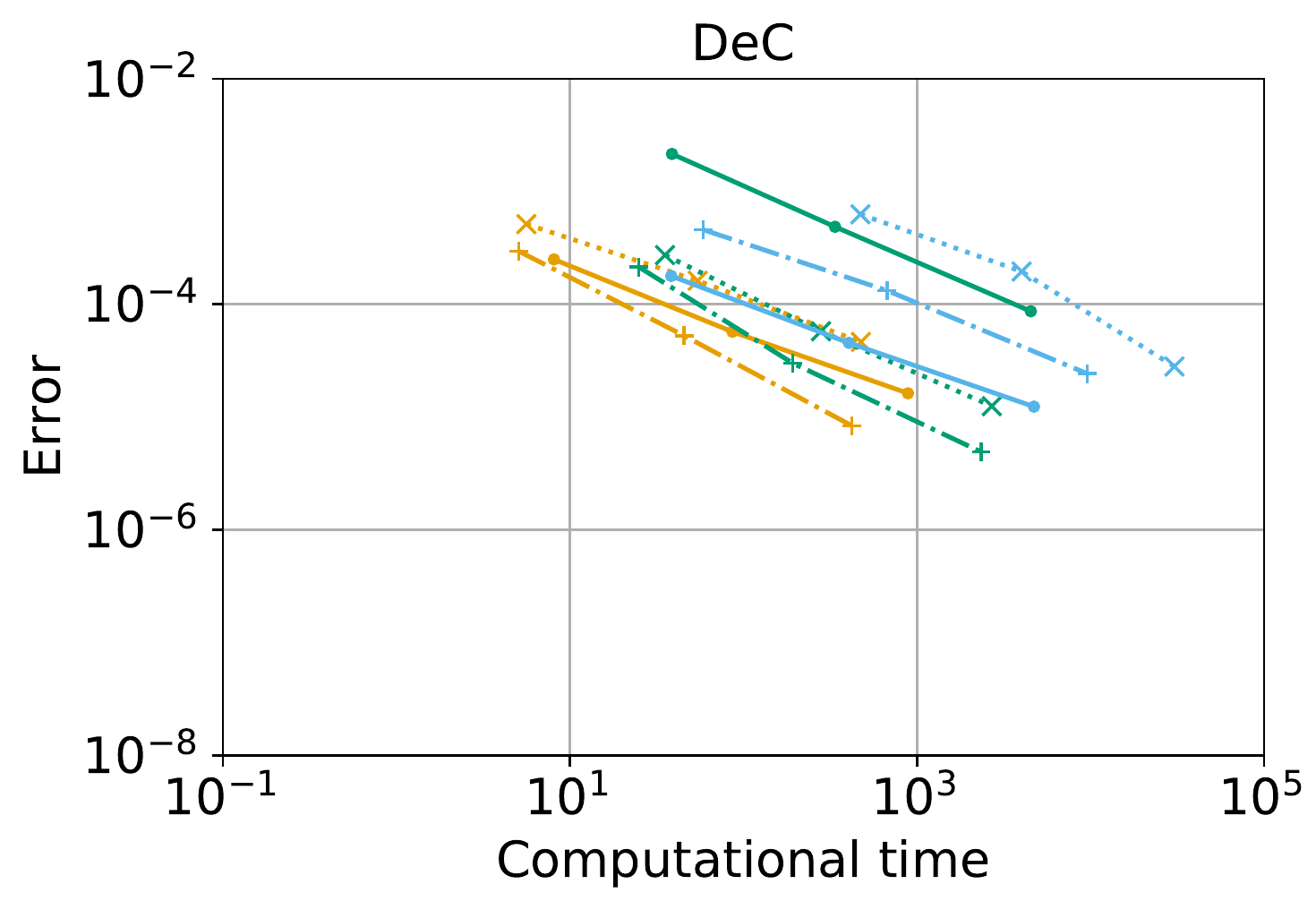}} \qquad
		\subfigure[\textit{Basic} elements with SSPRK]	{\includegraphics[height=0.18\textheight,trim={10mm 0 0 0}, clip]{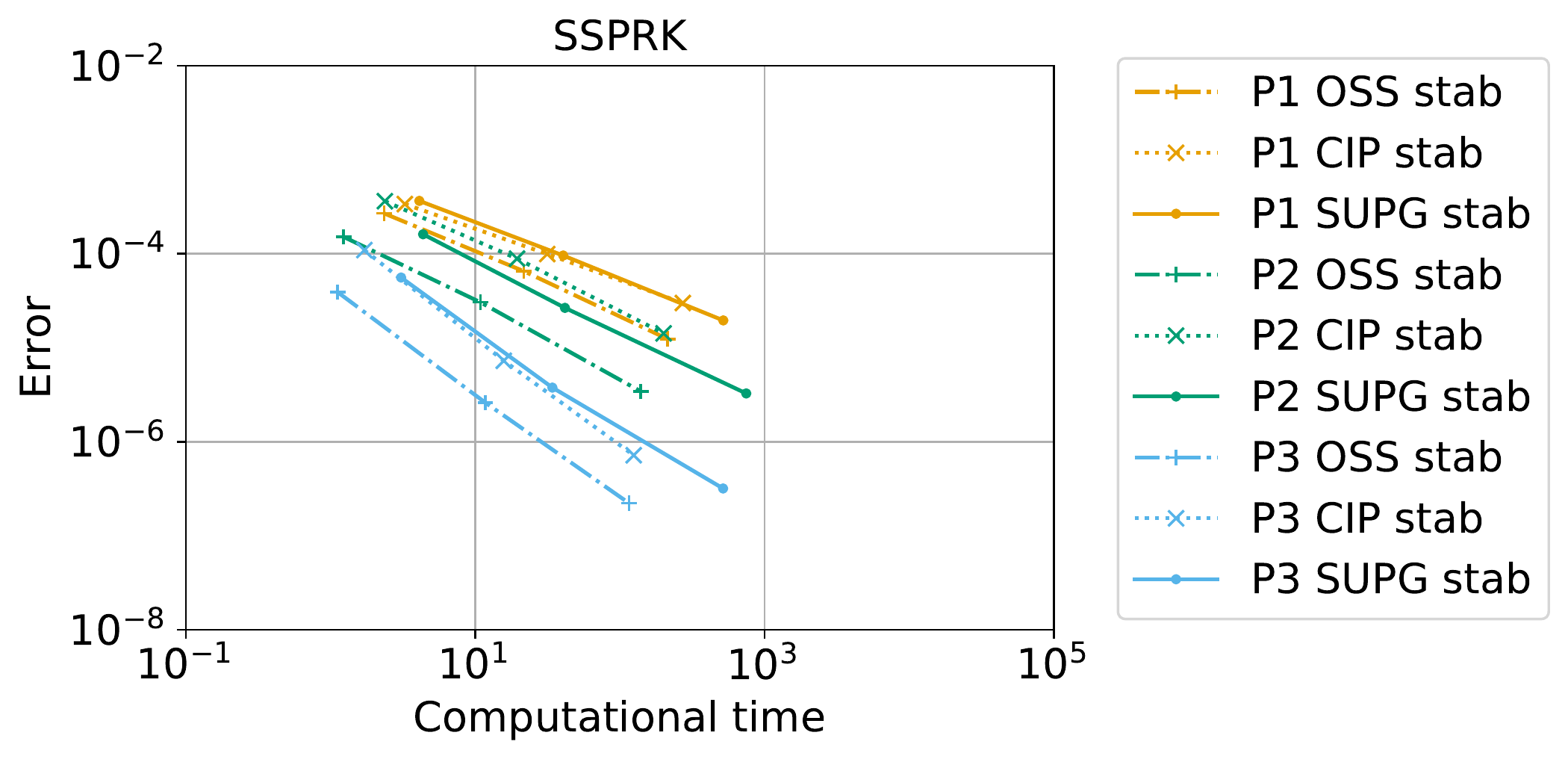}}
	\end{center}
	\caption{Error for linear advection problem \eqref{lin_adv_2D} with respect to computational time for all elements and stabilization techniques}
	\label{fig:2D_timeVsErrorLinAdv_meshcombined}
\end{figure}

We compare error and computational time for all methods presented above in \cref{fig:2D_timeVsErrorLinAdv_meshcombined}.
Looking at $\P_2$ and the $\P_3$ discretizations, as expected, the mass-matrix free combination, i.e., \textit{Cubature} elements with  SSPRK and OSS, gives smaller computational costs than other combinations with \textit{Basic} elements. Conversely, the SUPG technique increase the computational costs with respect to all other stabilizations for all schemes. That is why we will not use it for the next test.
The plots and all the errors are available at the repository \cite{TorloMichel2021git}.

\begin{remark}[Entropy viscosity]
	
	As remarked in \cref{sec:app_fourier_visco}, we can improve the stability of some schemes (\textit{Cubature} OSS) with extra entropy viscosity. Here,	we test the convergence rate on the \textit{T} mesh configuration, i.e., the one with more restrictive CFL conditions and most unstable. This test is performed using \textit{Cubature} $\TP_3$ elements, SSPRK and DeC time integration methods, and the OSS and the CIP stabilization techniques. We solve again problem \eqref{lin_adv_2D}. 
	
	Using formulation \eqref{eq:FV-visco_app}
	and tuning stability coefficient $\delta$, CFL and viscosity coefficient $c$ found in \cref{fig:fourier_cubP3_SSPRK_OSS_viso}, we obtain fourth order accurate schemes. These tuned coefficients, and the corresponding convergence orders are summarized in \cref{tab:conv_order_meshTvisco_P3cub_LinearAdvection-2D-def_init}.
	\begin{table}[H] 
		\small  
		\begin{center} 
			\begin{tabular}{| c | c || c | c | c || c | c | c | }  
				\hline 
				\multicolumn{2}{|c||}{Element $\&$ }  & \multicolumn{3}{|c||}{\textit{Cubature} $\TP_3$ OSS}  & \multicolumn{3}{|c|}{\textit{Cubature} $\TP_3$  CIP}  \\ \hline 
				\multicolumn{2}{|c||}{ Time scheme }  & \CFL  \,($\delta$) & $c$ & order  & \CFL\,($\delta$)  & $c$ & order  \\ \hline \hline 
				\multirow{2}{*}{ Cub.}                &  SSPRK & 0.15  (0.02) &  0.05  &  4.08 & 0.12  (0.0004) &  0.5  &  3.60  \\ 
				&  DeC  & 0.15  (0.02) &  0.05  &  4.09 &0.08  (0.001) &  0.2  &  3.76   \\ \hline
			\end{tabular} 
		\end{center} 
		\caption{Convergence order of methods using \textit{Cubature} $\TP_3$ elements and viscosity term \eqref{eq:FV-visco_app} with tuned parameters} \label{tab:conv_order_meshTvisco_P3cub_LinearAdvection-2D-def_init}
	\end{table}%
	
	Many other formulations of viscosity terms exist in literature and can ensure convergent methods of order $p+1$ (using $\P_p$ elements) \cite{2011JCoPh.230.4248G,KUZMIN2020104742,inbookLlobell2020}. The majority use a nonlinear evaluation of the parameter $\mu_K$, based on the local entropy production. \\
\end{remark}


\subsection{Shallow water equations}
\begin{table}[H] 
\small  
 \begin{center} 
		\begin{tabular}{| c | c || c | c | c || c | c | c | }  
	     \hline 
	     \multicolumn{2}{|c||}{Element $\&$ }  & \multicolumn{3}{|c||}{OSS}  & \multicolumn{3}{|c|}{CIP}  \\ \hline 
	     \multicolumn{2}{|c||}{ Time scheme }  & $\mathbb{P}_1$ & $\mathbb{P}_2$ & $\mathbb{P}_3$   & $\mathbb{P}_1$ & $\mathbb{P}_2$ & $\mathbb{P}_3$  \\ \hline \hline 
Basic              &  SSPRK & 1.94 & 2.98 & 4.25 & 2.15 & 2.52 & 4.11 \\ 
         \hline 
       \multirow{2}{*}{\centering  Cub.}               &  SSPRK & 1.03 & 3.17 & 3.59$^{**}$ & 1.39 & 2.57 &  /  \\ 
               &  DeC & 1.2 & 3.14 & 3.59$^{**}$ & 1.48 & 2.57 &  /  \\ 
         \hline 
Bern.              &  DeC & 1.28 & 3.14 &  3.15  & 1.36 & 2.73 &  2.66  \\ 
         \hline 
        \end{tabular} 
    \end{center} 
     \caption{Convergence order on shallow water for unstructured mesh, using coefficients obtained in \cref{tab:restrictive_param_LinearAdvection-2D-RES}.\newline
     $^{**}$ These values are found using only the $X$ mesh (see \cref{fig:fourier_combined_cubP3_problem}). \newline
     \textit{"/"} means that the scheme is clearly unstable.} \label{tab:conv_order_mesh_combined_nlsw-2D-RES}
\end{table}%

In this section we test the proposed schemes on the test case of \cref{sec:nlsw2D_results-meshX} with the unstructured mesh in \cref{fig:2d_mesh}.
Convergence orders are summarized in \cref{tab:conv_order_mesh_combined_nlsw-2D-RES}.
\begin{figure}[h!]
	\begin{center}
		\subfigure[\textit{Cubature} elements with DeC]{
			\includegraphics[height=0.18\textheight,trim={0 0 64mm 0}, clip]{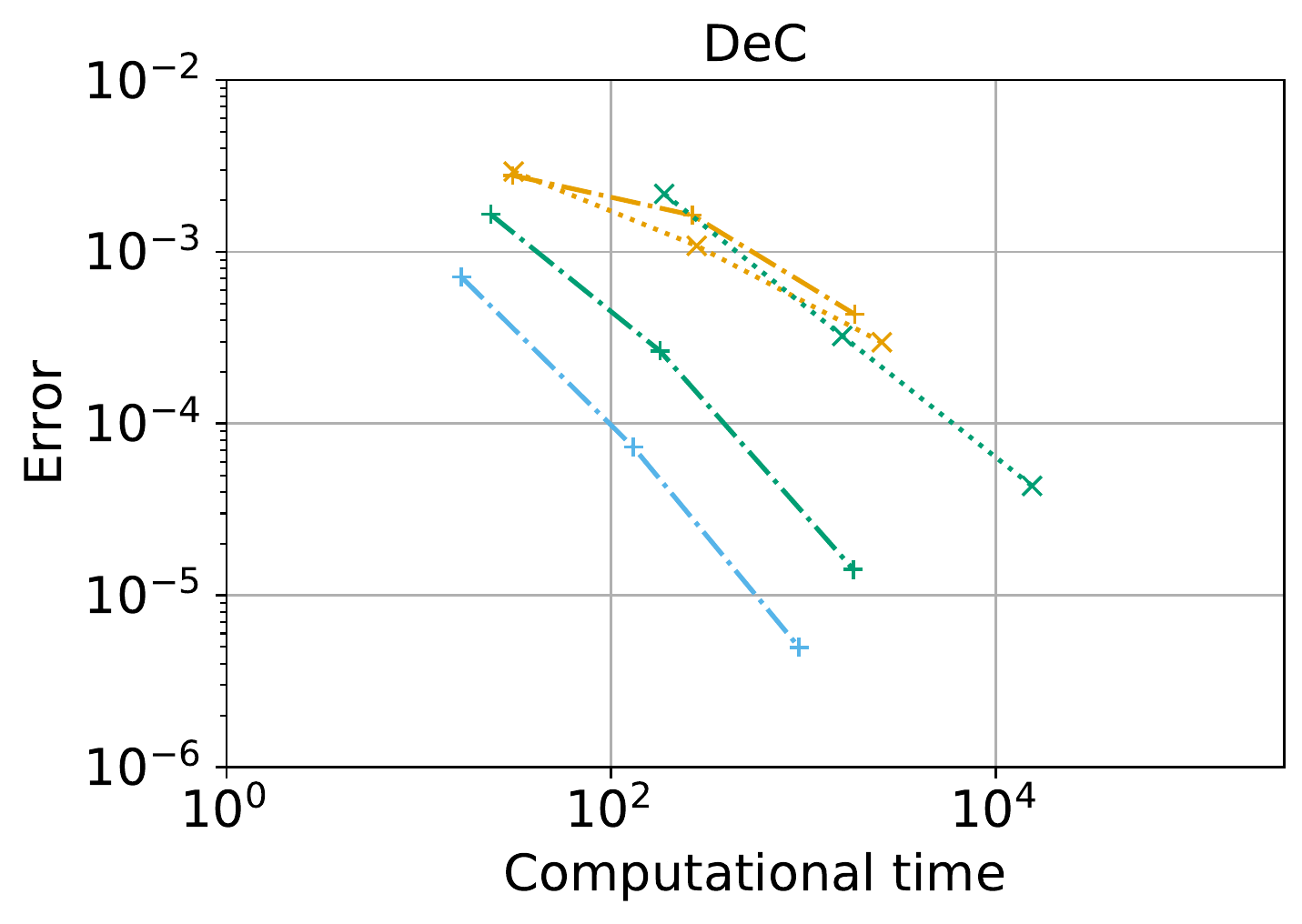}} \qquad
		\subfigure[\textit{Cubature} elements with SSPRK]{
			\includegraphics[height=0.18\textheight,trim={10mm 0 0 0}, clip]{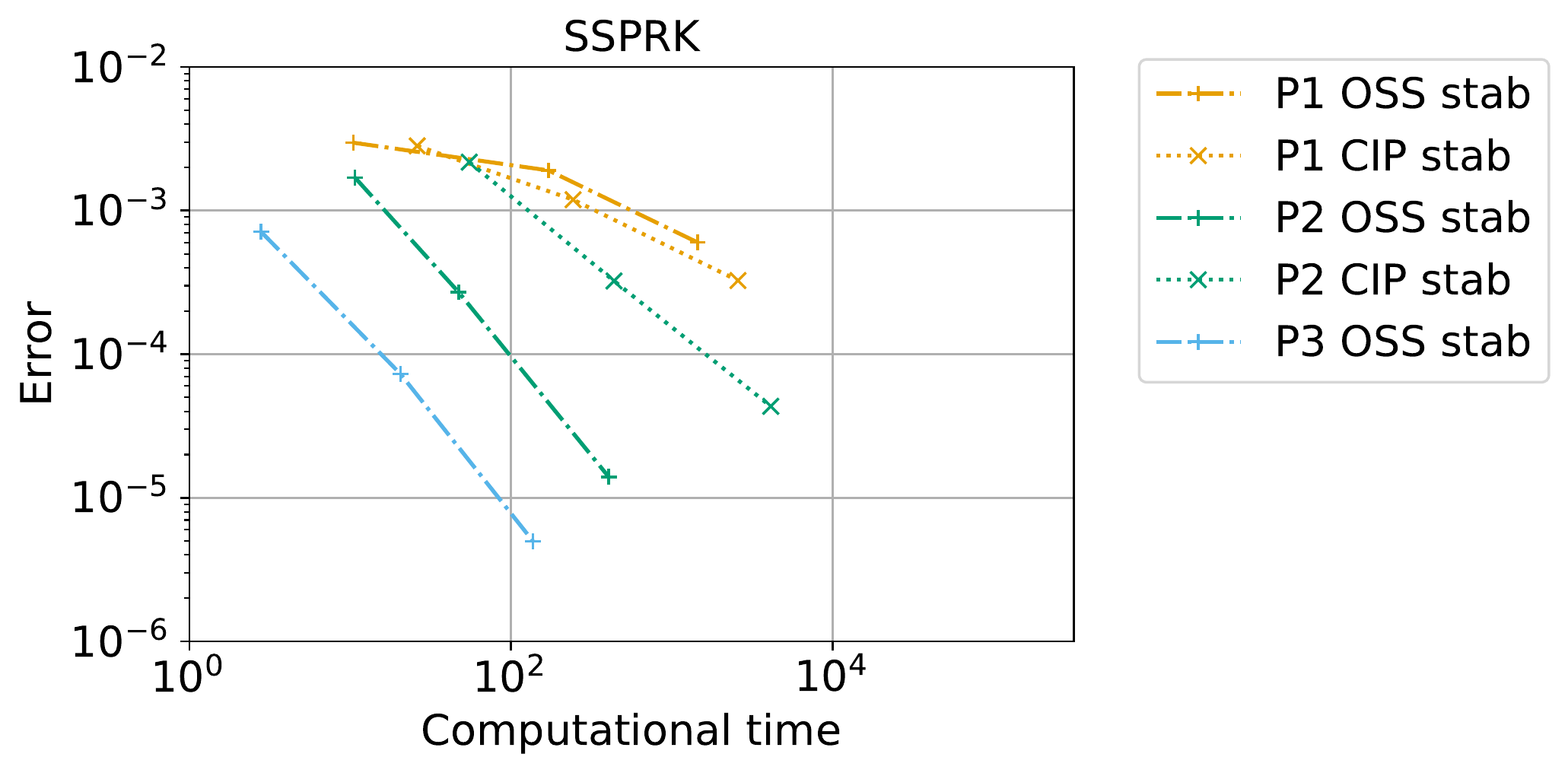}}\\
		\subfigure[\textit{Bernstein} elements with DeC]{
			\includegraphics[height=0.18\textheight,trim={0 0 64mm 0}, clip]{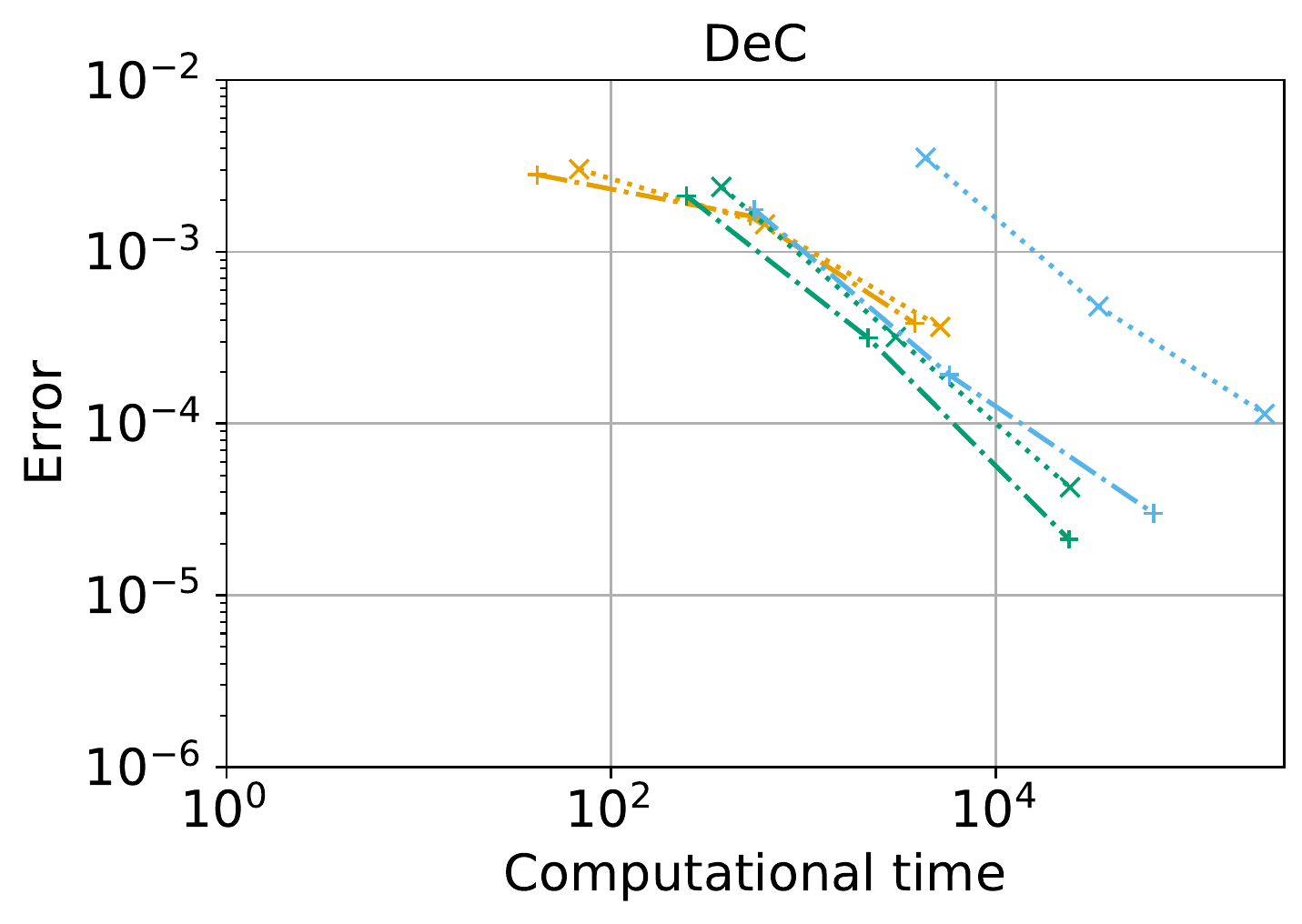}} \qquad
		\subfigure[\textit{Basic} elements with SSPRK]	{
			\includegraphics[height=0.18\textheight,trim={10mm 0 0 0}, clip]{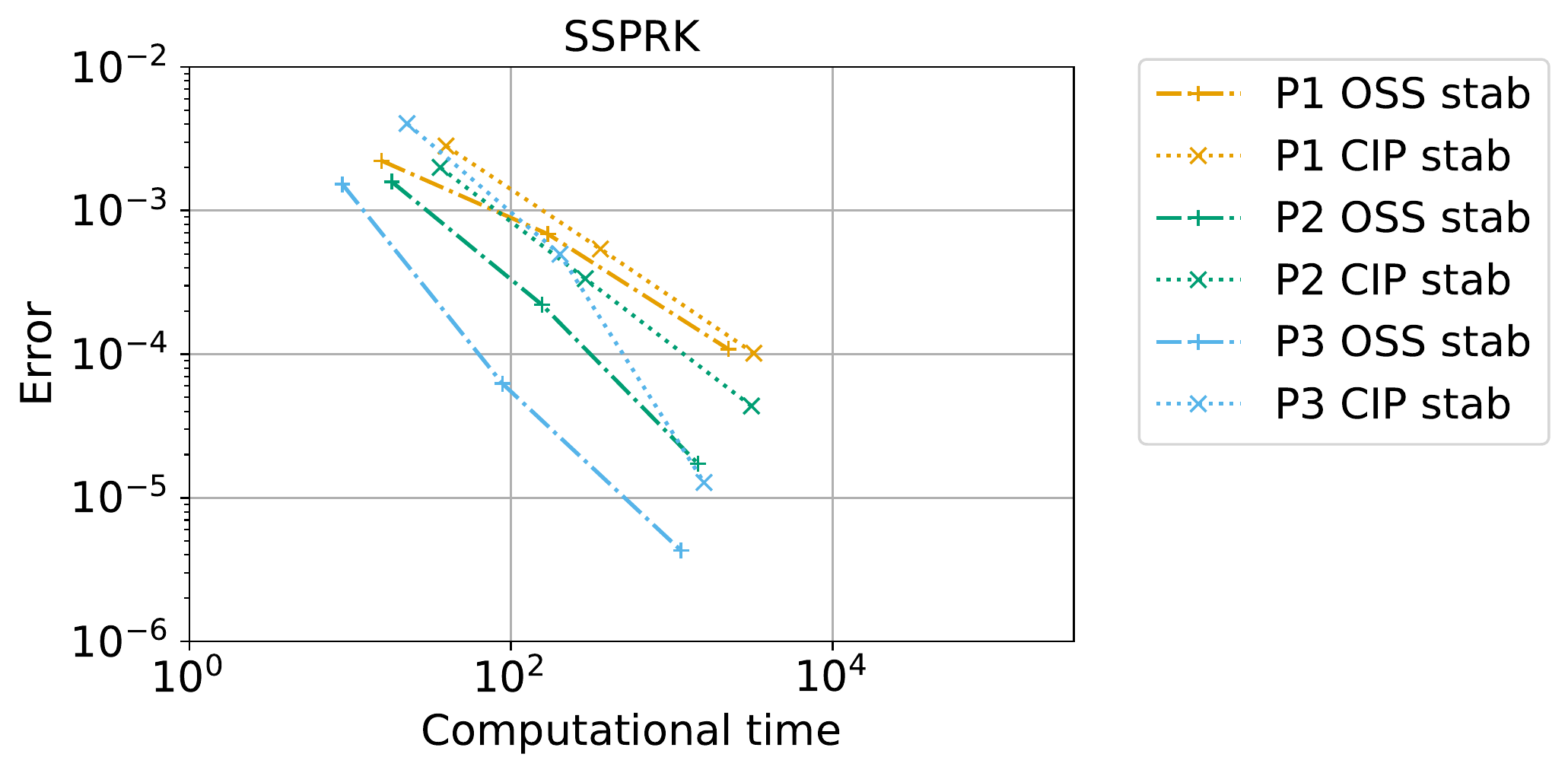}}
	\end{center}
	\caption{Error for shallow water problem \eqref{eq:num_test_SW2} with respect to computational time for all elements and stabilization techniques}
	\label{fig:2D_timeVsErrorSW}
\end{figure}
Also for the shallow water equations, we have results that resemble the ones of the structured mesh. There are small differences in the order of accuracy in both directions in different schemes. 
Comparing also the computational time of all the schemes in \cref{fig:2D_timeVsErrorSW}, we can choose what we consider the best numerical method for these test cases: \textit{Cubature} discretization with the OSS stabilization technique. 
This performance seems fully provided by the free mass-matrix inversion, as the CFLs for the OSS technique (with SSPRK scheme) is approximately the same between \textit{Basic} and \textit{Cubature} elements (see  \cref{tab:restrictive_param_LinearAdvection-2D-RES}). 

The plots and all the errors are available at the repository \cite{TorloMichel2021git}.



\subsection{Remark on the steady vortex case} \label{sec:numerical_nlsw2D-stedyVSunsteady}

For completeness we consider now a steady vortex, similarly to what reported in  \cite{paola_svetlana} for the isentropic Euler equations.
So, we consider again
the traveling vortex proposed in \cref{sec:nlsw2D_results-meshX} with $t_f =0.1s$. We compare the convergence orders between $u_c=0$ (steady case) and $u_c=0.6$ (unsteady case) in 
\cref{tab:conv_order_nlsw_2D-steady} and \cref{tab:conv_order_nlsw_2D-unsteady}. As we can see, in the steady case we obtain, \emph{without any additional viscous stabilization},   the expected convergence order for all schemes, in particular for the DeC with Bernstein polynomial function. These results agree with the ones in \cite{paola_svetlana}. Comparing with the unsteady case, all the other schemes reach similar order of accuracy as obtained in \cref{tab:conv_order_mesh_combined_nlsw-2D-RES}. Running the test with additional corrections in DeC scheme, as often suggested in \cite{paola_svetlana,DeC_2017}, does not improve the convergence order in the unsteady case (even with $K=50$). 

\begin{minipage}{0.54\textwidth}
\begin{table}[H] 
\small  
 \begin{center} 
		\begin{tabular}{| c | c || c | c | c || c | c | c | }  
	     \hline 
	     \multicolumn{2}{|c||}{Element $\&$ }  & \multicolumn{3}{|c||}{OSS}  & \multicolumn{3}{|c|}{CIP}  \\ \hline 
	     \multicolumn{2}{|c||}{ Time scheme }  & $\mathbb{P}_1$ & $\mathbb{P}_2$ & $\mathbb{P}_3$   & $\mathbb{P}_1$ & $\mathbb{P}_2$ & $\mathbb{P}_3$  \\ \hline \hline 
Basic              &  SSPRK & 2.31 & 2.67 & 3.89 & 1.97 & 2.64 & 3.62 \\ 
         \hline 
       \multirow{2}{*}{\centering  Cub.}              &  SSPRK & 2.05 & 3.2 & 3.56 & 1.79 & 2.83 &  /  \\ 
               &  DeC & 2.17 & 3.18 & 3.57 & 1.74 & 2.83 &  /  \\ 
         \hline 
Bern.              &  DeC & 2.33 & 3.28 & 3.65 & 1.85 & 3.0 & 3.63 \\ 
         \hline 
        \end{tabular} 
    \end{center} 
     \caption{Convergence order for steady vortex, $t_f=0.1s$. \newline
     \textit{"/"} means that the scheme is clearly unstable.} \label{tab:conv_order_nlsw_2D-steady}
\end{table}%
\vspace*{3mm}
\end{minipage}\hfill
\begin{minipage}{0.44\textwidth}
\begin{table}[H] 
\small  
 \begin{center} 
		\begin{tabular}{| c | c | c || c | c | c | }  
	     \hline 
	     \multicolumn{3}{|c||}{OSS}  & \multicolumn{3}{|c|}{CIP}  \\ \hline 
	      $\mathbb{P}_1$ & $\mathbb{P}_2$ & $\mathbb{P}_3$   & $\mathbb{P}_1$ & $\mathbb{P}_2$ & $\mathbb{P}_3$  \\ \hline \hline 
		2.34 & 2.68 & 3.86 & 1.94 & 2.53 & 3.61 \\ 
         \hline 
       2.03 & 3.13 & 3.57 & 1.74 & 2.7 &  /  \\ 
       2.13 & 3.09 & 3.57 & 1.71 & 2.7 &  /  \\ 
         \hline 
			2.33 & 3.19 & 2.87 & 1.75 & 2.77 & 2.76 \\ 
         \hline 
        \end{tabular} 
    \end{center} 
     \caption{Convergence order for unsteady vortex, $t_f=0.1s$. } \label{tab:conv_order_nlsw_2D-unsteady}
\end{table}%
\vspace*{3mm}
\end{minipage}
These results show that a numerical error appears in the spatio-temporal integration part of the solution \eqref{eq:L2}, which might be related to the fact that the high order derivatives are never penalized in our stabilizations and might produce some small oscillations.
%
\section{Conclusion}\label{sec:conclusion2D}

This work shows also that the stability results obtained in the one dimensional analysis \cite{michel2021spectral} can not be generalized for two dimensional problems on triangular meshes. In this direction, it could be interesting to perform the stability analysis on Cartesian quadrilateral meshes, to check whether in that situation the one dimensional results still hold true. \\
 
In the numerical test section, the order of accuracy found is not the expected one for all the methods, i.e., $p+1$ using $\P_p$ elements. For several cases, we reach only $p+1/2$ or $p$. Among the schemes that are stable and with the right order of accuracy, the method that uses \textit{Cubature} elements with OSS stabilization technique and SSPRK method of order 4 has proven to be the most accurate and less expensive. Secondly, comparing to the SUPG stabilization technique, very often used in the literature for hyperbolic system, we showed that other stabilization techniques such as CIP and OSS can provide the same accuracy and are cheaper in term of computational costs. \\

In this direction, it would be interesting to evaluate the stability of the CIP adding a additional penalty term on the jump of higher order derivatives as suggested in \cite{Burman2020ACutFEmethodForAModelPressure,burman2021weighted,paola_svetlana}. Moreover, it could be interesting to see the stability of \textit{Cubature} elements using higher degree polynomials. 
Another interesting point to explore is the loss of accuracy obtained using the DeC with \textit{Bernstein} third order polynomial basis functions for unsteady cases. \\

Finally, we provided a heuristic approach characterized by additional discontinuity capturing viscous operators
such as those proposed in \cite{2011JCoPh.230.4248G,KUZMIN2020104742}.  Even for smooth solutions, the very
small additional dissipation introduced by these terms is enough to stabilize some of the symmetric mass-matrix-free 
approaches, otherwise linearly unstable. This allows to obtain interesting schemes for practical purposes.

\section*{Acknowledgment}
This work was performed within the Ph.D. project of Sixtine Michel: ``Finite element method for shallow water equations'', supported by INRIA and the BRGM, co-funded by in INRIA--Bordeaux Sud--Ouest and the Conseil R\'egional de la Nouvelle Aquitaine.
Davide Torlo has been funded by a postdoctoral fellowship in the  team CARDAMOM in INRIA--Bordeaux Sud--Ouest and by a postdoctoral fellowship in SISSA.
R\'emi Abgrall has been supported by the Swiss National Foundation grant ``Solving advection dominated problems with high order schemes with polygonal meshes: application to compressible and incompressible flow problems'' under Grant Agreement No 200020\_175784. 

\section*{Declarations}
\textbf{Funding} 
Sixtine Michel was funded by in INRIA--Bordeaux Sud--Ouest and the Conseil R\'egional de la Nouvelle Aquitaine.
Davide Torlo has been funded by a postdoctoral fellowship in the  team CARDAMOM in INRIA--Bordeaux Sud--Ouest and by a postdoctoral fellowship in SISSA
R\'emi Abgrall has been supported by the Swiss National Foundation grant ``Solving advection dominated problems with high order schemes with polygonal meshes: application to compressible and incompressible flow problems'' under Grant Agreement No 200020\_175784. \\
\textbf{Conflicts of interest/Competing interests }
The authors certify that there is no actual or potential conflict of interest in relation to this article.\\
\textbf{Availability of data and material}
The images for all the parameters of the stability analysis and convergence plots are available at \cite{TorloMichel2021git}.\\

\appendix

\section{\textit{Cubature} elements, definition and construction} \label{sec:appendix_cohen_bf}
\begin{figure}
    \centering
    \hspace*{-1cm}
    \includegraphics[width=0.4\linewidth]{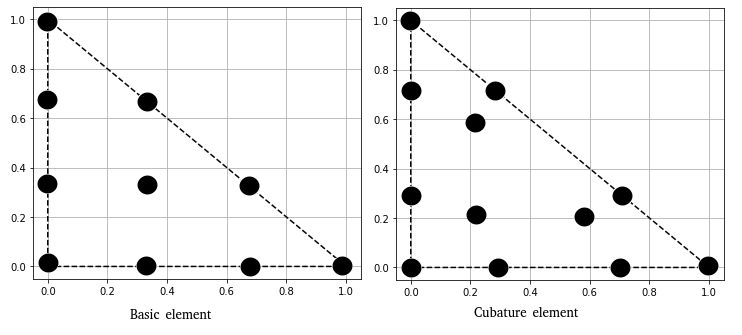}
    \caption{Comparison of two element of degree three: at left the classical one $\P_3$, at right the \textit{Cubature} one $\TP_3$.}
    \label{fig:cubature_mesh2}
\end{figure}

In this section we give a description of the \textit{Cubature} finite elements \cite{article_cubature_2006,article_cubature_2001}. In \cref{fig:cubature_mesh2} we show the $\TP_3$ example comparing the Lagrangian nodes of \textit{Basic} and \textit{Cubature} elements.
As defined in \cref{sec:discretization_cubature2D}, there are several requirements and optimization procedures in order to obtain the \textit{Cubature} elements. These elements are very import in our study because they permit to obtain diagonal mass matrix, and so they decrease considerably the time of computation. 
We describe for $p=1,2,3$ the basis functions of the \textit{Cubature} elements.

\subsection{\textit{Cubature} elements of degree 1}
%
The $\Tilde{\mathbb{P}}_1$ element contains $3$ degree of freedom. Their nodes are located at the vertices $v_1=(1,0,0)$, $v_2=(0,1,0)$ and $v_3=(0,0,1)$ of the triangle. 
\begin{itemize}
	\item At vertices of the triangle:
	\begin{equation*}
	\hspace*{-1cm}
    \phi_{v_i}(\lambda) = \lambda_i,
    \qquad \mbox{ for }i=1,2,3.
    \end{equation*}
\end{itemize}
Corresponding weights are $w_{v_i} = \frac{1}{3}$.

\subsection{\textit{Cubature} elements of degree 2}
%

The $\Tilde{\mathbb{P}}_2$ element contains $7$ degrees of freedom: three at the vertices $v_1$, $v_2$ and $v_3$ and three at the midpoint of the edges that we denote as  $e_{ij}=\frac{v_i+v_j}{2}$ for $(i,j)\in \lbrace (1,2),(2,3),(3,1) \rbrace$ and one at the  centroid point $G_\beta:=\frac{v_1+v_2+v_3}{3}$. Respectively, we have the following basis functions and weights:
\begin{itemize}
	\item At vertices of the triangle
	\begin{align*}
    &\phi_{v_i}(\lambda) = \lambda_i(2\lambda_i-1)+3\lambda_1 \lambda_2 \lambda_3, \text{ for }i  \in \llbracket 1,\dots, 3 \rrbracket,\\
        &w_{v}= \frac{1}{20};
    \end{align*}
	\item At edge midpoints 
	\begin{align*}
    &\phi_{e_{ij}}(\lambda) = 4 \lambda_i \lambda_j (1-3 \lambda_k), \text{ for all } i\neq j \neq k \neq i  \in \llbracket 1,\dots, 3 \rrbracket,\\
    &w_{e}=\frac{2}{15};
    \end{align*}    
	\item At the centroid
	\begin{align*}
    &\phi_{G_\beta}(\lambda) = 27 \lambda_1 \lambda_2 \lambda_3, \\ 
    &w_\beta = \frac{9}{20}.
    \end{align*}
\end{itemize}

\subsection{\textit{Cubature} elements of degree 3}
%
Following \cite{article_cubature_2001,article_cubature_2006} we derive the definitions of all the basis functions and points of \textit{Cubature} elements $\TP_3$. The notations are not uniform among different works, so we use the following one which can be used with all the different elements we have used in this work.

The space $\TP_3$ contains $12$ degrees of freedom: $3$ vertices $v_1$, $v_2$ and $v_3$, $6$ on edges: $e_{ij}^\alpha$ for $i,j \in \llbracket 1,\dots, 3 \rrbracket$ with $i\neq j$ defined by \begin{align*}
	&e_{ij}^\alpha = (\delta_{1i}\alpha +\delta_{1j}(1-\alpha),\delta_{2i}\alpha +\delta_{2j}(1-\alpha),\delta_{3i}\alpha +\delta_{3j}(1-\alpha))\\
	\text{ with }&\alpha = \frac{-15\sqrt{7}-21+\sqrt{168+174\sqrt{7}}}{2(-15\sqrt{7}-21)},
\end{align*} 
with $\delta_{ij}$ is the Kronecker delta and three internal points $G_i^\beta$ for $i \in \llbracket 1,\dots, 3 \rrbracket$, with 
\begin{align*}
G_i^\beta = \left(\beta \delta_{i1} + \frac{1-\beta}{2}(1-\delta_{i1}),\beta \delta_{i2} + \frac{1-\beta}{2}(1-\delta_{i2}), \beta \delta_{i3} + \frac{1-\beta}{2}(1-\delta_{i3})\right) \text{ with }\beta=\frac{1}{3}+\frac{2\sqrt{7}}{21},
\end{align*}
where $\alpha$ and $\beta$ are found through an optimization process \cite{article_cubature_2001,article_cubature_2006}. 
Let us start giving the definitions of the weights for the different types of points. 
We have that $w_v = \frac{1369+767\sqrt{7}}{120(859+395\sqrt{7})}$ is the weight for the vertices of the triangle, $w_\alpha= \frac{287+115\sqrt{7}}{40(173+49\sqrt{7})}$ is the weight on edges points, and $w_\beta $ the weight for barycentric points. \\
The weights corresponding to these types of points are $w_v = \frac{1369+767\sqrt{7}}{120(859+395\sqrt{7})}$, $w_\alpha = \frac{287+115\sqrt{7}}{40(173+49\sqrt{7})}$ and $w_\beta = \frac{21\sqrt{7}}{40(2\sqrt{7}+1)}$.
In order to simply the formulation of the basis functions, let us introduce some polynomials:
\begin{equation}
	p_i(\lambda) := \lambda_i \left( \sum_{l=1}^3 \lambda_l^2 - \frac{1-2\alpha + 2\alpha^2}{\alpha(1-\alpha)}\lambda_i(\lambda_j+\lambda_k)+ A_{i} \lambda_j \lambda_k \right), \qquad \text{with }j\neq i \neq k,
\end{equation}
with \begin{equation}
	\begin{split}
		A_{i} = &\left( w_v -\frac{1}{10}-\frac{1}{15} \left(1- \frac{1-2\alpha + 2 \alpha^2}{\alpha(1-\alpha)} \right) - \frac{1}{90}\frac{8}{\beta (1-\beta)^2 (3\beta-1)}\left(\sum_{l=1}^3 p_{i}(G_l)\right)      \right)  \frac{360}{6+\frac{8(1+\beta)  }{\beta (1-\beta) (3\beta-1)}} ;
	\end{split}
\end{equation} 
\begin{equation}
	p_{ij}(\lambda) :=  \frac{1}{\alpha(1-\alpha)(2\alpha-1)} \lambda_i \lambda_j (\alpha \lambda_i  - (1-\alpha)\lambda_j + (1-2\alpha) \lambda_k), \text{ with } i\neq j \neq k \neq i.
\end{equation}
We can then write the definition of the basis functions:
\begin{itemize}
	\item At vertices of the triangle
	\begin{align*}
    \phi_{v_i}(\lambda) = &p_{i}(\lambda)
     -\frac{8}{\beta (1-\beta)^2 (3\beta-1)}   \left( \sum_{l=1}^3 p_i(G_l) \left(\lambda_l - \frac{1-\beta  }{2}\right) \right)\prod_{l=1}^3\lambda_l, \text{ for }i \in \llbracket 1,\dots ,3 \rrbracket; 
    \end{align*}
	\item At the nodes on edges 
	\begin{align*}
    \phi_{e_{ij}^\alpha}(\lambda) = & p_{ij}(\lambda)- \frac{8}{\beta (1-\beta)^2 (3\beta-1)} \left( \sum_{l=1}^3 p_{ij}(G_l) \left(\lambda_l - \frac{1-\beta  }{2}\right) \right)\prod_{l=1}^3\lambda_l, \text{ for }i\neq j \in \llbracket 1,\dots ,3 \rrbracket;
    \end{align*}
	\item At the internal points
	\begin{equation*}
    \phi_{G_i^\beta}(\lambda) = \frac{8}{\beta(1-\beta)^2(3\beta-1)}  \left(\lambda_i-\frac{1-\beta}{2}\right)\prod_{l=1}^3\lambda_l, \text{ for }i \in \llbracket 1,\dots ,3 \rrbracket.
    \end{equation*}
\end{itemize}


\section{Time discretization coefficients}\label{sec:timeCoefficients}
In this appendix we introduce the time integration coefficients used in this work, to make the study fully reproducible. In  \cref{tab:Butcher} there are the RK coefficients, in \cref{tab:ButcherSSPRK} the SSPRK coefficients and in \cref{tab:DeCcoeff} the DeC coefficients.
\begingroup
\setlength{\tabcolsep}{6pt} 
\renewcommand{\arraystretch}{1.2} 
\begin{table}[h!]
	\centering
	\begin{tabular}{ c c  c }
		\begin{tabular}{| c |  c c  | }
			\hline 
			\multicolumn{3}{|c|}{\textit{RK2}} \\ \hline
			$\alpha$  & 1   &  \\ \hline
			$\beta$ & $\frac{1}{2}$ & $\frac{1}{2}$  \\ \hline
		\end{tabular} &
		\begin{tabular}{|c| c c c  | }
			\hline 
			\multicolumn{4}{|c|}{\textit{RK3}} \\ \hline
			$\alpha$  & $\frac{1}{2}$   &     & \\ 
			& -1 & 2  & \\ \hline
			$\beta$ & $\frac{1}{6}$ & $\frac{2}{3}$ & $\frac{1}{6}$ \\ \hline
		\end{tabular} &
		\begin{tabular}{| c | c  c c c | }
			\hline 
			\multicolumn{5}{|c|}{\textit{RK4}} \\ \hline
			$\alpha$ &  $\frac{1}{2}$   & & &   \\ 
			& 0 & $\frac{1}{2}$  & & \\ 
			& 0 & 0    & 1 &   \\ \hline 
			$\beta$ & $\frac{1}{6}$ & $\frac{1}{3}$ & $\frac{1}{3}$ & $\frac{1}{6}$ \\ \hline
		\end{tabular}
	\end{tabular}
	\caption{Butcher Tableau of RK methods}\label{tab:Butcher}
\end{table}
\begin{table}[h!]
	\begin{center}
		\begin{tabular}{| c  c c  |c c c  | }
			\hline 
			\multicolumn{6}{|c|}{\textit{SSPRK(3,2)} by \cite{shu-1988}} \\ \hline
			\multicolumn{3}{|c|}{$\gamma$} & \multicolumn{3}{|c|}{$\mu$} \\ \hline
			1 &  &        &    $\frac{1}{2}$  &   &   \\ 
			0 & 1 &   &   0 & $\frac{1}{2}$  & \\ 
			$\frac{1}{3}$ & 0 & $\frac{2}{3}$ &   0 & 0 & $\frac{1}{3}$ \\ \hline
			\multicolumn{6}{|c|}{CFL = 2.} \\ \hline
		\end{tabular} \qquad
		\begin{tabular}{| c  c c c  |  c c c c | }
			\hline 
			\multicolumn{8}{|c|}{\textit{SSPRK(4,3)} by \cite[Page 189]{Ruuth-2006}} \\ \hline
			\multicolumn{4}{|c|}{$\gamma$} & \multicolumn{4}{|c|}{$\mu$} \\ \hline
			1 &   &       &      & $\frac{1}{2}$ &  &     & \\ 
			0 & 1 &       &      & 0 & $\frac{1}{2}$ & &  \\ 
			$\frac{2}{3}$ & 0 & $\frac{1}{3}$ &      & 0 & 0 & $\frac{1}{6}$     & \\ 
			0   & 0 &  0  & 1    & 0 & 0 & 0    & $\frac{1}{2}$  \\ \hline
			\multicolumn{8}{|c|}{CFL = 2.} \\ \hline
		\end{tabular}\\
		\begin{tabular}{| c c c c c |}
			\hline 
			\multicolumn{5}{|c|}{\textit{SSPRK(5,4)} by \cite[Table 3]{Ruuth-2006} } \\ \hline
			\multicolumn{5}{|c|}{$\gamma$}   \\ \hline
			1 &   &  &     &                                                \\ 
			0.444370493651235 & 0.555629506348765 &   &    &                   \\ 
			0.620101851488403 & 0                 & 0.379898148511597 & &       \\ 
			0.178079954393132 & 0 &  0  & 0.821920045606868 &                 \\ 
			0 & 0 & 0.517231671970585 & 0.096059710526147 & 0.386708617503269  \\ \hline
			\multicolumn{5}{|c|}{$\mu$} \\ \hline
			0.391752226571890 &  & &   & \\
			0 & 0.368410593050371 & & &   \\
			0 & 0 & 0.251891774271694  &    &\\
			0 & 0 & 0    & 0.544974750228521  &\\
			0 & 0 & 0 & 0.063692468666290 & 0.226007483236906 \\ \hline
			\multicolumn{5}{|c|}{CFL = 1.50818004918983} \\ \hline
		\end{tabular}
		\caption{Butcher Tableau of SSPRK methods}\label{tab:ButcherSSPRK}
	\end{center}
\end{table}
\begin{table}[h!]
	\begin{center}
		\begin{tabular}{ |c|c|c c |}\hline
			\multicolumn{4}{|c|}{Order 2}\\ \hline
			m & $\beta^m$ & \multicolumn{2}{|c|}{$\rho^{m}_z$}\\ \hline
			1 & 1       & $\frac{1}{2}$ & $\frac{1}{2}$ \\ \hline
		\end{tabular}
		\quad
		\begin{tabular}{ |c|c|c c c |}\hline
			\multicolumn{5}{|c|}{Order 3}\\ \hline
			m & $\beta^m$ & \multicolumn{3}{|c|}{$\rho^{m}_z$}\\ \hline 
			1 & $\frac{1}{2}$      & $\frac{5}{24}$ & $\frac{1}{3}$  & $-\frac{1}{24}$ \\ 
			2 & 1         & $\frac{1}{6}$ & $\frac{2}{3}$  & $\frac{1}{3}$\\ \hline
		\end{tabular}
		\quad
		\begin{tabular}{ |c|c|c c c c|}\hline
			\multicolumn{6}{|c|}{Order 4}\\ \hline
			m & $\beta^m$ & \multicolumn{4}{|c|}{$\rho^{m}_z$}\\ \hline 
			1 & $\frac{1}{3}$      & $\frac{1}{8}$ & $\frac{19}{72}$  & $-\frac{5}{72}$  & $\frac{1}{72}$ \\ 
			2 & $\frac{2}{3}$       & $\frac{1}{9}$ & $\frac{4}{9}$  & $\frac{1}{9}$  & 0\\ 
			3 & 1         & $\frac{1}{8}$ & $\frac{3}{8}$  & $\frac{3}{8}$ & $\frac{1}{8}$\\ \hline
		\end{tabular}
		\caption{DeC coefficients for equispaced subtimesteps. }\label{tab:DeCcoeff}
	\end{center}
\end{table}
\endgroup
\section{Fourier analysis} \label{sec:appFourier}
In this section we collect all the plots and results that are essential to show the results of this work, but for structural reasons were not put in the main text.

\subsection{Mesh types and degrees of freedom}\label{sec:app_mesh4fourier}
\begin{figure}
	\centering
	\subfigure[\textit{X} mesh, \textit{Basic} elements \label{fig_Xmesh_P3}]{	\includegraphics[width=0.35\linewidth]{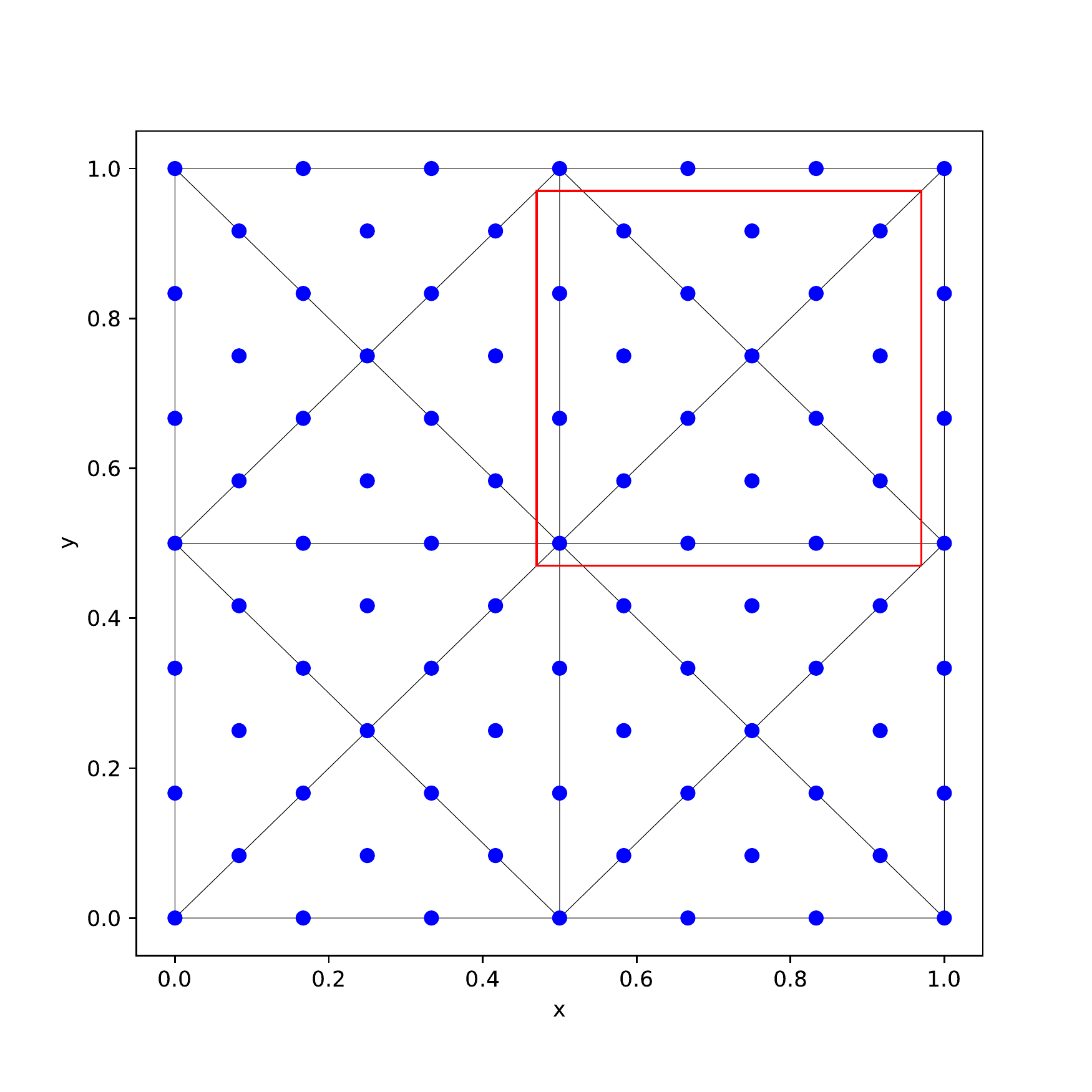}}
	\subfigure[\textit{X} mesh, \textit{Cubature} elements]{
		\includegraphics[width=0.35\linewidth]{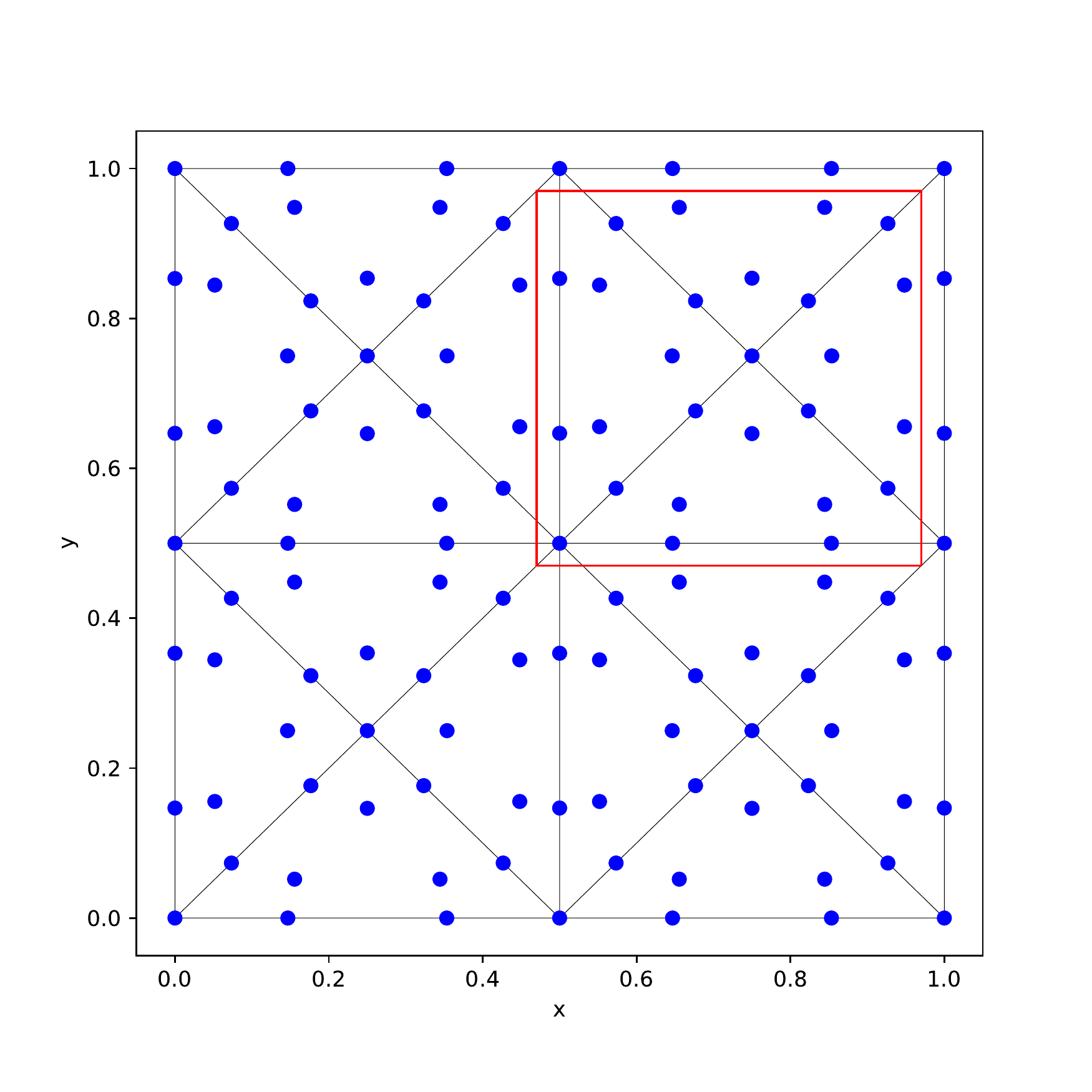}
	}
\subfigure[\textit{T} mesh, \textit{Cubature} elements]{
	\includegraphics[width=0.35\linewidth]{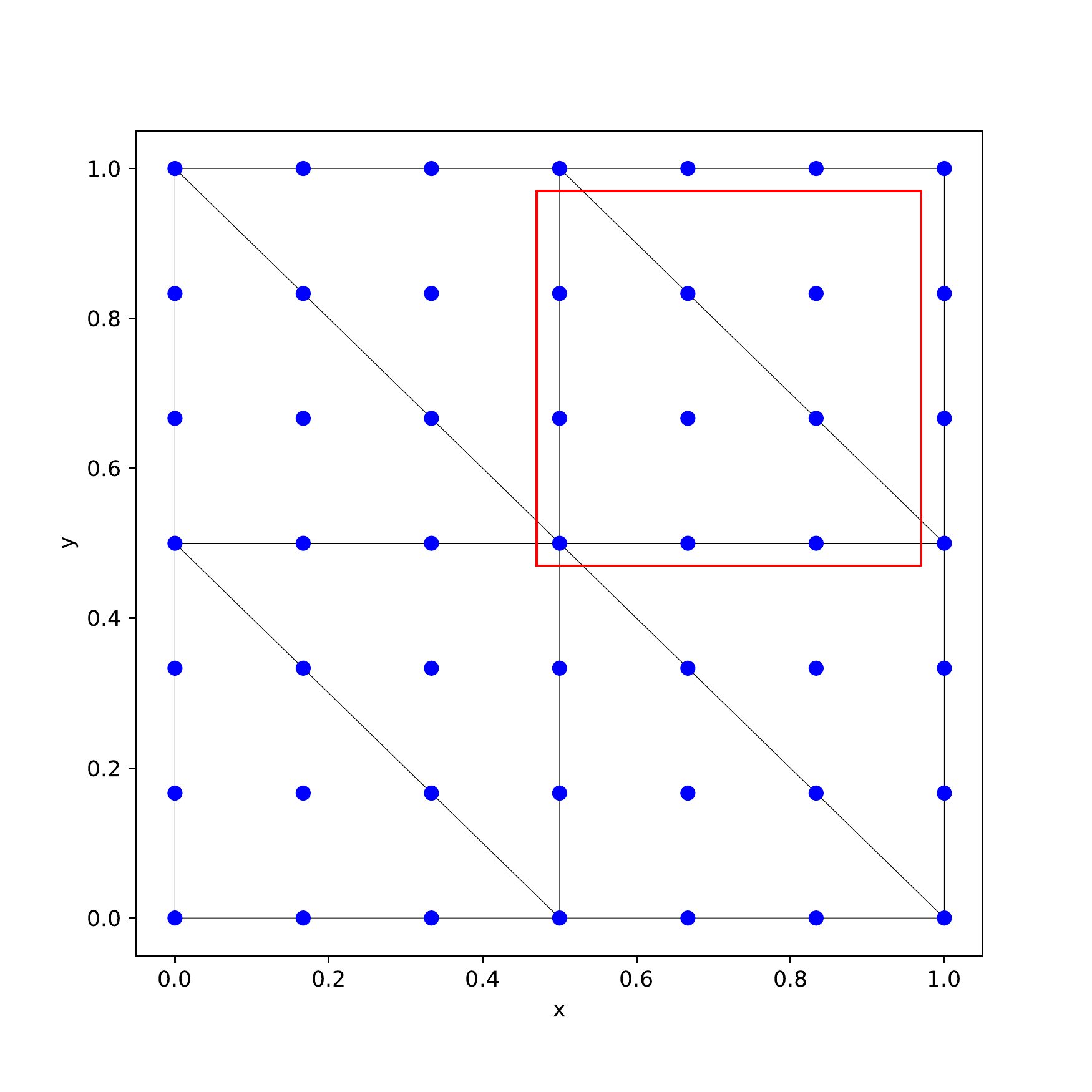}}
\subfigure[\textit{T} mesh, \textit{Cubature} elements]{
	\includegraphics[width=0.35\linewidth]{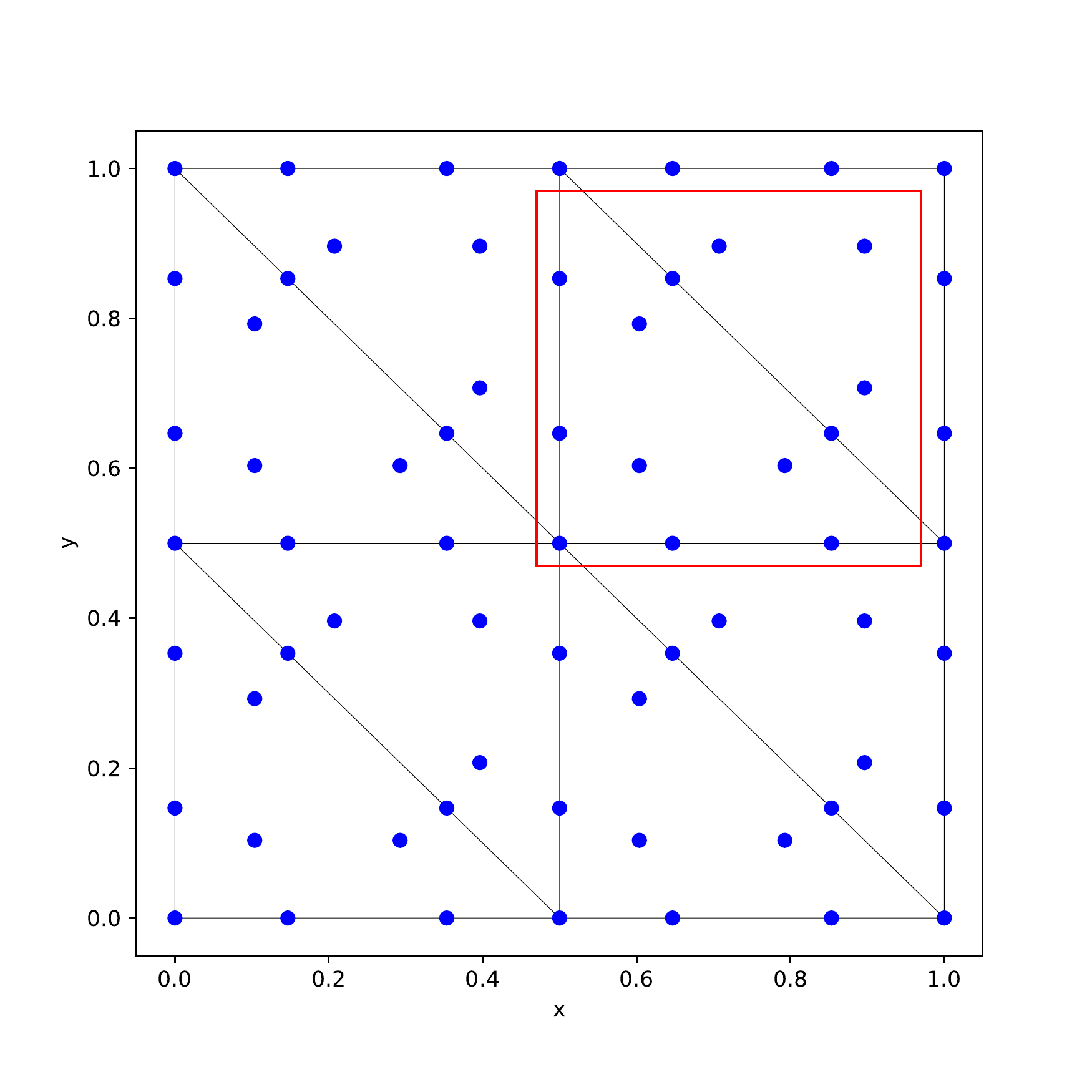}}
	\caption{Degrees of freedom and periodic unit for different mesh patterns and elements of degree 3 \label{fig:meshDofsP3}}
\end{figure}

We represent in Figure~\ref{fig:meshDofsP3} the mesh configurations used in the Fourier analysis and the degrees of freedom of the elements of degree 3. The red square represents the periodic elementary unit that contains the degrees of freedom of interest for the Fourier analysis.

\subsection{Fourier analysis results - Optimal Parameters} \label{sec:app_fourier_param}
In this section, we put the optimal values of the stability analysis of Section~\ref{sec:fourier_globalresults} after the modification proposed in Section~\ref{sec:fourier_rocket}. In Table~\ref{tab:new_param_meshX_LinearAdvection-2D-RES} we show the parameters for the \textit{X} mesh and in Table~\ref{tab:new_param_meshT_LinearAdvection-2D-RES} we show the parameters for the \textit{T} mesh.
\begin{table}[H] 
\small  
 \begin{center} 
		\begin{tabular}{| c | c || c | c | c | }  
	     \hline 
	     \multicolumn{2}{|c||}{Element $\&$ }  & \multicolumn{3}{|c|}{SUPG}  \\ \hline 
	     \multicolumn{2}{|c||}{ Time scheme }  & $\mathbb{P}_1$ & $\mathbb{P}_2$ & $\mathbb{P}_3$  \\ \hline \hline 
Basic              &  SSPRK  &  0.739 (0.127)   &  0.298 (0.058)   &  0.22 (0.026)  \\ 
         \hline 
       \multirow{2}{*}{ Cub.}               &  SSPRK  &  1.062 (0.28)   &  0.1 (0.1)$^*$   &  0.18 (0.04)$^*$  \\ 
               &  DeC  &  0.616 (0.28)   &  0.1 (0.04)$^*$   &  0.144 (0.04)  \\ 
         \hline 
Bern.              &  DeC  &  0.739 (0.298)   &  0.2 (0.2)$^*$   &  0.2 (0.153)$^*$  \\ 
         \hline 
        \end{tabular} 
    \end{center} 
 \begin{center} 
		\begin{tabular}{| c | c || c | c | c || c | c | c | }  
	     \hline 
	     \multicolumn{2}{|c||}{Element $\&$ }  & \multicolumn{3}{|c||}{OSS}  & \multicolumn{3}{|c|}{CIP}  \\ \hline 
	     \multicolumn{2}{|c||}{ Time scheme }  & $\mathbb{P}_1$ & $\mathbb{P}_2$ & $\mathbb{P}_3$   & $\mathbb{P}_1$ & $\mathbb{P}_2$ & $\mathbb{P}_3$  \\ \hline \hline 
Basic              &  SSPRK  &  0.403 (0.127)   &  0.298 (0.026)   &  0.22 (0.026)   &  0.403 (0.012)   &  0.298 (1.73e-03)   &  0.1 (1.00e-03)$^*$   \\ 
         \hline 
       \multirow{2}{*}{ Cub.}                &  SSPRK  &  0.58 (0.336)   &  0.379 (0.03)   &  0.248 (0.018)   &  0.58 (0.048)   &  0.06 (0.01)$^*$   &   /   \\ 
               &  DeC  &  0.379 (0.207)   &  0.248 (0.03)   &  0.162 (0.018)   &  0.379 (0.026)   &  0.06 (0.01)$^*$   &   /   \\ 
         \hline 
Bern.              &  DeC  &  0.173 (0.58)   &  0.036 (0.298)   &  0.015 (0.078)$^*$   &  0.173 (0.153)   &  0.012 (0.021)   &  0.002 (8.00e-03)$^*$   \\ 
         \hline 
        \end{tabular} 
    \end{center} 
     \caption{\textit{X} mesh: Optimized CFL and penalty coefficient $\delta $ in parenthesis.
     The symbol "/" means that the fourier analysis for the scheme results always in instability.
     The values denoted by $^*$ are not the optimal one, but they lay in a safer region, see Section~\ref{sec:fourier_rocket}.} \label{tab:new_param_meshX_LinearAdvection-2D-RES}
\end{table}%
\begin{table}[H] 
\small  
 \begin{center} 
		\begin{tabular}{| c | c || c | c | c | }  
	     \hline 
	     \multicolumn{2}{|c||}{Element $\&$ }  & \multicolumn{3}{|c|}{SUPG}  \\ \hline 
	     \multicolumn{2}{|c||}{ Time scheme }  & $\mathbb{P}_1$ & $\mathbb{P}_2$ & $\mathbb{P}_3$  \\ \hline \hline 
Basic              &  SSPRK  &  0.739 (0.127)   &  0.403 (0.026)   &  0.298 (0.012)  \\ 
         \hline 
       \multirow{2}{*}{ Cub.}              &  SSPRK  &  1.062 (0.28)   &  0.234 (0.078)   &  0.055 (0.153)  \\ 
               &  DeC  &  1.062 (0.127)   &  0.144 (0.078)   &  0.034 (0.153)  \\ 
         \hline 
Bern.              &  DeC  &  0.739 (0.298)   &  0.739 (0.153)   &  0.455 (0.153)  \\ 
         \hline
        \end{tabular} 
    \end{center} 
 \begin{center} 
		\begin{tabular}{| c | c || c | c | c || c | c | c | }  
	     \hline 
	     \multicolumn{2}{|c||}{Element $\&$ }  & \multicolumn{3}{|c||}{OSS}  & \multicolumn{3}{|c|}{CIP}  \\ \hline 
	     \multicolumn{2}{|c||}{ Time scheme }  & $\mathbb{P}_1$ & $\mathbb{P}_2$ & $\mathbb{P}_3$   & $\mathbb{P}_1$ & $\mathbb{P}_2$ & $\mathbb{P}_3$  \\ \hline \hline 
Basic              &  SSPRK  &  0.546 (0.127)   &  0.403 (0.058)   &  0.298 (0.012)   &  0.546 (0.026)   &  0.298 (7.39e-05)   &  0.298 (3.36e-05)   \\ 
         \hline 
       \multirow{2}{*}{ Cub.}              &  SSPRK  &  0.886 (0.336)   &  0.379 (0.048)   &   /   &  0.886 (0.048)   &  0.106 (7.85e-03)   &   /   \\ 
               &  DeC  &  0.58 (0.207)   &  0.379 (0.03)   &   /   &  0.58 (0.026)   &  0.045 (7.85e-03)   &   /   \\ 
         \hline 
Bern.              &  DeC  &  0.28 (0.58)   &  0.025 (0.153)   &  0.074 (0.078)   &  0.455 (0.078)   &  0.025 (5.46e-03)   &  0.017 (0.04)   \\ 
         \hline 
        \end{tabular} 
    \end{center} 
     \caption{\textit{T} mesh: Optimized CFL and penalty coefficient $\delta $ in parenthesis.
     The symbol "/" means that the fourier analysis for the scheme results always in instability.} \label{tab:new_param_meshT_LinearAdvection-2D-RES}
\end{table}%

\subsection{Fourier analysis results - stability area} \label{sec:fourier_results_appendix}
\begin{figure}
	\centering
	\subfigure[\textit{X} mesh \label{fig:App_fourier_global_basic_SSPRK_OSSX}]{
	\includegraphics[width=\linewidth]{image_dispersion_2D/meshX/cfl_vs_tau_angleLoop/error_disp_cfl_tau_2D_lagrange_SSPRK_OSS.pdf}}
	\subfigure[\textit{T} mesh\label{fig:App_fourier_global_basic_SSPRK_OSST}]{
		\includegraphics[width=\linewidth]{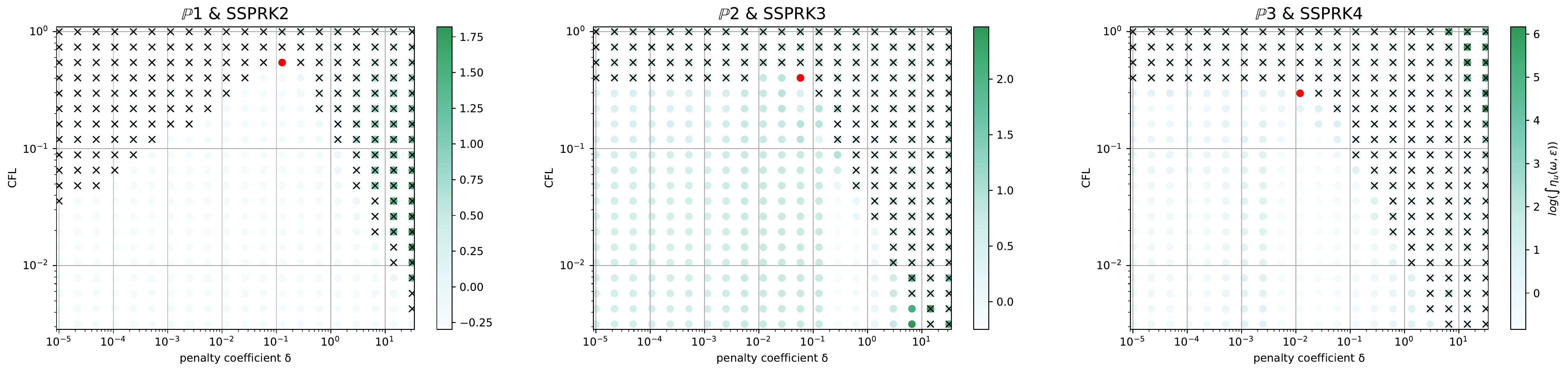}}
	\caption{$\log(\eta_u)$ values (blue scale) and stable area (unstable with black crosses), on $(\CFL,\delta)$  plane. The red dot denotes the optimal value. From left to right $\P_1$, $\P_2$, $\P_3$ \textit{Basic} elements with SSPRK scheme and OSS stabilization.}  \label{fig:App_fourier_global_basic_SSPRK_OSS}
\end{figure}

Finally, we present a comparison of stability area between the \textit{T} and the \textit{X} mesh. This comparison if perform as before, for all wave angles $\theta$. We choose as example the comparison using \textit{Basic} element, SSPRK time integration method and the OSS stabilization technique in \cref{fig:App_fourier_global_basic_SSPRK_OSS}.
The interested reader can access to results for all methods online \cite{TorloMichel2021git}.

\newpage
\medskip
\bibliographystyle{abbrv} 
\bibliography{bibliography}

\end{document}